\documentclass[11pt,letterpaper,reqno]{amsart}
\usepackage[margin={1.25in,1.25in}]{geometry}
\usepackage[T1]{fontenc}
\usepackage{amsmath,amssymb,mathrsfs,mathtools}
\usepackage{amsthm}
\usepackage{graphicx}
\usepackage{aliascnt}
\usepackage{comment}
\usepackage{tikz-cd}
\usepackage[dvipsnames]{xcolor}
\definecolor{linkBlue}{rgb}{.05,.25,.6}
\usepackage[colorlinks=true,linkcolor=linkBlue,citecolor=linkBlue,urlcolor=linkBlue]{hyperref}
\usepackage[nameinlink,noabbrev]{cleveref}
\theoremstyle{definition}
\newtheorem{definition}{Definition}[section]

\newaliascnt{example}{definition}

\aliascntresetthe{example}
\theoremstyle{plain}
\newaliascnt{theorem}{definition}
\newtheorem{theorem}[theorem]{Theorem}
\aliascntresetthe{theorem}
\newcommand{\namedtheoremname}{}
\newtheorem*{namedtheoreminner}{\namedtheoremname}
\newenvironment{namedtheorem}[1]{
    \renewcommand{\namedtheoremname}{#1}
    \begin{namedtheoreminner}
}{
    \end{namedtheoreminner}
}
\newaliascnt{proposition}{definition}
\newtheorem{proposition}[proposition]{Proposition}
\aliascntresetthe{proposition}
\newaliascnt{lemma}{definition}
\newtheorem{lemma}[lemma]{Lemma}
\aliascntresetthe{lemma}
\newaliascnt{corollary}{definition}
\newtheorem{corollary}[corollary]{Corollary}
\aliascntresetthe{corollary}
\newaliascnt{conjecture}{definition}
\newtheorem{conjecture}[conjecture]{Conjecture}
\aliascntresetthe{conjecture}
\newtheorem{claim}{Claim}[definition]
\theoremstyle{remark}
\newaliascnt{remark}{definition}
\newtheorem{remark}[remark]{Remark}
\aliascntresetthe{remark}

\crefname{claim}{claim}{claims}
\Crefname{claim}{Claim}{Claims}
\let\tilde\widetilde
\let\bar\overline
\let\hat\widehat
\newcommand{\de}{\text{d}}
\newcommand{\from}{\,:\,}
\newcommand{\into}{\hookrightarrow}
\newcommand{\acts}{\mathbin{{\text{\raisebox{.75em}{\rotatebox{-90}{$\circlearrowright$}}}}}}
\newcommand{\sm}{\smallsetminus}
\newcommand{\R}{\mathbb R}
\newcommand{\C}{\mathbb C}
\renewcommand{\P}{\mathbb P}
\newcommand{\Z}{\mathbb Z}

\newcommand{\Q}{\mathbb Q}
\newcommand{\CP}{\mathbb {CP}}
\newcommand{\cM}{\mathcal M}
\newcommand{\cX}{\mathcal X}
\newcommand{\cT}{\mathcal T}
\newcommand{\fb}{\mathfrak b}
\newcommand{\I}{\mathrm i}
\newcommand{\eps}{\varepsilon}
\renewcommand{\Re}{\operatorname{Re}}
\renewcommand{\Im}{\operatorname{Im}}
\DeclareMathOperator{\rk}{rk}
\DeclareMathOperator{\im}{im}
\DeclareMathOperator{\gr}{gr}
\DeclareMathOperator{\End}{End}
\DeclareMathOperator{\SL}{SL}
\DeclareMathOperator{\SU}{SU}
\DeclareMathOperator{\U}{U}
\DeclareMathOperator{\res}{res}
\DeclareMathOperator{\tr}{tr}

\DeclareMathOperator{\Id}{Id}
\DeclareMathOperator{\delbar}{\bar\partial}
\newcommand{\pa}{{\partial}}

\title[Periods in Families of Parabolic Higgs Bundle Moduli Spaces]{Hyperk\"ahler Periods in Families of Parabolic Higgs Bundle Moduli Spaces}

\author{Panagiotis Dimakis}
\email{pdimakis@umd.edu}

\author{Laura Fredrickson}
\email{lfredric@uoregon.edu}

\author{Arya Yae}
\email{aryayae@outlook.com}

\date{}
\begin{document}
\newif\ifFiguresThreeD
\FiguresThreeDtrue
\begin{abstract}
We consider the period map on families of parabolic Higgs bundle moduli spaces with varying parabolic weights $\alpha$ and Higgs field eigenvalues $\sigma$.  We assemble the moduli spaces into a smooth family $\pi\from\mathscr M\to\Lambda_0$ where $\Lambda_0$ encodes the admissible pairs $(\alpha,\sigma)$, and define parallel transport for homology along paths in $\Lambda_0$.  By generalizing a classical theorem of Duistermaat and Heckman, we prove that the period map is affine linear in $(\alpha,\sigma)$, with slope governed by the curvature class of a principal torus bundle over $\mathscr M$.  For the four-dimensional moduli spaces of types $E_6,E_7,E_8$, we compute these classes to get explicit formulas for the period maps.
\end{abstract}
\maketitle
\tableofcontents

\section{Introduction}\label{sec: Introduction}
We consider parabolic $\SL(r_\star,\C)$-Higgs bundles on a compact Riemann surface with fixed flag dimensions. The goal of this paper is to understand how the hyperk\"ahler periods $(\int_\eta \omega_1, \int_\eta \omega_2, \int_\eta \omega_3)$ of these Hitchin moduli spaces change as we vary the parabolic weights and the complex residues. We prove that, after canonically identifying homology classes, the three symplectic classes vary affine linearly with slopes given by a common characteristic class. We then compute these slopes for the four-dimensional moduli spaces of exceptional types $E_6,E_7,E_8$.

Let $C$ be a compact connected Riemann surface, $D=\{p_1,\ldots,p_n\}$ a divisor of distinct marked points on $C$, and $E$ a rank-$r_\star$ bundle. Fix an identification $\det E\cong\mathcal O_C$. We consider parabolic $\SL(r_\star,\C)$-Higgs bundles with flag dimensions
\[
    0<r_i^{(0)}<\cdots<r_i^{(\ell_i)}=r_\star,
    \qquad i\in\{1,\dots,n\}.
\]
We denote the flags and the parabolic weights associated to them by
\[
    0\subset F_i^{(0)}\subset\cdots\subset F_i^{(\ell_i)}=E_{p_i},
    \qquad
    \alpha_i^{(0)}>\cdots>\alpha_i^{(\ell_i)},
\]
where $\dim F_i^{(k)}=r_i^{(k)}$ and $\alpha_i^{(0)}-\alpha_i^{(\ell_i)}<1$. We allow partial flags, and denote the flag multiplicities $ m_i^{(k)}=r_i^{(k)}-r_i^{(k-1)}$ for $k=0,\ldots,\ell_i$, where $r_i^{(-1)}=0$. The residue of the holomorphic Higgs field $\varphi$ at $p_i$ is required to preserve the flag and act on the associated graded space $\gr F_i^{(k)}=F_i^{(k)}/F_i^{(k-1)}$ as scalar multiplication by $\sigma_i^{(k)}\in\C$. In other words, the complex residue $\res_{p_i}\varphi$ has prescribed Levi part $\bigoplus_k\sigma_i^{(k)}\Id$. The parabolic weights and complex residues satisfy the trace-free conditions for each $i\in\{1,\dots,n\}$
\[
    \sum_{k=0}^{\ell_i}m_i^{(k)}\alpha_i^{(k)}=0,
    \qquad
    \sum_{k=0}^{\ell_i}m_i^{(k)}\sigma_i^{(k)}=0.
\]

We write $\cM(\alpha,\sigma)$ for the moduli space of parabolic $\SL(r_\star,\C)$-Higgs bundles with these weights and residues.  We emphasize that only the Levi part of $\res_{p_i}\varphi$ is prescribed, while the upper triangular blocks are allowed to vary.  On its smooth stable locus, this moduli space carries a natural hyperk\"ahler structure \cite{BB04}, with K\"ahler forms $\omega_1,\omega_2,\omega_3$. Given an oriented cycle $\eta$ representing a class in $H_2(\cM(\alpha,\sigma),\Z)$, define its period triple by
\begin{equation}\label{eq: period triple}
    P_\eta=
    \left(\int_\eta\omega_1,\int_\eta\omega_2,\int_\eta\omega_3\right).
\end{equation}
To compare these triples as the weights and residues vary, we first need to identify homology classes in different moduli spaces. 

Consider the parameter space $\Lambda$ whose elements $\nu$ encode pairs $(\alpha,\sigma)$ satisfying the trace-free conditions above. We work over the generic parameter locus $\Lambda_0\subset\Lambda$ of \Cref{def: generic locus} so as to guarantee that $\cM_\nu=\cM(\alpha,\sigma)$ is a smooth hyperk\"ahler manifold. This $\Lambda_0$ is obtained by excluding the simultaneous weight and residue relations associated with a possible proper invariant subbundle. The first condition cuts out a real co-dimension one affine subspace inside $\Lambda$ and the second a complex co-dimension one subspace. Under the co-primality assumption \[\gcd\bigl(r_\star,\{m_i^{(k)}\}_{i,k}\bigr)=1,\] these excluded loci form a finite union of real co-dimension three affine subspaces inside $\Lambda$. In particular, this implies that $\Lambda_0$ is simply connected. 

In \Cref{sec: Universal Moduli Space} and \Cref{sec: The Analytic Family}, we assemble the fixed-parameter moduli spaces $\cM(\alpha,\sigma)$ into a smooth family $\mathscr M$ with a submersion $\pi:\mathscr M\to \Lambda_0$ whose fibers are the spaces $\cM(\alpha,\sigma)$.  The required canonical identification of homology classes follows from the stronger result that $\pi$ is in fact a fiber bundle. We prove this by showing that the family is locally topologically trivial over $\Lambda_0$, which implies that its integral homology groups form a local system compatible with the inclusions of the fibers into the total space. Thus a path $\gamma:[0,1]\to\Lambda_0$ gives a parallel transport identification
\[
    \mathrm{PT}_\gamma:
    H_2(\cM_{\gamma(0)},\Z)\longrightarrow H_2(\cM_{\gamma(1)},\Z)
\]
and dual parallel transport identification
\[
    \mathrm{PT}_\gamma^*:
    H^2(\cM_{\gamma(1)},\R)\longrightarrow H^2(\cM_{\gamma(0)},\R).
\]
Compatibility means that parallel transport commutes with the inclusions $i_\nu\from \cM_\nu\into\mathscr M$. Since $\Lambda_0$ is simply connected, these identifications depend only on the endpoints. The identifications of homology classes are canonical, even though the local trivializations used to construct them depend on choices.

To describe the variation of the three K\"ahler classes $[\omega_{1,\nu}]$, $[\omega_{2,\nu}]$, $[\omega_{3,\nu}]$ for $\nu\in\Lambda_0$, let $B$ be the boundary torus constructed in \Cref{sec: The Gauge Group}, with Lie algebra $\mathfrak b$.  This torus is a quotient of the gauge group, and has dimension $d=\sum_i\ell_i$.  On the other hand, the real dimension of $\Lambda_0$ is $3d$.  Using the duality pairing and the normalization of \eqref{eq: nu normalized definition}, we write the parameters as $\nu=(\nu_1,\nu_2,\nu_3)\in(\mathfrak b^*)^{\oplus3}$, where
\[
    \nu_1=2\pi\I\alpha,\quad
    \nu_2=4\pi\I\Re\sigma,\quad
    \nu_3=4\pi\I\Im\sigma.
\]
Here the weights and residues are regarded as block scalar matrices.

The simplest case we could hope for is that there is a common class
$c_\nu\in H^2(\cM_\nu,\R)\otimes\mathfrak b$, identified throughout $\Lambda_0$ via parallel transport, such that the Gauss--Manin directional derivative along a variation $\dot\nu$ of $\nu\in\Lambda_0$ is
\[
    D_{\dot\nu}[\omega_{s,\nu}]
    =-\left\langle c_\nu,\dot\nu\right\rangle
\]
for all $t$.  For a smooth path $\gamma(t):[0,1]\to \Lambda_0$, we can integrate to obtain
\begin{equation}\label{eq: intro affine periods}
    \mathrm{PT}_\gamma^*[\omega_{s,\gamma(1)}]
    =
    [\omega_{s,\gamma(0)}]
    -\langle c_{\gamma(0)},\gamma_s(1)-\gamma_s(0)\rangle,
    \qquad s=1,2,3.
\end{equation}
\Cref{thm: general period variation} proves precisely this formula. In particular, the same class governs all three variations.

There is a natural geometric candidate for the class $c$. Suppose that we can construct a principal $B$-bundle $q:\mathscr E\to\mathscr M$. Let $\theta$ be a connection and let $F_\theta$ denote its curvature, viewed as a $\mathfrak b$-valued two-form on $\mathscr M$. If $i_\nu:\cM_\nu\hookrightarrow\mathscr M$ is the inclusion, then the restricted curvature class
\[
    c_\nu:=i_\nu^*[F_\theta]
\]
is invariant under parallel transport. This gives a candidate class, but to identify it with the Gauss--Manin derivative of the symplectic classes we need additional geometric structure.

Specifically, suppose $\mathscr E$ carries three closed, $B$-invariant two-forms $\varpi_s$ which satisfy the Hamiltonian identity
\[
    \iota_{\zeta_b}\varpi_s=d\langle\widetilde\pi_s,b\rangle,
    \qquad b\in\mathfrak b,\quad s=1,2,3,
\]
where $\widetilde\pi=(\widetilde\pi_1,\widetilde\pi_2,\widetilde\pi_3)$ is the lifted parameter map $\mathscr E\to\Lambda_0$ and $\zeta_b$ is the vector field generated by $b$.  Furthermore, assume $\varpi_s$ descends to $\omega_s$ under the quotient by $B$.  We can then show
\[
    \varpi_s+d\langle\theta,\widetilde\pi_s\rangle
\]
descends to a closed two-form $\eta_s$ on $\mathscr M$. Restriction to the fiber over $\nu$ and calculating the exterior derivative yields
\[
    [\omega_{s,\nu}]
    =i_\nu^*[\eta_s]-\langle c_\nu,\nu_s\rangle.
\]
Both $i_\nu^*[\eta_s]$ and $c_\nu$ are constant under parallel transport, so this identity gives \eqref{eq: intro affine periods}.  Thus the Hamiltonian identities realize the curvature class $c$ as the Gauss--Manin derivative, while compatibility of homology transport with the inclusions $i_\nu$ makes that slope constant.

\bigskip

The process described above is almost identical to the mechanism behind the theorem of Duistermaat and Heckman \cite[Theorem~1.1]{DH82}. In its classical torus formulation, a Hamiltonian action of a compact torus $G$ on a finite-dimensional symplectic manifold $M$, with moment map
\[
    \mu:M\longrightarrow\mathfrak g^*,
\]
gives affine variation of the reduced symplectic classes on a suitable region of regular moment-map values. In our setting, the analogous hyperk\"ahler quotient involves an infinite-dimensional configuration space and an infinite-dimensional gauge group. Moreover, the reduced moduli spaces $\cM(\alpha,\sigma)$ are noncompact. 

\Cref{sec: a generalized DH theorem} develops a general framework for dealing with situations of this type. The first observation is that, under suitable hypotheses, partial reduction by a subgroup $H$ leaves a finite-dimensional space with an action of the residual group $B=G/H$. The second observation is that the Duistermaat--Heckman proof applies once a geometric and a topological hypothesis are satisfied. The geometric hypothesis essentially is the existence of the principal bundle, closed forms, and Hamiltonian identities described above, while the topological hypothesis is the existence of a homology local system compatible with the inclusions into the total space.
\begin{remark}
As in the proof in \cite{DH82}, the existence of a complete Ehresmann connection on the total space guarantees that the topological hypothesis is satisfied. We illustrate this through an example, treating the multi Taub-NUT family in \Cref{sec: multi Taub NUT}.  Proving the existence of a complete Ehresmann connection is highly non-trivial both generally and in particular for the case treated in this paper.
\end{remark}

In \Cref{sec: Universal Moduli Space}, we produce the structures required by the geometric hypothesis of \Cref{thm: DH after partial reduction} in the case of parabolic Higgs bundles. First we define the extended configuration space $\mathscr C$ allowing varying boundary data, together with the appropriate gauge group $G$. Then we extend the usual $L^2$ K\"ahler forms to a triple of closed $2$-forms $\tilde\omega_1,\tilde\omega_2,\tilde\omega_3$ on $\mathscr C$, and produce a map
\[
    \mu=(\mu_1,\mu_2,\mu_3)\from\mathscr C\to(\mathfrak g^*)^{\oplus3}
\]
satisfying a Hamiltonian identity with respect to the $\tilde\omega_s$ and the $G$-action. The moment map is obtained by recasting Hitchin's equations using distributional derivatives and regarding them as bounded linear functionals on the Lie algebra $\mathfrak g$. Taking level sets, the interior part $\mu^\mathrm{int}$ gives the usual Hitchin equations, and the boundary part $\mu^\mathrm{bdr}$ gives the weight and residue prescriptions.  We identify $\Lambda\subset(\mathfrak g^*)^{\oplus3}$ using point-mass distributions supported on $D$.

The gauge group $G$ contains a normal subgroup $H$ such that $G/H$ is a finite-dimensional \emph{boundary torus} $B$. The forms $\tilde\omega_s$ descend to closed forms $\varpi_s$ on the partial quotient $\mathscr E=(\mu^\mathrm{int})^{-1}(0)/H$. The latter retains a residual action of $B$, and the boundary moment map descends to a moment map $\tilde\pi\from\mathscr E\to\Lambda\subset(\mathfrak g^*)^{\oplus3}$ for the $B$-action which reads off the parameters.  Taking one more quotient gives the moduli family $\pi\from\mathscr M\to\Lambda$.

The existence of these closed 2-forms $\varpi_s$ and the Hamiltonian identity for $\tilde\pi$ is precisely the geometric hypothesis of \Cref{thm: DH after partial reduction}.  The difficulty in defining them is that the usual $L^2$ forms on $\mathscr C$ have divergent integrals when evaluated on deformations with varying boundary data.  To account for this, we define the truncated forms $\omega_s^\eps$ by integrating over the complement of small balls $B_\eps(p_i)$ in $C$.  Subtracting logarithmically divergent boundary terms $\log(1/\eps)b_s(u,v)$ obtained by pairing the parameter variations, we get a well-defined form
\[
    \tilde\omega_s=\lim_{\eps\to0}\left(\omega_s^\eps+\log(\eps)b_s\right).
\]
The fact that they descend to the forms $\varpi_s$ on $\mathscr E$ follows from the Hamiltonian identity for $\mu^\mathrm{int}$ with respect to the $H$-action.

\begin{remark}
    Although we only consider regular singularities in this paper, the distributional formulation suggests an extension to irregular singularities.  In a local coordinate centered at $p_i$, the eigenvalues of the higher order singular terms $A_nz^{-n}\,\de z$ are prescribed by allowing $\Lambda$ to contain Levi-central matrices times distributions of the form $\partial_z^{n-1}\delta_{p_i}$, after appropriately shrinking the gauge Lie algebra near the divisor.
\end{remark}

\Cref{sec: The Analytic Family} supplies the smoothness assumptions of the affine variation theorem.  For fixed Sobolev regularity $k\geq2$ and weight $\delta>0$, we consider the extended deformation complex for the action of $G$ on $\mathscr M$, whose cohomology provides a tangent space model for $\mathscr M^{\mathrm{irr}}$.  We use Biquard--Boalch's splitting of the fixed-parameter deformation complex \cite{BB04} to produce a splitting for the extended deformation complex.  This allows us to construct local coordinates on Coulomb slices in $\mathscr C$ orthogonal to the gauge orbits, which in turn define smooth local charts on $\mathscr M^{\mathrm{irr}}$.  We then give local trivializations for $q\from\mathscr E^{\mathrm{irr}}\to\mathscr M^{\mathrm{irr}}$, showing it is a principal bundle, and prove that $\pi\from\mathscr M^{\mathrm{irr}}\to\Lambda$ is a submersion.

Let $\mathscr M_\delta$ and $\mathscr E_\delta$ be the spaces constructed for a given choice of Sobolev weight $\delta>0$.  The smoothness arguments require the Sobolev weight to be smaller than the smallest gap between parabolic weights at each puncture, considered as elements of $\R/\Z$. To allow the parabolic weights to vary within the generic locus, we show that the families $\{\mathscr M_\delta\}_\delta$, $\{\mathscr E_\delta\}_\delta$ form directed systems and consider the direct limits
\[
    \mathscr M=\lim_{\delta\to0}\mathscr M_\delta\qquad\mathscr E=\lim_{\delta\to0}\mathscr E_\delta.
\]

Since the later sections use the algebraic description for Higgs bundles, we provide a bridge between the algebraic and analytic perspectives. The existence of harmonic metrics for stable parabolic Higgs bundles translates the algebraic moduli family into our analytic moduli family. To turn this into a diffeomorphism, we control the harmonic metrics along $1$-parameter families of Higgs bundles with varying boundary data, following Kim and Wilkin's proof in \cite{KW18} for varying parabolic weights with fixed residues.

\Cref{sec: Homology Transport} establishes local triviality of the family $\pi\from\mathscr M\to\Lambda$ over $\Lambda_0\subset\Lambda$, where semistability implies stability. $\Lambda_0$ is obtained from $\Lambda$ by removing a codimension-3 subset. Due to the noncompactness of the fibers of $\pi$, a horizontal lift of a path in $\Lambda_0$ could a-priori escape to infinity in finite time. We construct lifts locally for small open sets in $\Lambda_0$, and show that they give a well-defined local system for the homology groups $H_2(\mathcal M_\nu,\C)$. This proves the topological hypothesis of the affine variation theorem. 

Working in the algebraic setting, we begin with fixed generic parabolic weights, only allowing the Higgs field residues to vary in the vector space $V$ of possible residue data.  Denote this family $\pi_\alpha\from\mathscr M_\alpha^\mathrm{alg}\to V$.  We use the $\C^\times$-action to choose lifts of paths in $V$ whose images under the Hitchin map $\mathrm{Hit}\from\mathscr M_\alpha^\mathrm{alg}\to\mathcal B$ stay bounded. Properness of the joint map $(\pi_\alpha,\mathrm{Hit})$, established by Lee and Lee in \cite{LL25}, ensures lifted paths remain in a compact set.  We extend these local trivializations to the full generic parameter locus $\Lambda_0$ by comparing stability conditions as the weights vary.  The resulting integral homology local systems have trivial monodromy because $\Lambda_0$ is simply connected.  This completes the hypotheses needed for \Cref{thm: DH after partial reduction}, which we use to prove the statement for Higgs bundles \Cref{thm: general period variation}.

\Cref{sec: affine fiber alternative} calculates the period maps explicitly in the four-dimensional moduli spaces of exceptional types $E_6,E_7,E_8$ over $C=\C\P^1$, extending the surjectivity result of \cite{FMSW26} for the $D_4$ case.  \Cref{thm: DH after partial reduction} states that the period map is affine-linear in $\nu\in\Lambda_0$.  Choosing a basis for the cohomology of $\cM_\nu=\mathcal M(\alpha,\sigma)$ and transporting throughout $\Lambda_0$, we can express the period map as
\[
    \mathcal P(\nu)=\mathbf A\nu+\mathbf b
\]
for a matrix $\mathbf A$ whose entries are determined by the curvature class $c=c(\mathscr E)$, and a constant class $\mathbf b$.  To determine $\mathbf A$, it suffices to work in a single reference fiber, since parallel transport commutes with restriction to fibers (see \eqref{eq: DH homology compatibility}). We pick parabolic weights sufficiently small to admit an embedding from a star-shaped Nakajima quiver variety.

For a star-shaped quiver $Q$ and dimension vector $\vec r$ determined by the flag dimensions, there is, for suitable $\alpha$, an embedding
\[
    \cT:\cX(\beta,\tau)\hookrightarrow\cM(\alpha,\sigma)
\]
from a Nakajima quiver variety, where
\[
    \beta_i^{(k)}=\alpha_i^{(k)}-\alpha_i^{(k+1)},
    \qquad
    \tau_i^{(k)}=\sigma_i^{(k)}-\sigma_i^{(k+1)}.
\]
This embedding was constructed by Godinho--Mandini in the strongly parabolic rank~2 case \cite{GM11}, and generalized in \cite{RS21,Yae26}. Its image is the locus $\cM_0$ of Higgs bundles with trivial underlying holomorphic structure. Its complement is a Hitchin section $\mathcal S_p\cong\C$.The map $\cT$ intertwines the holomorphic symplectic forms up to a nonzero scalar depending on conventions. The quiver description determines the restriction of the boundary characteristic class to $\cM_0$. To determine the class on the full moduli space, we identify the homology class missing from this open subset.

The ALE topology of the quiver variety, pushed forward into $\cM$, gives us $n$ homology generators $\eta_j$ with intersection matrix equal to the negative $E_n$ Cartan matrix, all of which pair trivially with the class $[F]$ represented by a smooth Hitchin fiber. In particular, the second Betti numbers are $7,8,9$ for $E_6,E_7,E_8$, respectively.  We compute the characteristic class on the quiver generators using the determinants of the tautological flag bundles $F_i^{(k)}$.  Under the quiver embedding, their Chern classes are related to the tautological vertex bundles for the quiver variety, whose pairings with the homology generators are well understood.  In particular, we have $\langle c_\nu,[F]\rangle=0$, so the fiber periods $\langle[\omega_{s,\nu}],[F]\rangle$ are constant in $\nu$.

Transporting the quiver generators $\eta_v$ together with $[F]$ gives an integral homology basis throughout $\Lambda_0$.  Here $v$ runs over the vertices of the finite $E_n$ diagram obtained by deleting a rank-one vertex from the affine quiver (see \Cref{fig: E_n flag arrangements}).  Let $\dot\nu$ be a variation of $\nu$ encoding a variation $(\dot\alpha,\dot\sigma)$ of $(\alpha,\sigma)$, with $(\dot\beta,\dot\tau)$ the first variation of the corresponding quiver parameters extended to the central vertex by
\[
    \dot\beta_\star=\sum_i\dot\alpha_i^{(\ell_i)},
    \qquad
    \dot\tau_\star=\sum_i\dot\sigma_i^{(\ell_i)}.
\]
Then \Cref{prop: affine period slopes} states the period slopes
\[
    d\mathcal P_{\eta_v}(\dot\nu)
    =
    \left(-4\pi^2\,d\dot\beta_v,\,
          -8\pi^2\,\Re d\dot\tau_v,\,
          -8\pi^2\,\Im d\dot\tau_v\right).
\]

This theorem tells us the matrix $\mathbf A$ in the period formula $\mathcal P(\nu)=\mathbf A\nu+\mathbf b$.  To find $\mathbf b$, we first calculate the $\eta_v$ periods in the limit $(\alpha,0)\to(0,0)$ using the $\U(1)$-action on $\mathcal M(\alpha,0)$ and its Morse--Bott moment map.  Then, we work near a pair $(\alpha,0)$ where the moduli space becomes an orbifold, applying the singular elliptic curve description of \cite{Gro14} to prove $\bigl\langle[\omega_1],[F]\bigr\rangle=4\pi^2$ and $\bigl\langle[\Omega],[F]\bigr\rangle=0$.
\begin{namedtheorem}{\Cref{thm: periods}}
    Let $\nu\in\Lambda_0$ encode a weights and residues $(\alpha,\sigma)$, and define $(\beta,\tau)$ as in \Cref{def: quiver parameters}.  Let $\mathfrak B$ be the admissible basis for $H_2(\mathcal M(\alpha,\sigma),\C)$ from \Cref{sec: period variation}, indexed by the affine Dynkin vertex set $\tilde I$ with distinguished vertex $a=v_1^{(0)}$.  For all $v\in\tilde I$,
    \begin{equation*}\tag{\ref{eq: calibrated quiver periods}}
        \int_{\eta_v}\omega_1=4\pi^2\delta_{va}-4\pi^2\beta_v,
        \qquad
        \int_{\eta_v}\Omega=-8\pi^2\tau_v.
    \end{equation*}
\end{namedtheorem}
Finally, \Cref{thm: almost-global surjectivity} makes this an almost-global surjectivity statement onto the space of allowable Torelli parameters in the corresponding family of ALG spaces.

\subsection{Relation to Other Work}

The standard treatment of moduli spaces of parabolic Higgs bundles fixes a marked Riemann surface $(C,D)$, flags at the markings, as well as parabolic weights and complex residues compatible with the flags. The question of constructing well-behaved moduli spaces of Higgs bundles where at least one of the above structures is allowed to vary has recently received considerable attention \cite{CFW24, HSZ25, CTW25, ACGGT25, Hit26}. 

The particular case of allowing the weights and residues to vary studied in the current manuscript is at least partially motivated by its intimate relationship to the \textbf{Modularization conjecture}:
\begin{conjecture}
Every $\mathrm{ALG}/\mathrm{ALG}^{\star}$ hyperk\"ahler gravitational instanton arises as the $L^2$-metric of a moduli space of possibly wild parabolic Higgs bundles. 
\end{conjecture}
In \cite{FMSW26}, the authors computed the period triples of the $12$-parameter family of $\mathrm{SU}(2)$-Hitchin moduli spaces on the four-punctured sphere and showed that the values realize an open dense subset of the admissible Torelli data of the $12$-parameter family of ALG-$D_4$ gravitational instantons, thus essentially proving the Modularization conjecture modulo showing that the $L^2$-metrics have finite energy $\int_\mathcal{M} |\mathrm{Rm}|^2 < \infty$. The failure of actual surjectivity is because for particular values of the weight the associated moduli space in the literature becomes zero-dimensional and as a result does not correspond to an instanton; however, there should be a different four-dimensional Hitchin moduli space obtained by quotienting by a smaller gauge group which one can instead use.

The main application of the current paper is to establish the analogous result for $E_6,E_7,E_8$ which, modulo showing that the associated $L^2$-metrics have finite energy, essentially settles the Modularization conjecture for all ALG metrics associated to moduli spaces of tame parabolic Higgs bundles. The analytic tools we develop are influenced by the ones appearing in \cite{CFW24}, which was the first paper carefully treating the moduli space problem when the complex residues are allowed to vary. 

At this point it makes sense to explain how our treatment differs from \cite{CFW24} where the authors construct the framed moduli space where the residues are allowed to vary and only allow full flags. In our paper we only retain a determinant unitary framing which gives rise to the real torus bundle $\mathscr E$. We suspect that both constructions can lead to \emph{incomplete} hyperk\"ahler metrics. One example of this is the moduli space with extended Hitchin base $\C_\sigma$ parameterizing $q=-\det \varphi =(z^2 + \sigma) \de z^2$ on $\CP^1$ with an irregular singularity at infinity. The eigenvalues at infinity in coordinate $\widetilde{z}=z^{-1}$ are 
\[
    \lambda_\pm
    =\pm\left(
        \frac{1}{\widetilde z^3}
        +\frac{\sigma}{2\widetilde z}
        +O(\widetilde z)
    \right)\de\widetilde z.
\]
so one sees that changing $\sigma$ corresponds to changing the residue of the the eigenvalue associated to $z=\infty$, while leaving the higher-order polar part fixed. Gaiotto--Moore--Neitzke implicitly conjecture that this case should recover the well-known Ooguri-Vafa metric. Progress towards this perspective appears in the work of Tulli \cite{Tul25} and Nackan \cite{Nac24}. We conjecture that after fixing a compact subset $K\subset \Lambda\times \mathcal B$ of the product of the parameter space times the Hitchin base, its pre-image in $\mathscr E$ should also admit a hyperk\"ahler metric.

\subsection{Weighted Sobolev Spaces}\label{sec: Weighted Sobolev Spaces} In this subsection we state our conventions for the weighted Sobolev spaces appearing in the paper and some estimates used explicitly or implicitly in various places. 

Write $C^\circ=C\sm D$, and choose disjoint coordinate disks centered at the punctures, with $z=re^{\I\theta}$ and $B_\eps(p_i)=\{|z|<\eps\}$ for $0<\eps\leq1$.  Write $B_\eps^\circ(p_i)=B_\eps(p_i)\sm\{p_i\}$ for the punctured disk.  Fix once and for all a smooth metric $g_C$ compatible with the complex structure on $C$.  In a local coordinate $z$, write $g_C=\kappa(z)|\de z|^2$.

Let $\hat C\to C$ be the real oriented blowup along $D$.  Its $b$-vector fields are the smooth vector fields tangent to the boundary, locally generated by $r\partial_r,\partial_\theta$ near the boundary and generating the usual tangent sheaf away from it.  Together with smooth functions they generate a differential-operator algebra $\mathcal D_b$.

For a sheaf of sections $\mathscr S$ of a vector bundle on $\hat C$ with smooth fiber metric and connection, let $\mathcal H_b^q(U,\mathscr S)$ be the Sobolev completion of its smooth sections of finite norm over $U\subset\hat C$, using covariant derivatives of order up to $q$ with respect to $b$-vector fields and the surface measure $\de\mathrm{vol}_{g_C}$.  Thus the measure near a puncture is $\kappa r\,\de r\,\de\theta$, rather than the cylindrical measure.

Choose positive functions $\rho_i$, smooth on $C^\circ$, equal to $r$ on $B_1^\circ(p_i)$ and to $1$ outside slightly larger disjoint coordinate disks, and set $\rho=\prod_i\rho_i$.  For $\delta\in\R$, define
\begin{equation}\label{eq: weighted Sobolev norm}
    \mathcal H_\delta^q(U,\mathscr S)=\rho^{\delta+1}\mathcal H_b^q(U,\mathscr S),\qquad
    \|u\|_{\mathcal H_\delta^q}=\|\rho^{-\delta-1}u\|_{\mathcal H_b^q}.
\end{equation}
The exponent is chosen so that $r^a\in \mathcal H_\delta^q(B_\eps^\circ(p_i))$ if and only if $a>\delta$.
We sometimes omit $U,\mathscr S$ when the context is clear, and use the same notation for a puncture neighborhood and its preimage in $\hat C$.  For sheaves of differential forms, we work with the $b$-cotangent bundle locally generated by $\de r/r,\de\theta$, and write $\Omega_{\hat C}^q$ for its sheaf of smooth exterior $q$-forms.

We will use these spaces to define configuration space and gauge group in terms of model singular parts and decaying remainders.  In a smooth local frame near $p_i$, the norm is equivalent to
\begin{equation}\label{eq: cylindrical Sobolev norm}
    \|u\|_{\mathcal H_\delta^q}^2\sim
    \sum_{a+b\leq q}\int_0^1\int_0^{2\pi}
        r^{-2\delta}\bigl|(r\partial_r)^a\partial_\theta^b u\bigr|^2
        \,\frac{\de r}{r}\,\de\theta.
\end{equation}
These spaces will be defined in terms of sections of vector bundles equipped with a background metric, giving the norm used above.

Much of our analytic construction is inspired by \cite{BB04}, but we use $b$-derivatives in place of their weighted covariant derivatives.  At order zero, our $\mathcal H_\delta^0$ agrees with their $L_\delta^2$ on functions, while our coefficient convention on 1-forms gives $\mathcal H_\delta^0(\Omega_{\hat C}^1)=L_{\delta-1}^2(\Omega^1)$, with equivalent norms.

\noindent\textbf{Useful Estimates:}
    \begin{itemize}
        \addtolength{\itemsep}{4pt}
        \item If $q>m+1$, then the derivatives $(r\partial_r)^a\partial_\theta^b u$ with $a+b\leq m$ are continuous and
        \begin{equation}\label{eq: weighted Sobolev embedding}
            \sum_{a+b\leq m}\sup_{B_1^\circ(p_i)}
                r^{-\delta}\bigl|(r\partial_r)^a\partial_\theta^b u\bigr|
            \leq C\|u\|_{\mathcal H_\delta^q}.
        \end{equation}

        \item If $\delta\geq\delta'$ and $q\geq q'$, there is a continuous inclusion
        \begin{equation}\label{eq: weighted Sobolev inclusion}
            \mathcal H_\delta^q\into\mathcal H_{\delta'}^{q'}.
        \end{equation}
        An operator $P\in\mathcal D_b$ of order $j\leq q$ acts continuously from $\mathcal H_\delta^q$ to $\mathcal H_\delta^{q-j}$.

        \item For $k\geq2$, pointwise products, tensor contractions, and wedge products satisfy
        \begin{equation}\label{eq: weighted Sobolev multiplication}
            \|uv\|_{\mathcal H_{\delta+\eta}^k}
            \leq C\|u\|_{\mathcal H_\delta^k}\|v\|_{\mathcal H_\eta^k}.
        \end{equation}
        Smooth coefficients on $\hat C$ act as bounded multipliers.

        \item With $\de V_b=\rho^{-2}\de\mathrm{vol}_{g_C}$, weighted Cauchy--Schwarz gives
        \begin{equation}\label{eq: weighted Cauchy Schwarz}
            \int_{C^\circ}|\langle u,v\rangle|\,\de V_b
            \leq \|u\|_{\mathcal H_\delta^0}\|v\|_{\mathcal H_{-\delta}^0}.
        \end{equation}

        \item For $\delta>0$ and $f\in\mathcal H_\delta^q(B_1^\circ(p_i))$, we have
        \begin{equation}\label{eq: radial integration}
            \int_0^r f(s,\theta)\,\frac{\de s}s\in\mathcal H_\delta^q(B_1^\circ(p_i)).
        \end{equation}

        \item For a conformal metric locally of the form $g_C=\kappa(\de r^2+r^2\de\theta^2)$,
        \begin{equation}\label{eq: weighted Hodge star}
            \star\frac{\de r}{r}=\de\theta,\qquad
            \star\de\theta=-\frac{\de r}{r},\qquad
            \star\left(\frac{\de r}{r}\wedge\de\theta\right)
                =\frac{1}{\kappa r^2},\qquad
            \star1=\kappa r^2\frac{\de r}{r}\wedge\de\theta.
        \end{equation}

    \end{itemize}

\subsection*{Acknowledgements}

We thank Rafe Mazzeo, Andy Neitzke, Richard Wentworth, and Max Zimet for inspiring conversations over the years which refined the ideas presented in this paper. 
LF is partially supported by NSF DMS-2005258.

\smallskip\noindent
\emph{Use of AI.}
We used ChatGPT's \emph{Astra 6} and DeepSeek V4.1 Flash to help digest existing literature, correct typographical mistakes, and write code for our figures. The ideas are our own. The authors are responsible for the mathematical
content and the accuracy of the references.

\section{A Generalized Duistermaat--Heckman Theorem}\label{sec: a generalized DH theorem}

Consider a family $\{\mathcal M_\alpha\}_\alpha$ of moduli spaces for $\alpha$ in some parameter space $\Lambda$, with each $\mathcal M_\alpha$ equipped with a symplectic form $\omega_\alpha$.  The goal of this section is to describe how the cohomology class of $\omega_\alpha$ depends on the choice of $\alpha$, which we do by adapting Duistermaat and Heckman's theorem \cite[Theorem~1.1]{DH82} to have suitably weak hypotheses that can apply to noncompact quotients of infinite-dimensional configuration spaces.  Namely, we give a geometric hypothesis and a topological hypothesis.

Assume the family of moduli spaces can be joined into a smooth family $\pi\from\mathscr M\to\Lambda$, so that $\mathcal M_\alpha=\pi^{-1}(\alpha)$.  Under ideal circumstances, the family arises from a single quotient of a single configuration space $\mathscr C$ by an action $G\acts\mathscr C$ which is Hamiltonian with respect to a symplectic form $\omega$ on $\mathscr C$.  In our convention, the moment map $\mu\from\mathscr C\to\mathfrak g^*$ satisfies
\begin{equation}\label{eq: DH Hamiltonian identity}
    \iota_{\zeta_\xi}\omega=\de\langle\mu,\xi\rangle
\end{equation}
where $\zeta_\xi$ is any vector field on $\mathscr C$ generated by some $\xi\in\mathfrak{g}$.  In this situation, $\Lambda\subset\mathfrak g^*$ and $\mu$ descends to $\pi$ after restricting to $\mu^{-1}(\Lambda)$ and taking quotients.  Also ideally, $\pi\from\mathscr M\to\Lambda$ admits a complete Ehresmann connection, allowing one to lift paths in $\Lambda$ to paths in $\mathscr M$.  This gives identifications between moduli spaces used to compare the classes $[\omega_\alpha]$.  \Cref{sec: multi Taub NUT} demonstrates this ideal situation in the well-studied multi-Taub-NUT family.

In the case of parabolic Higgs bundle moduli spaces, existence of such an $\omega$ and completeness of an Ehresmann connection are substantial obstacles.  This motivates our careful setup for \Cref{thm: DH after partial reduction}, which softens both hypotheses.  The main features are local systems for homology, and a procedure for reducing to finite dimensions.

\subsection{Affine Period Variation Theorem}

Throughout this subsection, we work over $\mathbb F=\R$ or $\mathbb F=\C$ and take all Lie groups and Lie algebras to be respectively real or complex.  We also assume all differential forms are $\mathbb F$-valued.  We begin with a definition pertaining to the topological hypothesis of \Cref{thm: DH after partial reduction}.

\begin{definition}\label{def: compatible local system}
    Let $\pi\from\mathscr M\to\Lambda$ be a continuous family of smooth manifolds $\mathcal M_\alpha=\pi^{-1}(\alpha)$, with inclusions $i_\alpha\from\mathcal M_\alpha\into\mathscr M$.  A \emph{homology local system} for this family consists of the vector spaces $H_2(\mathcal M_\alpha,\mathbb F)$ and parallel transport isomorphisms
    \[
        \mathrm{PT}_\gamma\from H_2(\mathcal M_{\alpha_*},\mathbb F)\longrightarrow H_2(\mathcal M_\alpha,\mathbb F)
    \]
    for paths $\gamma$ in $\Lambda$ from $\alpha_*$ to $\alpha$.  Transport depends only on homotopy relative to endpoints, respects concatenation, and is the identity on constant paths.  The homology local system is \emph{compatible} with $\pi$ if, for every path $\gamma$ in $\Lambda$, the following diagram commutes:
    \begin{equation}\label{eq: DH homology compatibility}
        \begin{tikzcd}[ampersand replacement=\&,column sep=large,row sep=large]
            H_2(\mathcal M_{\alpha_*},\mathbb F) \arrow[rr,"\mathrm{PT}_\gamma"] \arrow[dr,"(i_{\alpha_*})_*"']
                \&\& H_2(\mathcal M_\alpha,\mathbb F) \arrow[dl,"(i_\alpha)_*"] \\
                \& H_2(\mathscr M,\mathbb F) \&
        \end{tikzcd}
    \end{equation}
    Thus a cycle and its transport are homologous in the total space $\mathscr M$.  We write $\mathrm{PT}_\gamma^*$ for the dual transport on degree-two cohomology.
\end{definition}

For many moduli spaces of interest, one can take partial quotients to reduce from infinite dimensions to to finite dimensions.  The first advantage of this reduction is that we can interpret the coefficients of the affine period map in terms of characteristic classes of a finite-dimensional principal bundle $E$ over $\mathscr M$.  The second advantage is that the geometric hypothesis of \Cref{thm: DH after partial reduction} can be purely stated in finite-dimensional terms.  \Cref{lem: partial quotient} below shows how the ideal situation described above automatically gives the structure needed for geometric hypothesis.

Let $G\acts (M,\omega)$ with moment map $\mu$ as in this section's introduction.  Fix a finite-dimensional parameter subspace $\Lambda\subseteq(\mathfrak g^*)^G$, and a closed normal Lie subgroup $H\subseteq G$ whose Lie algebra annihilates $\Lambda$, such that $B=G/H$ is finite dimensional.  Set $E=\mu^{-1}(\Lambda)/H$, and let $H^{\mathrm{eff}}$ denote the quotient of $H$ by the kernel of its action on $\mu^{-1}(\Lambda)$.  Suppose $\mu^{-1}(\Lambda)$ is smooth and $p\from\mu^{-1}(\Lambda)\longrightarrow E$ is a smooth principal $H^{\mathrm{eff}}$-bundle.

\begin{lemma}\label{lem: partial quotient}
    In the above setting, the restriction of $\omega$ to $\mu^{-1}(\Lambda)$ descends to a closed $B$-invariant $2$-form $\varpi$ on the partial quotient $E$, and $\mu$ descends to $\tilde\pi\from E\to\Lambda\subset\mathfrak{b}^*$ such that
    \begin{equation}\label{eq: DH partial quotient identity}
        \iota_{\zeta_b}\varpi=\de\langle\tilde\pi,b\rangle,\qquad b\in\mathfrak b.
    \end{equation}
    Moreover, $\tilde\pi^{-1}(\alpha)/B=\mu^{-1}(\alpha)/G=\mathcal{M}_{\alpha}$.
\end{lemma}
\begin{proof}
    Let $j\from\mu^{-1}(\Lambda)\into M$.  For $\xi\in\mathfrak h$, the function $\langle\mu,\xi\rangle$ vanishes on $\mu^{-1}(\Lambda)$ because $\mathfrak h$ annihilates $\Lambda$, so $\iota_{\zeta_\xi} j^*\omega=0$.  Thus $j^*\omega$ is $H$-basic and descends to a closed $B$-invariant form $\varpi$.  Since $\mu$ is constant on $G$-orbits in $\mu^{-1}(\Lambda)$, it descends to the map $\tilde\pi$.  For a lift $\xi\in\mathfrak g$ of $b\in\mathfrak b$,
    \begin{equation}
        p^*(\iota_{\zeta_b}\varpi)
        =\iota_{\zeta_\xi} j^*\omega
        =p^*\de\langle\tilde\pi,b\rangle.
    \end{equation}
    The pairing is independent of the lift, again because $\mathfrak h$ annihilates $\Lambda$.  The quotient identification follows by taking subsequent quotients by $H$ and then by $G/H$.
\end{proof}

The conclusion of the lemma is the setting in which we state \Cref{thm: DH after partial reduction}.  Having reduced, everything should be taken to be finite-dimensional over $\mathbb F=\R$ or $\mathbb F=\C$.

\begin{theorem}\label{thm: DH after partial reduction}
Let $\mathscr M$ be a paracompact manifold, let $B$ be a Lie group, and fix a coadjoint-invariant subspace $\Lambda\subseteq(\mathfrak b^*)^B$.  Suppose $q\from E\to\mathscr M$ is a principal $B$-bundle, and $\tilde\pi\from E\to\Lambda$ is a $B$-invariant map descending to a map $\pi\from\mathscr M\to\Lambda$.

Suppose there is an open region $\Lambda_0\subseteq\Lambda$ such that the restrictions $E_0=\tilde\pi^{-1}(\Lambda_0)$ and $\mathscr M_0=\pi^{-1}(\Lambda_0)$ are smooth, and all the maps above are smooth after restriction.  Assume further that $\pi\from\mathscr M_0\to\Lambda_0$ is a submersion, and the fibers $\mathcal M_\alpha$ have symplectic forms $\omega_\alpha$.  Write $i_\alpha$, $j_\alpha$, and $q_\alpha$ for the inclusions and submersion in the following diagram:
\begin{equation}\label{eq: partial quotient diagram}
    \begin{tikzcd}[ampersand replacement=\&,column sep=large,row sep=large]
        \tilde\pi^{-1}(\alpha) \arrow[d,two heads,"q_\alpha"] \arrow[r,hook,"j_\alpha"]
            \& E_0 \arrow[d,two heads,"q"'] \arrow[dr,"\tilde\pi"] \& \\
        \mathcal M_\alpha \arrow[r,hook,"i_\alpha"]
            \& \mathscr M_0 \arrow[r,two heads,"\pi"'] \& \Lambda_0 .
    \end{tikzcd}
\end{equation}
We additionally assume the following:
    \begin{itemize}
        \item The \emph{geometric hypothesis}: there is a closed, $B$-invariant, $\mathbb F$-valued two-form $\varpi$ on $E_0$ satisfying
        \begin{equation}\label{eq: Hamiltonian identity}
            \iota_{\zeta_b}\varpi=\de\langle\tilde\pi,b\rangle,\qquad b\in\mathfrak b,
        \end{equation}
        where $\zeta_b$ is the vector field generated by $b$.  For each $\alpha\in\Lambda_0$, the restriction of $\varpi$ to $\tilde\pi^{-1}(\alpha)$ descends to $\omega_\alpha$.  
        \item The \emph{topological hypothesis}: the family $(\mathcal M_\alpha)_{\alpha\in\Lambda_0}$ admits a homology local system compatible with $\pi\from\mathscr M_0\to\Lambda_0$.
    \end{itemize}
    Then there is a class $c\in H^2(\mathscr M_0,\mathbb F)\otimes_{\mathbb F}\Lambda^*$ such that, for every $\alpha_0,\alpha_1\in\Lambda_0$ and every path $\gamma$ in $\Lambda_0$ from $\alpha_0$ to $\alpha_1$,
    \begin{equation}\label{eq: main theorem as parallel transport}
        \mathrm{PT}_\gamma^*[\omega_{\alpha_1}]
        =[\omega_{\alpha_0}]
        -\langle c_{\alpha_0},\alpha_1-\alpha_0\rangle.
    \end{equation}
    Here $c$ is the curvature class of the principal $B$-bundle $q\from E_0\to\mathscr M_0$ and $c_\alpha=i_\alpha^*c$.
\end{theorem}

\begin{proof}
    This proof is an adaptation of \cite{DH82}.  The strategy is to modify $\varpi$ by an exact form so that it descends to $\mathscr M_0$.  Choose a smooth principal connection $\theta$ on $q\from E_0\to\mathscr M_0$ satisfying $\theta(\zeta_b)=b$.  Such a connection exists by standard partition of unity arguments.

    Let $r\from\mathfrak b\to\Lambda^*$ be the evaluation map, and set $L=r\circ\theta$.  Since $\Lambda$ is coadjoint-invariant, $L$ is a $B$-invariant one-form satisfying $L(\zeta_b)=r(b)$.  Cartan's formula gives
    \begin{equation}
        \iota_{\zeta_b}\de L
        =\mathcal L_{\zeta_b}L-\de(L(\zeta_b))
        =-\de(r(b))=0.
    \end{equation}
    Hence $\de L$ descends to a closed $\Lambda^*$-valued two-form $\Omega$ on $\mathscr M_0$ satisfying
    \begin{equation}\label{eq: DH projected curvature}
        q^*\Omega=\de L=\de(r\circ\theta)=r(F_\theta),
    \end{equation}
    where $F_\theta=\de\theta+\frac12[\theta\wedge\theta]$ is the curvature of $\theta$, and the last equality uses the fact that $r$ vanishes on Lie brackets by adjoint-invariance.  Changing $\theta$ changes $L$ by a basic one-form and hence $\Omega$ by an exact form.  Thus the projected curvature class $c(E_0)=[\Omega]$ is independent of $\theta$.

    The modified 2-form $\varpi+\de\langle L,\tilde\pi\rangle$ is closed and $B$-invariant.  We also have
    \begin{equation}
        \iota_{\zeta_b}\de\langle L,\tilde\pi\rangle
        =-\de\langle L(\zeta_b),\tilde\pi\rangle
        =-\de\langle\tilde\pi,b\rangle
        =-\iota_{\zeta_b}\varpi,
    \end{equation}
    which implies $\varpi+\de\langle L,\tilde\pi\rangle$ is basic and descends to a closed two-form $\beta$ on $\mathscr M_0$ with
    \[
        q^*\beta=\varpi+\de\langle L,\tilde\pi\rangle.
    \]
    To restrict this identity to the fiber over $\alpha$, consider the commutative square
    \begin{equation}\label{eq: DH fiber restriction diagram}
        \begin{tikzcd}[ampersand replacement=\&,column sep=large,row sep=large]
            \tilde\pi^{-1}(\alpha) \arrow[d,two heads,"q_\alpha"'] \arrow[r,hook,"j_\alpha"]
                \& E_0 \arrow[d,two heads,"q"] \\
            \mathcal M_\alpha \arrow[r,hook,"i_\alpha"']
                \& \mathscr M_0 .
        \end{tikzcd}
    \end{equation}
    The restriction of $\de\langle L,\tilde\pi\rangle$ to $\tilde\pi^{-1}(\alpha)$ is $j_\alpha^*\langle\de L,\alpha\rangle$, and by the geometric hypothesis $j_\alpha^*\varpi=q_\alpha^*\omega_\alpha$.  Thus
    \begin{equation}
        j_\alpha^*q^*\beta=q_\alpha^*\omega_\alpha+j_\alpha^*\langle \de L,\alpha\rangle=q_\alpha^*\omega_\alpha+j_\alpha^*q^*\langle\Omega,\alpha\rangle.
    \end{equation}
    The diagram \eqref{eq: DH fiber restriction diagram} says $j_\alpha^*q^*=q_\alpha^*i_\alpha^*$.  Since $q_\alpha$ is a surjective submersion, $q_\alpha^*$ is injective, so the above implies
    \begin{equation}\label{eq: basic identity}
        \omega_\alpha
        =i_{\alpha}^*\beta
        -\langle i_{\alpha}^*\Omega,\alpha\rangle.
    \end{equation}

    Let $c=c(E_0)$ be the curvature class, and write $c_\alpha=i_\alpha^*c=i_\alpha^*[\Omega]$.  By the topological hypothesis, parallel transport commutes with the family of maps $(i_\alpha)_*\from H_2(\mathcal M_\alpha,\mathbb F)\to H_2(\mathscr M_0,\mathbb F)$ (see \eqref{eq: DH homology compatibility}), and duality gives the analogous statement for cohomology.  For a path $\gamma$ from $\alpha_0$ to $\alpha_1$ in $\Lambda_0$,
    \begin{equation}\label{eq: DH transported classes}
        \mathrm{PT}_\gamma^*i_{\alpha_1}^*[\beta]=i_{\alpha_0}^*[\beta],
        \qquad
        \mathrm{PT}_\gamma^*c_{\alpha_1}=c_{\alpha_0}.
    \end{equation}
    Applying these identities to \eqref{eq: basic identity} gives
    \[
        \begin{aligned}
            \mathrm{PT}_\gamma^*[\omega_{\alpha_1}]
                &=i_{\alpha_0}^*[\beta]-\langle c_{\alpha_0},\alpha_1\rangle\\
                &=[\omega_{\alpha_0}]+\langle c_{\alpha_0},\alpha_0\rangle
                    -\langle c_{\alpha_0},\alpha_1\rangle\\
                &=[\omega_{\alpha_0}]-\langle c_{\alpha_0},\alpha_1-\alpha_0\rangle,
        \end{aligned}
    \]
    as required.
\end{proof}

A stronger, sufficient condition for the topological hypothesis is the existence of a complete Ehresmann connection.  An Ehresmann connection on $\pi\from\mathscr M_0\to\Lambda_0$ consists of horizontal tangent spaces complementary to $\ker\de\pi$, giving uniquely-defined lifts of vector fields on $\Lambda_0$ to $\mathscr M_0$.  It is complete if horizontal lifts of complete vector fields remain complete.

\begin{corollary}\label{cor: complete Ehresmann}
    Under the geometric hypothesis of \Cref{thm: DH after partial reduction}, suppose $\pi\from\mathscr M_0\to \Lambda_0$ admits a smooth and complete Ehresmann connection.  Then $\pi$ is a fiber bundle and satisfies the topological hypothesis of the theorem.
\end{corollary}
\begin{proof}
    Fix a convex neighborhood $U$ of a parameter $\alpha_*$.  Horizontal transport along straight segments from $\alpha_*$ gives a smooth trivialization over $U$.  By homotopy invariance, the induced transport on homology depends only on the homotopy class of the path relative to its endpoints, so the fiber homology groups form a local system.  Along any path, horizontal transport gives a homotopy between the inclusion of the initial fiber and the transported inclusion, so the induced homology transport is compatible with the fiber inclusions.
\end{proof}

\begin{remark}[Hyperk\"ahler reduction]\label{rem: hyperkahler DH}
    For a family of hyperk\"ahler moduli spaces, let $\Lambda_0\subseteq \Lambda^{\oplus3}$, write $\underline\alpha=(\alpha_1,\alpha_2,\alpha_3)$, and suppose the partial quotient $E$ has three closed $B$-invariant forms $\varpi_1,\varpi_2,\varpi_3$ satisfying $\iota_{\zeta_b}\varpi_s=\de\langle\tilde\pi_s,b\rangle$, where $(\tilde\pi_1,\tilde\pi_2,\tilde\pi_3)$ is the parameter map on the partial quotient $E$, and the induced forms $\omega_{s,\underline\alpha}$ are the three K\"ahler forms on the fibers.  Apply the same proof to each $\varpi_s+\de\langle L,\tilde\pi_s\rangle$, where now the topological hypothesis must hold over the tripled parameter space.  Along a path $\gamma$ from $\underline\alpha_*$ to $\underline\alpha$, this gives
    \begin{equation}\label{eq: DH hyperkahler variation}
        \mathrm{PT}_\gamma^*[\omega_{s,\underline\alpha}]
        =[\omega_{s,\underline\alpha_*}]
        -\langle c,\alpha_s-\alpha_{*,s}\rangle,
        \qquad s=1,2,3,
    \end{equation}
    with the same curvature class $c=c(E)$ for all three forms.  A complete Ehresmann connection on the family over $\Lambda_0$ again suffices for the topological hypothesis.
\end{remark}

\subsection{Multi-Taub--NUT Example}\label{sec: multi Taub NUT}
The multi-Taub--NUT spaces are the ALF-$A_k$ gravitational instantons given by the Gibbons--Hawking ansatz.  They are noncompact and hence outside the scope of Duistermaat--Heckman's original proof, but their periods can be directly shown to vary affine-linearly with the relative positions of their $k+1$ centers.  We recall their hyperk\"ahler quotient construction (see \cite[Section~3]{GR97} for details), which is already finite dimensional, so no partial reduction is needed.  We apply \Cref{thm: DH after partial reduction} by proving completeness of a natural Ehresmann connection.

Fix $k\geq0$ and $\lambda>0$, and let
\[
    E=\mathbb H^{k+1}\times(\R^3\times S^1),\qquad B=T^{k+1},
\]
with quaternion coordinates $q_j$, Cartesian coordinates $\mathbf y=(y_1,y_2,y_3)$, and $\tau\in\R/2\pi\Z$.  Equip $E$ with the product flat hyperk\"ahler metric $g_E$, using $\lambda|\de\mathbf y|^2+\lambda^{-1}\de\tau^2$ on $\R^3\times S^1$.  Use $J_1,J_2,J_3$ given by left multiplication by $\I,\mathrm j,\mathrm k$ on the quaternion factors, and choose the structure on $\R^3\times S^1$ so that $\iota_{\pa_\tau}\varpi_s=\de y_s$.  The $j$th factor of $B$ acts by $(q_j,\tau)\mapsto(q_j e^{-\I\theta},\tau+\theta)$, fixing the other coordinates.  The hyperk\"ahler moment map is given by the triple
\begin{equation}\label{eq: Taub NUT moment map}
    \tilde\pi=(\tilde\pi_1,\tilde\pi_2,\tilde\pi_3)\from E\longrightarrow(\mathfrak b^*)^{\oplus3},
    \qquad \tilde\pi_s=(\tilde\pi_s^{(0)},\ldots,\tilde\pi_s^{(k)}),
\end{equation}
whose components satisfy
\[
    \tilde\pi_1^{(j)}\I+\tilde\pi_2^{(j)}\mathrm j+\tilde\pi_3^{(j)}\mathrm k
    =\tfrac12q_j\,\I\,\bar q_j+\mathbf y,
\]
where we identify $\mathbf y=y_1\I+y_2\mathrm j+y_3\mathrm k\in\mathbb H$.

Write $\underline\alpha=(\alpha_1,\alpha_2,\alpha_3)$ with $\alpha_s=(\alpha_s^{(0)},\ldots,\alpha_s^{(k)})$, and set $\alpha^{(j)}=(\alpha_1^{(j)},\alpha_2^{(j)},\alpha_3^{(j)})$.  Over
\[
    \Lambda_0=\left\{\underline\alpha\in(\mathfrak b^*)^{\oplus3}:
        \alpha^{(j)}\ne\alpha^{(\ell)}\text{ for }j\ne\ell\right\},
\]
the action on $E_0=\tilde\pi^{-1}(\Lambda_0)$ is free, giving a principal bundle $q\from E_0\to\mathscr M_0=E_0/B$.  The fibers of $\pi\from\mathscr M_0\to\Lambda_0$ are the multi-Taub--NUT spaces, with Gibbons--Hawking centers $\alpha^{(j)}\in\R^3$, the points over which the circle fibers collapse.

Lifting vectors in $\Lambda$ to vectors in $E_0$ requires a smooth choice of right inverses to $\de\tilde\pi$.  At a point $x\in E_0$, consider the operator $DD^\dagger$ where $D=\de\tilde\pi_x\from T_xE\to(\mathfrak b^*)^{\oplus3}$.  If $DD^\dagger$ is invertible, then $D^\dagger(DD^\dagger)^{-1}$ gives a natural choice of right inverse.  As $x$ varies, this defines an Ehresmann connection.  To bound the lifts and prove completeness, we must bound the eigenvalues of $DD^\dagger$ below over compact sets in $\Lambda$.  Let $u_j$ generate the $j$th $S^1$ factor of $B$.  The moment-map identity $\de\tilde\pi_s^{(j)}=\iota_{u_j}\varpi_s=g_E(J_su_j,-)$ gives
\[
    \begin{aligned}
        (DD^\dagger)_{(s,j),(r,\ell)}&=g_E(J_su_j,J_ru_\ell)=\delta_{sr}A_{j\ell},\\
        A_{j\ell}&=g_E(u_j,u_\ell)=\delta_{j\ell}|q_j|^2+\lambda^{-1},
    \end{aligned}
\]
which shows $DD^\dagger$ consists of three diagonal blocks, each equal to $A=(A_{j,\ell})$.

We now find a lower bound on the eigenvalues of $A$, uniform over an arbitrary compact $K\subset\Lambda_0$.  There is a constant $\delta>0$ such that $|\alpha^{(j)} - \alpha^{(\ell)}| \geq \delta$ for $j\ne\ell$ throughout $K$.  Since $|q_j|^2=2|\alpha^{(j)}-\mathbf y|$, at most one $|q_j|^2$ can be smaller than $\delta$.  Relabel the coordinates so that $|q_0|$ is smallest.  Then,
\begin{equation}\label{eq: Taub NUT coercivity}
    v^TAv=\sum_{j=0}^k|q_j|^2v_j^2
    +\lambda^{-1}\Bigl(\sum_{j=0}^kv_j\Bigr)^2 \geq \delta \sum_{j=1}^k v_j^2 + \lambda^{-1}\Bigl(\sum_{j=0}^kv_j\Bigr)^2.
\end{equation}
The right hand side is positive definite in $v$, so $A\geq c_K\operatorname{Id}$ on $\tilde\pi^{-1}(K)$ for some $c_K>0$.  With this estimate, Gr\"onwall's inequality shows that all paths in $\Lambda$ have well-defined lifts, which proves completeness.  By \Cref{cor: complete Ehresmann} and \Cref{rem: hyperkahler DH}, \Cref{thm: DH after partial reduction} gives affine period variation with respect to $\underline\alpha$.

\section{Universal Geometric Structure}\label{sec: Universal Moduli Space}

In this section we describe the geometric construction of the moduli family $\pi\from\mathscr M\to\Lambda$ of parabolic Higgs bundles and its natural principal torus bundle $\mathscr E$.  We endow $\mathscr E$ with an invariant triple of extended Hitchin forms, and recover the moduli spaces $\mathcal M(\alpha,\sigma)$ together with their k\"ahler forms through Hamiltonian reduction.  We then show that our construction satisfies the geometric hypothesis of \Cref{thm: DH after partial reduction}.

As in the introduction, we fix a compact connected Riemann surface $C$, a divisor of distinct marked points $D=\{p_1,\ldots,p_n\}$, and a rank-$r_\star$ bundle $E$ over $C$, together with a fixed identification $\det E\cong\mathcal O_C$.  We denote $C^\circ=C\sm D$ and $E^\circ$ the restriction to $C^\circ$.  We also fix parabolic filtrations over each point in $D$, and a background metric $h_0$ on $E$.

\subsection{The Gauge Group}\label{sec: The Gauge Group}

To describe the boundary gauge transformations, consider the associated graded space of the parabolic flag
\[
    \gr F_i=\bigoplus_{k=0}^{\ell_i}\gr F_i^{(k)},
\]
where $\gr F_i^{(k)}=F_i^{(k)}/F_i^{(k-1)}$.  Recall $r_i^{(k)}=\dim F_i^{(k)}$ and $m_i^{(k)}=\dim\gr F_i^{(k)}$.  Consider the complex Levi group
\[
    L_i=\left\{g\in\prod_{k=0}^{\ell_i}\mathrm{GL}(\gr F_i^{(k)})\ \Big|\ \det(g)=1\right\},
\]
and the real compact Levi group
\[
    K_i=\left\{g\in\prod_{k=0}^{\ell_i}U(\gr F_i^{(k)})\ \Big|\ \det(g)=1\right\}.
\]
We denote the respective Lie algebras by $\mathfrak l_i$ and $\mathfrak k_i$.  Note the latter has center
\begin{equation}
    \mathfrak z(\mathfrak k_i)
    =\left\{\bigoplus_{k=0}^{\ell_i}a_i^{(k)}\Id_{\gr F_i^{(k)}}:
    a_i^{(k)}\in\I\R,\quad
    \sum_{k=0}^{\ell_i}m_i^{(k)}a_i^{(k)}=0\right\}.
\end{equation}
Taking block determinants gives a homomorphism from $K_i$ onto the \emph{$i$th boundary torus}
\[
    T_i=\left\{(u_0,\dots,u_{\ell_i})\in U(1)^{\ell_i+1}\ \Big|\ \prod_{k=0}^{\ell_i}u_k=1\right\},
\]
and taking traces identifies $\mathfrak z(\mathfrak k_i)$ with the Lie algebra of $T_i$.

\bigskip

We now specify the gauge transformations allowed at a puncture.  Full details will be given in \Cref{sec: The Analytic Family}.  Choose once and for all orthonormal trivializations
\[
    E|_{B_1(p_i)}\cong B_1(p_i)\times\gr F_i,
\]
which we call the \emph{analytic model frames}, since they are needed to define our function spaces.

Fix an integer $k\geq2$ and a sufficiently small positive weight $\delta$.  See our weighted Sobolev space conventions in \Cref{sec: Weighted Sobolev Spaces}.  We use the $h_0$-unitary gauge group
\begin{equation}
    G=\left\{g\in \mathcal H^{k+1}_{\mathrm{loc}}(C^\circ,\SU(E))\ \Big|\ g-g_i\in \mathcal H_\delta^{k+1}\text{ on }B_1^\circ(p_i)\text{ for some }g_i\in K_i\right\},
\end{equation}
The Lie algebra $\mathfrak g$ consists of sections $\xi$ with limits $\xi_i=\lim_{z\to p_i}\xi(z)\in\mathfrak k_i$ and $\xi-\xi_i\in \mathcal H_\delta^{k+1}$ on each punctured disk $B_1^\circ(p_i)$, so elements can be extended to $C$.

For $a_i\in\mathfrak z(\mathfrak k_i)$, the point distribution $a_i\delta_{p_i}$ given by
\[
    \langle a_i\delta_{p_i},\xi\rangle=\langle a_i,\xi_i\rangle_{h_0}
\]
defines an element of $\mathfrak g^*$.  We use this to identify $\bigoplus_i\mathfrak z(\mathfrak k_i)\subset\mathfrak g^*$.

For $g\in G$, let $g_i^{(k)}$ denote the induced actions of $g_i$ on the associated graded $\gr F_i^{(k)}$.  Let $C_{r_\star}\subset G$ be the subgroup of $r_\star$th roots of unity (as scalar endomorphisms).  Then
\begin{equation}\label{eq: definition of H}
    H=C_{r_\star}\cdot\left\{g\in G\ \Big|\ \det g_i^{(k)}=1\text{ for every }i,k\right\}
\end{equation}
is a subgroup of $G$, and the quotient $B:=G/H$ is the \emph{effective boundary torus} $(T_1\times\cdots\times T_n)/C_{r_\star}$.  The Lie algebra $\mathfrak h$ of $H$ satisfies $\mathrm{Ann}(\mathfrak h)=\bigoplus_i\mathfrak z(\mathfrak k_i)$ under the identification $\bigoplus_i\mathfrak z(\mathfrak k_i)\subseteq\mathfrak g^*$ above.  The short exact sequence
\[
    0\longrightarrow\mathfrak b^*\longrightarrow\mathfrak g^*\longrightarrow\mathfrak h^*\longrightarrow0
\]
therefore shows $\bigoplus_i\mathfrak z(\mathfrak k_i)=\mathfrak b^*\subset\mathfrak g^*$.

\subsection{Distribution-Valued Moment Maps}\label{sec: Distribution-Valued Moment Maps}

We define the extended configuration space of Higgs bundles, where the parabolic weights and Levi-central residues are allowed to vary within an admissible parameter locus to be specified later.  Write $\alpha_i$ and $\sigma_i$ for the block scalar matrices holding the weights and residues, respectively, so that $\I\alpha_i\in\mathfrak z(\mathfrak k_i)$ and $\sigma_i\in\mathfrak z(\mathfrak l_i)$.  Let $\mathscr C$ consist of $\mathcal H^k_{\mathrm{loc}}$ pairs $(A,\Phi)$ with fixed determinant, which can locally be written as
\begin{equation}\label{eq: configuration boundary conditions}
    A=\I\alpha_i\,\de\theta + \mathcal H_\delta^{k},\qquad
    \Phi=-2\I\left(\Re\sigma_i\,\de\theta+\Im\sigma_i\,\frac{\de r}r\right)+\mathcal H_\delta^{k}
\end{equation}
for some admissible $\alpha,\sigma$.  Our Sobolev space convention for 1-forms means each remainder has the form $a\,\de r/r+b\,\de\theta$ for functions $a,b\in\mathcal H_\delta^k(B_\eps^\circ(p_i))$.

We shall interpret Hitchin's equations in terms of distributions on $C$, so that they can be packaged together with weight and residue specifications in a moment map.  For a locally integrable 1-form $\beta$ on $U\subseteq C$, its distributional derivative $\de\beta$ pairs with a smooth compactly supported test section $\xi\in C_c^\infty(U,\mathfrak{su}(E,h_0))$ by
\[
    \langle\de\beta,\xi\rangle
    =\int_U\langle\beta\wedge\de\xi\rangle_{h_0}.
\]
This agrees with the usual derivative away from $D$, but also retains contributions at the punctures.
\begin{definition}\label{def: distributional and boundary moment maps}
    The \emph{distributional moment map} is the triple $\mu=(\mu_1,\mu_2,\mu_3)\from\mathscr C\to(\mathfrak g^*)^3$ defined by
    \begin{equation}\label{eq: extended moment maps}
        \mu_1=F_A-\tfrac12[\Phi\wedge\Phi],\qquad
        \mu_2=-\de_A\Phi,\qquad
        \mu_3=-\de_A(\star\Phi),
    \end{equation}
    interpreted as distributional 2-forms on $C$.  The \emph{boundary moment map} $\mu^\mathrm{bdr}\from\mathscr C\to(\mathfrak b^*)^3$ is the part of $\mu$ supported on $D$.  The remainder $\mu^\mathrm{int}=\mu-\mu^\mathrm{bdr}$ is the \emph{interior moment map}.
\end{definition}

The interior part $\mu^\mathrm{int}$ maps into $\mathcal H_\delta^{k-1}(C^\circ,\Omega_{\hat C}^2\otimes\mathfrak{su}(E,h_0))$, whose elements give bounded functionals on $\mathfrak g$ by the weighted Cauchy--Schwarz inequality \eqref{eq: weighted Cauchy Schwarz}.  The name ``moment map'' is justified below.  First, we describe $\mu^\mathrm{bdr}$ explicitly.

\begin{lemma}\label{lem: boundary moment map}
    For $(A,\Phi)\in\mathscr C$ with weights $\alpha$ and residues $\sigma$, the boundary moment map is
    \begin{equation}
        \mu^\mathrm{bdr}(A,\Phi)=\left(
            2\pi\I\sum_i\alpha_i\delta_{p_i},\,
            4\pi\I\sum_i\Re\sigma_i\delta_{p_i},\,
            4\pi\I\sum_i\Im\sigma_i\delta_{p_i}\right).
    \end{equation}
\end{lemma}
\begin{proof}
    Under distributional differentiation, we have
    \begin{equation}
        \de(\de\theta)=2\pi\delta_{p_i},\qquad
        \de\left(\frac{\de r}{r}\right)=0.
    \end{equation}
    Applying these identities to the local models in \eqref{eq: configuration boundary conditions} gives the part of $\mu$ supported on $D$.  The derivative terms in $\mu_1,\mu_2,\mu_3$ are respectively $\de A$, $-\de\Phi$, and $-\de(\star\Phi)$.  Using $\star(\de r/r)=\de\theta$ and $\star\de\theta=-\de r/r$, the leading-order $\de\theta$ coefficients of $A$, $-\Phi$, and $-\star\Phi$ are respectively $\I\alpha_i$, $2\I\Re\sigma_i$, and $2\I\Im\sigma_i$.  The nonlinear terms $[A\wedge A]$, $[\Phi\wedge\Phi]$, $[A\wedge\star\Phi]$, and $[A\wedge\Phi]$ in Hitchin's equations have size $O(r^{\delta-2})$ and are therefore integrable,\footnote{
        Decompose each term as its limiting model plus a remainder, which respectively belong to $O(r^{-1})$ and $O(r^{\delta-1})$ by definition of $\mathscr C$ and the weighted Sobolev embedding \eqref{eq: weighted Sobolev embedding}.  The Lie bracket of two model terms vanishes by the Levi-central condition.  The remaining Lie bracket terms are $O(r^{\delta-2})$.
    }  with no point-mass contribution.
\end{proof}

In light of the lemma, we define $\nu=(\nu_1,\nu_2,\nu_3)$ where
\begin{equation}\label{eq: nu normalized definition}
    \nu_{1,i}=2\pi\I\alpha_i,\qquad
    \nu_{2,i}=4\pi\I\Re\sigma_i,\qquad
    \nu_{3,i}=4\pi\I\Im\sigma_i,\qquad
    \nu_s=\sum_i\nu_{s,i}\delta_{p_i}.
\end{equation}

\begin{definition}
    We call the distributions $\nu_s\in\mathfrak b^*\subset\mathfrak g^*$ the \emph{boundary parameters} or simply \emph{parameters}.  The \emph{parameter space} is the subset $\Lambda\subset(\mathfrak b^*)^3$ of those $\nu$ whose parabolic weights remain strictly ordered.
\end{definition}

We choose $\mathscr C$ to exclude parabolic weight collisions, so that the boundary moment map is $\mu^\mathrm{bdr}\from\mathscr C\to\Lambda$.  Note this is $G$-invariant by the Levi-central condition, and hence so is $\mu$.

\begin{definition}[Extended parabolic moduli space]\label{def: extended parabolic moduli space}
    Let $\mathcal Z=\{(A,\Phi)\in\mathscr C:\mu^\mathrm{int}(A,\Phi)=0\}$ be the solution space of Hitchin's equations on $C^\circ$.  The \emph{extended parabolic moduli space} is $\mathscr M=\mathcal Z/G$.  Using the subgroup $H\subset G$ from \eqref{eq: definition of H}, define the partial quotient $\mathscr E=\mathcal Z/H$, which comes equipped with a quotient map $q\from\mathscr E\to\mathscr M$ and a residual action of the effective boundary torus $B$.

    On $\mathcal Z$, the total moment map equals its boundary part and induces maps $\pi\from\mathscr M\to\Lambda$ and $\tilde\pi\from\mathscr E\to\Lambda$.  In \Cref{sec: The Analytic Family} we will prove that, on the irreducible loci, the quotients are smooth, $q$ is a principal $B$-bundle, and $\pi$ is a submersion.
\end{definition}
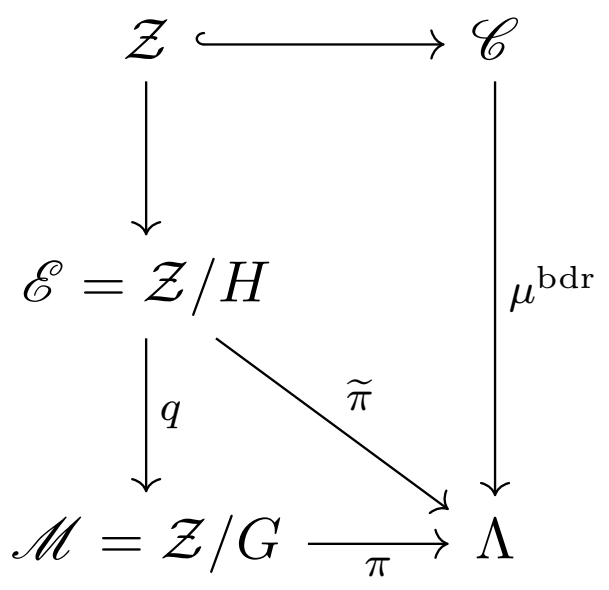
\begin{figure}[!ht]
    \centering
    \begin{tikzcd}[row sep=large]
    \mathcal Z \arrow[r,hook] \arrow[d] & \mathscr C \arrow[dd,"\mu^{\mathrm{bdr}}"] \\
    \mathscr E=\mathcal Z/H \arrow[rd,"\tilde\pi"] \arrow[d,"q"] \\
    \mathscr M=\mathcal Z/G \arrow[r,"\pi"'] & \Lambda
    \end{tikzcd}
    \caption{The partial quotient and the extended parabolic moduli space.}
    \label{fig: quotient diagram}
\end{figure}

This fits the setup for \Cref{thm: DH after partial reduction}.  The next subsection establishes the geometric hypothesis.

\subsection{Closed Forms on the Partial Quotient}\label{sec: The Closed Forms on the Partial Quotient}

We now extend the usual K\"ahler forms on the fixed-parameter moduli spaces to closed forms on the partial quotient $\mathscr E=\mathcal Z/H$, and to prove the Hamiltonian identity \eqref{eq: Hamiltonian identity}.  We in fact construct the extended 2-forms on the full configuration space $\mathscr C$, and show they descend to the desired forms.

If one evaluates the usual $L^2$ pairings on Higgs bundle deformations with varying boundary data, the integral diverges.  We address this by subtracting off a logarithmic term coming from the boundary pairing to get an integrable expression.  On tri-symplectic leaves, the resulting forms agree with the usual $L^2$ pairings.  \Cref{thm: DH after partial reduction} does not require nondegeneracy, and indeed we suspect that our forms degenerate on certain deformations varying boundary data.

Set $C_\eps=C\setminus\bigcup_{i=1}^n B_\eps(p_i)$ for $\eps>0$, and let $u=(\dot A,\dot\Phi)$ and $w=(\dot A',\dot\Phi')$ in $T_{(A,\Phi)}\mathscr C$ be Higgs bundle deformations.  Define the truncated Hitchin forms by
\begin{align}
    \omega_1^\eps(u,w)&=\int_{C_\eps}
    \langle\dot A\wedge\dot A'\rangle_{h_0}
    -\langle\dot\Phi\wedge\dot\Phi'\rangle_{h_0},\label{eq: truncated Hitchin forms}\\
    \omega_2^\eps(u,w)&=\int_{C_\eps}
    \langle\dot A'\wedge\dot\Phi\rangle_{h_0}
    -\langle\dot A\wedge\dot\Phi'\rangle_{h_0},\notag\\
    \omega_3^\eps(u,w)&=\int_{C_\eps}
    \langle\star\dot A\wedge\dot\Phi'\rangle_{h_0}
    -\langle\star\dot A'\wedge\dot\Phi\rangle_{h_0}.\notag
\end{align}
As it stands, the forms $\omega_s^\eps$ may diverge as $\eps\to0$, so we must isolate the divergent parts.

\begin{lemma}\label{lem: leading order of omega_s}
    Let $\dot\nu_{s,i}(u)=\de\mu^\mathrm{bdr}_{s,i}(u)$ and $\dot\nu_{s,i}(w)=\de\mu^\mathrm{bdr}_{s,i}(w)$ be the parameter variation along $u$ and $w$.  For $(s,j,k)$ a cycle of $(1,2,3)$, near $p_i$ the integrand defining $\omega_s^\eps(u,w)$ is
    \begin{equation}\label{eq: leading order of omega_s integrand}
        \frac{1}{4\pi^2}\bigl(\langle\dot \nu_{j,i}(u),\dot \nu_{k,i}(w)\rangle_{h_0}
            -\langle\dot \nu_{j,i}(w),\dot \nu_{k,i}(u)\rangle_{h_0}\bigr)\,r^{-1}\de r\wedge\de\theta + \mathcal H_\delta^k.
    \end{equation}
\end{lemma}
\begin{proof}
Near $p_i$ we write
\[
    \dot A=\dot A_{i,r}\de r+\dot A_{i,\theta}\de\theta,\qquad
    \dot\Phi=\dot\Phi_{i,r}\de r+\dot\Phi_{i,\theta}\de\theta.
\]
To highlight the hyperk\"ahler symmetry between connection and Higgs field, we set
\[
    (u_{0,i},u_{1,i},u_{2,i},u_{3,i})=(-r\dot A_{i,r},\dot A_{i,\theta},-\dot\Phi_{i,\theta},-r\dot\Phi_{i,r})
\]
and similarly for $w$.  Then for $(s,j,k)$ a cycle of $(1,2,3)$, the integrand of $\omega_s^\eps$ near $p_i$ is
\begin{equation}\label{eq: truncated Hitchin forms symmetric}
    \bigl(\langle u_{s,i},w_{0,i}\rangle_{h_0}-\langle w_{s,i},u_{0,i}\rangle_{h_0}
    +\langle u_{j,i},w_{k,i}\rangle_{h_0}-\langle w_{j,i},u_{k,i}\rangle_{h_0}\bigr)\,r^{-1}\de r\wedge\de\theta,
\end{equation}

By definition of $\mathscr C$, the remainders $u_{s,i}-\frac{\dot\nu_{s,i}(u)}{2\pi}$ and $w_{s,i}-\frac{\dot\nu_{s,i}(w)}{2\pi}$ belong to $\mathcal H_\delta^k$ for $s=1,2,3$, and $u_{0,i},w_{0,i}\in\mathcal H_\delta^k$ as well.  Thus the expression in \eqref{eq: leading order of omega_s integrand} indeed captures the leading order of the integrand.  Multiplying these remainders against the constant matrices $\dot\nu_{s,i}(u)$, $\dot\nu_{s,i}(w)$ preserves their regularity, and the Sobolev multiplication estimate \eqref{eq: weighted Sobolev multiplication} controls the products of the remainders with each other.
\end{proof}

Integrating \eqref{eq: leading order of omega_s integrand} over the annuli $\eps\leq r\leq1$ is gives $-\log(\eps)b_s(u,w)$, where
\begin{equation}\label{eq: boundary counterterms}
    b_s(u,w)=\frac{1}{2\pi}\sum_i\bigl(
    \langle\dot \nu_{j,i}(u),\dot \nu_{k,i}(w)\rangle_{h_0}
    -\langle\dot \nu_{j,i}(w),\dot \nu_{k,i}(u)\rangle_{h_0}\bigr).
\end{equation}
Subtracting this term leaves an integrand belonging to $\mathcal H_\delta^{k}$.

\begin{definition}\label{def: extended Hitchin forms}
    For $s=1,2,3$, the \emph{extended Hitchin form} is
    \begin{equation}\label{eq: extended Hitchin forms}
        \tilde\omega_s=\lim_{\eps\to0}\bigl(\omega_s^\eps+\log(\eps)b_s\bigr).
    \end{equation}
\end{definition}

By weighted Cauchy--Schwarz \eqref{eq: weighted Cauchy Schwarz}, we have
\begin{equation}\label{eq: renormalized form convergence}
    \bigl|\bigl[(\omega_s^{\eps}+\log(\eps)b_s)-(\omega_s^{\eps'}+\log(\eps')b_s)\bigr](u,w)\bigr|
    \leq C_\delta\eps^\delta
\end{equation}
for $C$ uniform in $(A,\Phi)$, $0<\eps'\leq\eps$, and $u,w$ in bounded subsets of the tangent space.  This proves convergence.  The forms $\omega_s^\eps$ and $b_s$ are translation-invariant, so their limit $\tilde\omega_s$ is smooth and closed.  The resulting forms do not depend on the analytic frames, but they do depend on the local coordinates $(r,\theta)$ near the punctures.

\bigskip

It remains to check that these forms descend to the quotient $\mathcal Z/H$ and that the residual torus action is Hamiltonian.  We do this by showing that the distributional expressions \eqref{eq: extended moment maps} are moment maps for the extended Hitchin forms.

For $\xi\in\mathfrak g$, the infinitesimal action on $(A,\Phi)$ is $\zeta_\xi=(-\de_A\xi,[\xi,\Phi])$.  We can now state the Hamiltonian identity:

\begin{lemma}\label{lem: Hamiltonian identity}
On the configuration space $\mathscr C$, the extended moment maps satisfy
\begin{equation}\label{eq: parabolic Hamiltonian identity}
    \iota_{\zeta_\xi}\tilde\omega_s=\de\langle\mu_s,\xi\rangle.
\end{equation}
\end{lemma}
\begin{proof}
    It suffices to prove the identity for the limit of the truncated forms $\omega_s^\eps$, as $\zeta_\xi=(-\de_A\xi,[\xi,\Phi])$ has no boundary variation and therefore $b_s(\zeta_\xi,w)$ from \eqref{eq: boundary counterterms} vanishes.  Let $w=(\dot A,\dot\Phi)\in T_{(A,\Phi)}\mathscr C$.  The linearizations on $C^\circ$ are
    \begin{align}
            D\mu_1(w)&=\de_A\dot A-[\Phi\wedge\dot\Phi],\\
            D\mu_2(w)&=-\de_A\dot\Phi-[\dot A\wedge\Phi],\notag\\
            D\mu_3(w)&=-\de_A(\star\dot\Phi)-[\dot A\wedge\star\Phi].\notag
    \end{align}
    Substituting $\zeta_\xi$ into \eqref{eq: truncated Hitchin forms} and integrating by parts gives
    \begin{align}
        \omega_1^\eps(\zeta_\xi,w)
        &=\int_{C_\eps}\langle\xi,D\mu_1(w)|_{C^\circ}\rangle_{h_0}
            -\int_{\partial C_\eps}\langle\xi,\dot A\rangle_{h_0},\label{eq: omegas contracted}
            \\
        \omega_2^\eps(\zeta_\xi,w)
        &=\int_{C_\eps}\langle\xi,D\mu_2(w)|_{C^\circ}\rangle_{h_0}
            +\int_{\partial C_\eps}\langle\xi,\dot\Phi\rangle_{h_0},
            \notag\\
        \omega_3^\eps(\zeta_\xi,w)
        &=\int_{C_\eps}\langle\xi,D\mu_3(w)|_{C^\circ}\rangle_{h_0}
            +\int_{\partial C_\eps}\langle\xi,\star\dot\Phi\rangle_{h_0}.
            \notag
    \end{align}
    Near $p_i$, write
    \begin{equation*}
        u_{1,i}=\dot A_{i,\theta},\qquad
        u_{2,i}=-\dot\Phi_{i,\theta},\qquad
        u_{3,i}=-r\dot\Phi_{i,r}
    \end{equation*}
    as in \eqref{eq: truncated Hitchin forms symmetric}.  The integrals over $C_\eps$ in \eqref{eq: omegas contracted} converge absolutely as $\eps\to0$, matching their respective $D\mu_s^\mathrm{int}(w)$.  Since the inner boundary is oriented by $-\de\theta$, the three boundary integrals become
    \begin{equation}
        \sum_i\int_0^{2\pi}\langle\xi(\eps,\theta),u_{s,i}(\eps,\theta)\rangle_{h_0}\,\de\theta
    \end{equation}
    for the respective $s=1,2,3$.  These integrals read off the leading-order coefficients $\sum_i\langle\xi_i,\dot\nu_{s,i}(w)\rangle_{h_0}$ with an $O(\eps^\delta)$ correction, giving $\langle D\mu_s^\mathrm{bdr}(w),\xi\rangle$.  This proves the identity.
\end{proof}

Since $\mu^\mathrm{int}=0$ gives Hitchin's equations and $\mu^\mathrm{bd}=\nu$ fixes the boundary parameters, Hamiltonian reduction for the extended moment map recovers the fixed-parameter moduli space
\begin{equation}
    \cM_{\nu}
    =\bigl(\mu_1^{-1}(\nu_1)
    \cap\mu_2^{-1}(\nu_2)
    \cap\mu_3^{-1}(\nu_3)\bigr)/G.
\end{equation}
As usual, the combination $\mu_2+\I\mu_3=0$ gives the holomorphicity condition $2\delbar_A\varphi=0$, which is now packaged together with the Levi-central specification for the residues of $\Phi$.  The remaining $\mu_1=0$ gives $F_A+[\Phi,\Phi]=0$ together with the parabolic weight specification.

The forms $\tilde\omega_s$ are $G$-invariant, and the Hamiltonian identity of \Cref{lem: Hamiltonian identity} implies they descend to closed forms $\varpi_1,\varpi_2,\varpi_3$ on $\mathscr E=\mathcal Z/H$.  Moreover, $\mu^\mathrm{bdr}$ descends to a triple $\tilde\pi$ of moment maps for the torus action $B\acts\mathscr E$.  Upon restriction to a fiber $\tilde\pi^{-1}(\nu)$, the forms $\varpi_s$ agree with the usual $L^2$ Hitchin forms.  This proves the geometric hypothesis of \Cref{thm: DH after partial reduction}.

\section{Smooth Structure of the Moduli Family}\label{sec: The Analytic Family}

In this section, we rigorously construct the quotients in \Cref{sec: Universal Moduli Space} and show that their irreducible loci are smooth.  We also show $q\from\mathscr E^{\mathrm{irr}}\to\mathscr M^{\mathrm{irr}}$ is a principal $B$-bundle, and $\pi\from\mathscr M^\mathrm{irr}\to\Lambda$ is a submersion.  We then relate our analytic construction to the algebraic moduli family using harmonic metrics.

Throughout this section we restrict to Higgs bundles whose parabolic weights are strictly decreasing, so that
\begin{equation}\label{eq: ordered weights}
    0<\alpha_i^{(k)}-\alpha_i^{(k')}<1
    \quad\text{for }k<k'.
\end{equation}
We use $\nu=(\nu_1,\nu_2,\nu_3)$ to encode the parabolic weights and residue parameters $(\alpha,\sigma)$ as in \eqref{eq: nu normalized definition}, and we use $\Lambda$ for the space of parameters.

\subsection{Extended Deformation Complex}\label{sec: Extended Deformation Complex}

We first define the analytic configuration space and restate the gauge group's definition, and show the action is well-defined.  We consider the extended deformation complex describing the tangent space of $\mathscr M$, and create a splitting using Biquard--Boalch's inverse to the gauge Laplacian \cite{BB04}.  We then construct local slices with smooth coordinates near irreducible Higgs bundles in $\mathscr M$.  We use the local slices to establish smoothness for $\mathscr M^\mathrm{irr}$, to show that $q\from\mathscr E^{\mathrm{irr}}\to\mathscr M^{\mathrm{irr}}$ is a principal $B$-bundle, and to show the parameter map $\pi\from\mathscr M^\mathrm{irr}\to\Lambda$ is a submersion.  Finally, we use direct limits to allow $\delta\to0$.

Fix the $\SL(r_\star,\C)$-bundle $E\to C$ with its parabolic filtrations, background metric $h_0$ and analytic model frames, with $h_0$ normalized so that the trivialization $\det E\cong\mathcal O_C$ has unit norm in $h_0$.  The configuration space and gauge action below are defined using $\mathcal H_\delta^k$ Sobolev completions for any given $k\geq2$ and $\delta>0$, though we will eventually specify an upper bound on $\delta$ for later results.

Let $\mathcal A_{\mathrm{loc}}^k(E,h_0)$ be the affine space of $\mathcal H^k_{\mathrm{loc}}$ unitary connections inducing the fixed connection on $\det E$.  We define the configuration space
\[
    \mathscr C\subseteq\mathcal A_{\mathrm{loc}}^k(E,h_0)
    \times \mathcal H^k_{\mathrm{loc}}(C^\circ,\Omega_{\hat C}^1\otimes\mathfrak{su}(E,h_0))
\]
consisting of pairs $(A,\Phi)$ such that, on each $B_1^\circ(p_i)$,
\[
    A-\frac{\nu_{1,i}}{2\pi}\de\theta\in \mathcal H_\delta^k(B_1^\circ(p_i),\,\Omega_{\hat C}^1\otimes\mathfrak{su}(E,h_0))
\]
\[
    \Phi+\frac{1}{2\pi}\left(\nu_{3,i}\frac{\de r}{r}+\nu_{2,i}\de\theta\right)\in \mathcal H_\delta^k(B_1^\circ(p_i),\,\Omega_{\hat C}^1\otimes\mathfrak{su}(E,h_0))
\]
for some parameter triple $\nu$ in the admissible locus $\Lambda$.

As before, we define the gauge group
\[
    G=\left\{g\in \mathcal H^{k+1}_{\mathrm{loc}}(C^\circ,\SU(E,h_0))\ \Big|\
    g-g_i\in \mathcal H_\delta^{k+1}\text{ on }B_1^\circ(p_i)\text{ for some }
    g_i\in K_i\right\}
\]
where $K_i$ is the compact Levi group defined in \Cref{sec: Universal Moduli Space}.  Multiplication is smooth by the weighted Sobolev multiplication estimate \eqref{eq: weighted Sobolev multiplication}.  Since $g$ is unitary, $g^{-1}=g^{\dagger}$, and $g^{-1}-g_i^{-1}=(g-g_i)^{\dagger}\in \mathcal H_\delta^{k+1}$, so inversion is smooth.  The Lie algebra of $G$ is
\[
    \mathfrak g=\left\{\xi\in \mathcal H^{k+1}_{\mathrm{loc}}(C^\circ,\mathfrak{su}(E,h_0))\ \Big|\
    \xi-\xi_i\in \mathcal H_\delta^{k+1}\text{ on }B_1^\circ(p_i)
    \text{ for some }\xi_i\in\mathfrak k_i\right\}.
\]
The \emph{asymptotic model} of $(A,\Phi)$ near $p_i$ is the leading order part
\begin{equation}\label{eq: asymptotic model}
    A_{\mathrm{mod}}=\frac{\nu_{1,i}}{2\pi}\de\theta,\qquad
    \Phi_{\mathrm{mod}}=-\frac{1}{2\pi}
        \left(\nu_{3,i}\frac{\de r}{r}+\nu_{2,i}\de\theta\right),
\end{equation}
and the asymptotic model of $g$ and $\xi$ are the limits $g_i\in K_i$ and $\xi_i\in\mathfrak k_i$.

Near $p_i$, we use the Sobolev norm \eqref{eq: cylindrical Sobolev norm} on functions.  We apply this to sections of $E$ using $h_0$ and the model frames.  For differential forms, we work in the $b$-cotangent frame $\de r/r,\de\theta$.  For compactly supported deformations $u=(\dot A,\dot\Phi)$, $w=(\dot A',\dot\Phi')$, the pairing is
\[
    (u,w)=-\int_C\langle\star\dot A\wedge\dot A'\rangle_{h_0}
        +\langle\star\dot\Phi\wedge\dot\Phi'\rangle_{h_0}.
\]
The Hodge star uses the metric $g_C$, which, in a holomorphic coordinate $z=re^{\I\theta}$, has the form
\[
    g_C=\kappa\,\de\bar z\,\de z
    =\kappa(\de r^2+r^2\de\theta^2),\qquad
    \de\mathrm{vol}=\kappa r\,\de r\,\de\theta.
\]
See \eqref{eq: weighted Hodge star} for Hodge star formulas.  These pairings extend from the spaces of compactly supported sections to the Sobolev completions, and are used to define formal adjoints below.

\begin{lemma}\label{lem: gauge action}
The gauge group $G$ acts smoothly on $\mathscr C$, preserving the weight and residue parameters.
\end{lemma}
\begin{proof}
Away from the punctures, this follows from ordinary Sobolev multiplication.  On $B_1^\circ(p_i)$, write $A=A_{\mathrm{mod}}+a$, $\Phi=\Phi_{\mathrm{mod}}+\psi$, and $g=g_i+u$, where $a,\psi\in\mathcal H_\delta^k$ and $u\in\mathcal H_\delta^{k+1}$.  Since $g_i$ is constant and commutes with $A_\mathrm{mod}$ and $\Phi_\mathrm{mod}$,
\[
    \begin{aligned}
        g\cdot A-A_{\mathrm{mod}}
        &=gag^{-1}+\bigl([u,A_{\mathrm{mod}}]-\de u\bigr)g^{-1},\\
        g\cdot\Phi-\Phi_{\mathrm{mod}}
        &=g\psi g^{-1}+[u,\Phi_{\mathrm{mod}}]g^{-1}.
    \end{aligned}
\]
The Lie brackets involving $u$ belong to $\mathcal H_\delta^{k+1}$, and $r\partial_ru,\partial_\theta u\in \mathcal H_\delta^k$.  The weighted multiplication estimate \eqref{eq: weighted Sobolev multiplication} shows products of the remainders lie in $\mathcal H_\delta^k$, with smooth dependence on $g,a,\psi,\nu$.  Since $g_i$ is Levi-central, the action preserves the boundary data and decay conditions.
\end{proof}

\bigskip

Let $\mathscr C_{\nu_*}=(\mu^\mathrm{bdr})^{-1}(\nu_*)$ be the \emph{fixed parameter configuration space}.  The tangent space is
\begin{equation}\label{eq: variation space}
    T_{x}\mathscr C_{\nu_*}=\ker\de\mu_{x}^\mathrm{bdr}=\mathcal H_\delta^k\bigl(C^\circ,
    (\Omega_{\hat C}^1\oplus\Omega_{\hat C}^1)\otimes\mathfrak{su}(E,h_0)\bigr).
\end{equation}

Fix an irreducible solution $x=(A_*,\Phi_*)$ with $\mu(x)=\nu_*$.  The usual fixed-parameter deformation complex is
\begin{equation}\label{eq: classic deformation complex}
    0\longrightarrow \mathfrak g\xrightarrow{d_0}T_{x}\mathscr C_{\nu_*}
    \xrightarrow{d_1}\mathcal H_\delta^{k-1}(C^\circ,\Omega_{\hat C}^2\otimes\mathfrak{su}(E,h_0))^{\oplus3}\longrightarrow0,
\end{equation}
with $d_0\xi=(-\de_{A_*}\xi,[\xi,\Phi_*])$ the infinitesimal action and $d_1$ the restriction of $\de\mu^{\mathrm{int}}_{x}$.  The quotient $\ker d_1/\operatorname{im}d_0$ is the tangent space to the fixed-parameter moduli space $\mathscr M_{\nu_*}$.  The adjoint $d_0^*$ satisfies
\[
    (d_0\xi,u)=(\xi,d_0^*u),
\]
and maps $T_{x}\mathscr C_{\nu_*}$ to $\mathcal H_{\delta-2}^{k-1}(C^\circ,\mathfrak{su}(E,h_0))$.  The following lemma is due to Biquard and Boalch:

\begin{lemma}\label{lem: gauge laplacian}
    Fix an admissible parameter $\nu_*\in\Lambda$ with parabolic weights $\alpha_i^{(k)}$.  Choose $\delta>0$ satisfying
    \begin{equation}\label{eq: Sobolev weight bound}
        \delta<\min_{i,k}\left(\alpha_i^{(k)}-\alpha_i^{(k+1)}\right),\qquad
        \delta<\min_i\left(1-\alpha_i^{(0)}+\alpha_i^{(\ell_i)}\right).
    \end{equation}
    Then the gauge Laplacian
    \begin{equation}\label{eq: gauge Laplacian}
        \mathcal L_{(A_*,\Phi_*)}=d_0^*d_0\from \mathfrak g\longrightarrow \mathcal H_{\delta-2}^{k-1}(C^\circ,\mathfrak{su}(E,h_0))
    \end{equation}
    is an isomorphism at any irreducible $(A_*,\Phi_*)\in\mu^{-1}(\nu_*)$.
\end{lemma}
\begin{proof}
    In the tame case, \cite[Lemma~5.1]{BB04} is our situation.  Their gauge Lie algebra is equal to ours, with equivalent norms for sufficiently small $\delta>0$.  The gauge Laplacian has the same domain and codomain.\footnote{
        The off-diagonal norms are equivalent to $\mathcal H_\delta^2$, and on Levi blocks \cite[Equation~(3.4)]{BB04} separates the constant limit.  Weighted elliptic regularity upgrades the domain and codomain orders to $k+1$ and $k-1$, respectively.
    }
\end{proof}

\begin{remark}\label{rem: parabolic gaps}
We emphasize that the condition on $\delta$ depends on $\nu_*$.  The quantities $\alpha_i^{(k)}-\alpha_i^{(k+1)}$ and $1-\alpha_i^{(0)}+\alpha_i^{(\ell_i)}$ in \eqref{eq: Sobolev weight bound} are the gaps between parabolic weights in $\R/\Z$, and the condition \eqref{eq: ordered weights} defining $\Lambda$ says these gaps are positive.
\end{remark}

The inverse of the gauge Laplacian allows one to project deformations in $T_{x}\mathscr C_{\nu_*}$ to the orthogonal complement of gauge orbits, placing them in \emph{Coulomb gauge}.  Biquard--Boalch also provide a useful right inverse for $d_1$:
\begin{lemma}\label{lem: moment map right inverse}
Let
\[
    S_{x}=\ker\bigl(d_0^*\colon T_{x}\mathscr C_{\nu_*}\longrightarrow
        \mathcal H_{\delta-2}^{k-1}(C^\circ,\mathfrak{su}(E,h_0))\bigr)
\]
be the fixed-parameter Coulomb space at $x$.  Under the hypotheses of \Cref{lem: gauge laplacian}, the map
\[
    d_1|_{S_{x}}\from S_{x}\longrightarrow \mathcal H_\delta^{k-1}(C^\circ,\Omega_{\hat C}^2\otimes\mathfrak{su}(E,h_0))^{\oplus3}
\]
admits a bounded right inverse $\mathcal G$.
\end{lemma}

Now we define the extended deformation complex for the universal configuration space.  We begin by enlarging the codomain of $d_1$ to include boundary variations:
\begin{equation}\label{eq: image de mu definition}
    \mathscr I_x
        =\mathcal H_\delta^{k-1}\bigl(C^\circ,\Omega_{\hat C}^2\otimes\mathfrak{su}(E,h_0)\bigr)^{\oplus3}+T_{\nu_*}\Lambda
            \subset(\mathfrak g^*)^{\oplus3}.
\end{equation}
As before, the second summand consists of distributions supported on $D$, so that pairing $\mathscr I_x\times\mathfrak g^{\oplus3}\to\C$ is defined using integration over $C$.  We equip this space with the Banach norm
\[
    \|\eta+\nu\|_{\mathscr I_x}
    =\|\eta\|_{\mathcal H_\delta^{k-1}}+\|\nu\|,
\]
using any fixed norm on $(\mathfrak b^*)^{\oplus3}$.

Going forward, $d_0$ shall refer to the extended deformation complex
\begin{equation}
    0\longrightarrow\mathfrak g\xrightarrow{d_0}T_{x}\mathscr C
    \xrightarrow{\de\mu_{x}}\mathscr I_x
    \longrightarrow0
\end{equation}
where $\de\mu=\de\mu^{\mathrm{int}}+\de\mu^{\mathrm{bdr}}$ is the linearization of the extended moment map.  Surjectivity is shown below.  We define the extended Coulomb space by
\begin{equation}\label{eq: extended Coulomb space}
    \tilde S_{x}=\ker d_0^*
        =\{u\in T_{x}\mathscr C:(d_0\xi,u)=0
            \text{ for every }\xi\in\mathfrak g\}.
\end{equation}
Although $u$ need not have finite $L^2$ norm, its pairing with $d_0\xi$ converges because of the decay condition on $\xi$ on the off-diagonal blocks in the model frame.  Under the hypotheses of \Cref{lem: gauge laplacian}, the orthogonal projection
\begin{equation}\label{eq: proj to extended Coulomb space}
    \Pi_{x}=\Id-d_0\mathcal L_{x}^{-1}d_0^*\from T_{x}\mathscr C\to\tilde S
\end{equation}
gives the extended Coulomb splitting
\begin{equation}\label{eq: extended Coulomb splitting}
    T_{x}\mathscr C=\operatorname{im}d_0\oplus\tilde S_{x},
    \qquad
    \tilde S_{x}\cap T_{x}\mathscr C_{\nu_*}=S_{x}.
\end{equation}

\begin{proposition}\label{prop: extended moment map right inverse}
    Under the hypotheses of \Cref{lem: gauge laplacian}, the restriction of the linearized moment map to the extended Coulomb space
    \[
        \de\mu_{x}\from\tilde S_{x}\longrightarrow\mathscr I_x
    \]
    admits a bounded right inverse $\tilde{\mathcal G}$.  In particular,
    $\operatorname{im}\de\mu_{x}=\mathscr I_x$.
\end{proposition}
\begin{proof}
    Given a vector in $\dot\nu\in T_{\nu_*}\Lambda$, one can construct a lift $(\dot A,\dot\Phi)\in T_{x}\mathscr C$ using bump functions and the asymptotic model \eqref{eq: asymptotic model}.  Fixing the bump functions defines a continuous right inverse $m\from T_{\nu_*}\Lambda\longrightarrow T_{x}\mathscr C$ for the unrestricted $\de\mu_{x}$.  Composing with the projection $\Pi_{x}$ from \eqref{eq: proj to extended Coulomb space} gives
    \[
        \tilde m=\Pi_{x}\circ m\from T_{\nu_*}\Lambda\longrightarrow\tilde S_{x}.
    \]
    We still have $\de\mu^{\mathrm{bdr}}_{x}\circ\tilde m=\Id$, since the term $d_0\mathcal L_{x}^{-1}d_0^*$ in \eqref{eq: proj to extended Coulomb space} is killed by $\de\mu_{x}^\mathrm{bdr}$.

    Let $f+\dot{\nu}\in\mathscr I_x$, where the summands respect the decomposition \eqref{eq: image de mu definition}.  The lift $\tilde m(\dot\nu)$ captures the boundary variation, but need not infinitesimally preserve Hitchin's equations.  However, $\tilde m(\dot\nu)-\mathcal G\circ\de\mu_{x}^{\mathrm{int}}(\tilde m(\dot\nu))$
    does:
    \begin{equation}\label{eq: adjusted lift preserves Hitchin equations}
        \de\mu_{x_0}^\mathrm{int}\left(\tilde m(\dot\nu)-\mathcal G\circ\de\mu_{x}^{\mathrm{int}}(\tilde m(\dot\nu))\right)=\de\mu_{x_0}^\mathrm{int}\left(\tilde m(\dot\nu)\right)-\de\mu_{x_0}^\mathrm{int}\left(\tilde m(\dot\nu)\right)=0.
    \end{equation}
    Here we use $\mathcal G$ from \Cref{lem: moment map right inverse}, recalling that $d_1$ is the restriction of $\de\mu^\mathrm{int}$.  Define
    \begin{equation}\label{eq: def of G tilde}
        \tilde{\mathcal G}(f+\dot{\nu})
            =\mathcal G(f)+\tilde m(\dot\nu)-\mathcal G\circ\de\mu_{x}^{\mathrm{int}}(\tilde m(\dot\nu))\in\tilde S_{x}.
    \end{equation}
    Since $\mathcal G(f)$ is tangent to $\mathscr C_{\nu_*}=(\mu^\mathrm{bdr})^{-1}(\nu_*)$, we have
    \begin{align}
        \de\mu_x\circ\tilde{\mathcal G}(f+\dot\nu)
            &=(\de\mu^\mathrm{int}+\de\mu^\mathrm{bdr})(\mathcal G(f)+\tilde m(\dot\nu)-\mathcal G\circ\de\mu_{x}^{\mathrm{int}}(\tilde m(\dot\nu)))
                \notag\\
            &=\de\mu^\mathrm{int}\circ\mathcal G(f)+\de\mu^\mathrm{bdr}(\tilde m(\dot\nu)-\mathcal G\circ\de\mu_{x}^{\mathrm{int}}(\tilde m(\dot\nu)))
                \notag\\
            &=f+\dot\nu,
    \end{align}
    where the second equality used \eqref{eq: adjusted lift preserves Hitchin equations}, and the third used $\de\mu_x^\mathrm{int}\circ\mathcal G=\Id$.
\end{proof}

\subsection{Local Slices}\label{sec: Local Slices}

When working with the universal moduli space $\mathscr M$, care must be taken to handle the dependence of $\delta$ on $\nu\in\Lambda$.  To address this, we fix $0<\delta\ll1$ and define $\Lambda_\delta\subset\Lambda$ consisting of those $\nu$ whose parabolic weight gaps in $\R/\Z$ are all greater than $\delta$ (see \Cref{rem: parabolic gaps}). Let $\mathscr C_\delta^\mathrm{irr}\subset(\mu^\mathrm{bdr})^{-1}(\Lambda_\delta)$ be the configuration space of irreducible Higgs bundles with parameters in $\Lambda_\delta$.  Recall the subgroup $H\subset G$ defined in \eqref{eq: definition of H} and the effective boundary torus $B=G/H$.  Set $\mathscr E_\delta^\mathrm{irr}=\mathscr C_\delta^\mathrm{irr}/H$ and $\mathscr M_\delta^\mathrm{irr}=\mathscr C_\delta^\mathrm{irr}/G$.

\begin{proposition}\label{prop: local regularity}
Fix $0<\delta\ll1$ and let $\Lambda_\delta$ be as above.  The extended moduli space of irreducible Higgs bundles $\mathscr M_\delta^{\mathrm{irr}}$ has a smooth structure, and the moment map descends to a submersion $\pi\from\mathscr M_\delta^{\mathrm{irr}}\to\Lambda_\delta$.  For the partial quotient $\mathscr E_\delta=\mu^{-1}(\Lambda_\delta)/H$, the restriction $q\from\mathscr E_\delta^{\mathrm{irr}}\to\mathscr M_\delta^{\mathrm{irr}}$ is a principal $B$-bundle.
\end{proposition}
\begin{proof}
Let $\mathcal Z_\delta=(\mu^\mathrm{int})^{-1}(0)$ be the solution space to Hitchin's equations, and let $x_*=(A_*,\Phi_*)\in\mathcal Z_\delta$ be an irreducible solution satisfying $\mu(x_*)=\nu_*\in\Lambda_\delta$.  Consider the affine Coulomb slice $\mathcal C=x_*+\tilde S_{x_*}\subset\mathscr C_\delta$, where $\tilde S_{x_*}=\ker\de_0^*\subset T_{x_*}\mathscr C$ is the extended Coulomb space from \eqref{eq: extended Coulomb space}.  At the point $x_*$, the slice $\mathcal C$ is orthogonal to the $G$-orbit (see \Cref{fig: local Coulomb slice}).

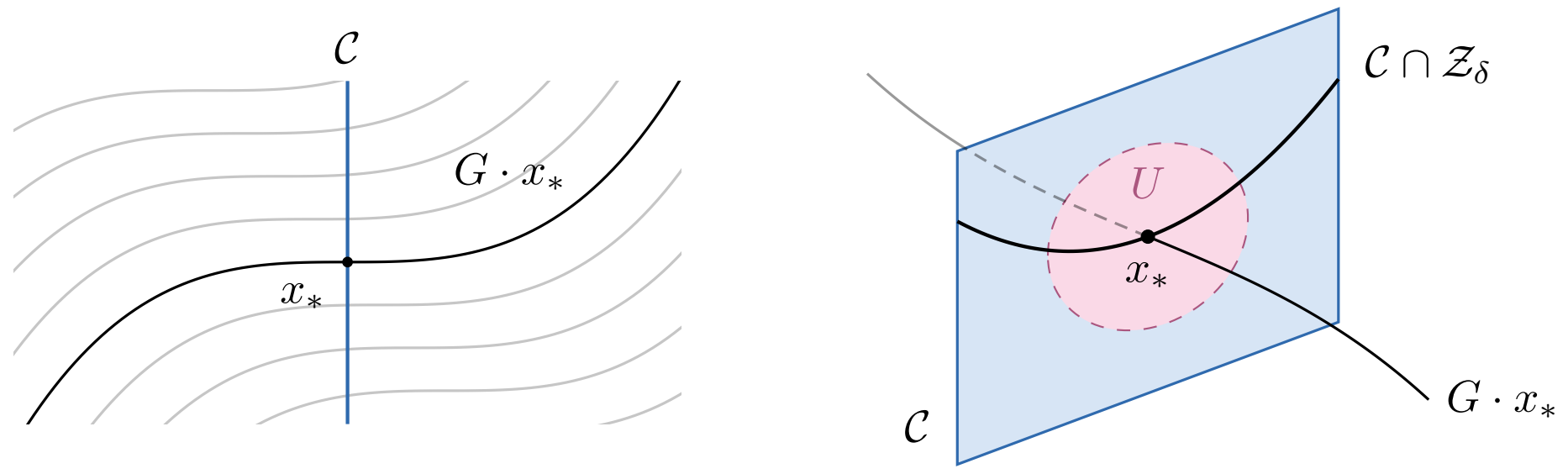
\begin{figure}[tbp]
    \centering
    \begin{tikzpicture}[x=1cm,y=1cm,font=\small]
        \definecolor{transportBlue}{RGB}{55,125,205}
        \begin{scope}[scale=0.75]
        \begin{scope}
            \clip (-3.5,-1.7) rectangle (3.5,1.9);
            \foreach \s in {-6,-5,-4,-3,-2,-1,1,2,3,4,5,6} {
                \draw[semithick,gray!45,domain=-6:6,samples=160]
                    plot ({\x-0.35*\s},{0.045*\x*\x*\x+0.45*\s});
            }
            \draw[semithick,domain=-3.5:3.5,samples=100]
                plot (\x,{0.045*\x*\x*\x});
        \end{scope}
        \draw[thick,transportBlue!85!black] (0,-1.7) -- (0,1.9) node[above,text=black] {$\mathcal C$};
        \node[above right,xshift=4mm,yshift=2mm] at (0.4,0.3) {$G\cdot x_*$};
        \fill (0,0) circle (1.6pt) node[below left=1pt] {$x_*$};
        \end{scope}
        \begin{scope}[shift={(6.3,0.2)},x={(-0.75cm,-0.28cm)},y={(0cm,0.85cm)},z={(-0.85cm,0.35cm)}]
            \filldraw[fill=transportBlue!20,draw=transportBlue!85!black,semithick]
                (-2,-1.45) -- (2,-1.45) -- (2,1.45) -- (-2,1.45) -- cycle;
            \node[below left,yshift=2mm] at (2.1,-1) {$\mathcal C$};
            \filldraw[fill=magenta!18,draw=magenta!65!black,dashed,domain=0:360,samples=100]
                plot ({1.05*cos(\x)},{0.85*sin(\x)+0.17*cos(\x)}) -- cycle;
            \node[text=magenta!65!black] at (0,0.5) {$U$};
            \draw[thick,domain=-2:2,samples=80]
                plot (\x,{0.20*\x*\x});
            \node[right] at (-2.1,0.9) {$\mathcal C\cap\mathcal Z_\delta$};
            \draw[semithick,domain=0:-2.6,samples=60]
                plot (0,{0.025*\x*\x*\x},\x)
                node[right] {$G\cdot x_*$};
            \draw[semithick,opacity=0.4,densely dashed,domain=0:1.692515,samples=60]
                plot (0,{0.025*\x*\x*\x},\x);
            \draw[semithick,opacity=0.4,domain=1.692515:2.6,samples=40]
                plot (0,{0.025*\x*\x*\x},\x);
            \fill (0,0,0) circle (1.6pt) node[below = 2pt] {$x_*$};
        \end{scope}
    \end{tikzpicture}
    \caption{The  affine Coulomb slice $\mathcal C$ is orthogonal to the gauge orbit at $x_*$.  On the right, $\mathcal C$ is shown as a plane containing a neighborhood $U$ of $x_*$ with $\mathcal C\cap\mathcal Z_\delta$ passing through.}
    \label{fig: local Coulomb slice}
\end{figure}

The tangent space of the intersection $\mathcal C\cap\mu^{-1}(\nu_*)$ at $x_*$ is the space of fixed-parameter harmonic deformations
\[
    \mathscr H=\tilde S_{x_*}\cap\ker\de\mu_{x_*}\subset T_{x_*}\mathscr C_{\nu_*},
\]
which is finite dimensional by \cite[Lemma~5.2]{BB04}.  The right inverse $\tilde{\mathcal G}$ of $\de\mu_{x_*}$ from \Cref{prop: extended moment map right inverse} gives the splitting
\begin{equation}
    \tilde S_{x_*}
    =\mathscr H\oplus\tilde{\mathcal G}(\mathscr I_{x_*}).
\end{equation}
The Banach submersion theorem therefore gives local coordinates $(a,b)$ on a neighborhood $U\subset \mathcal C$ of $x_*$, where $a\in\mathscr H$, $b\in\mathscr I_{x_*}$, and $\mu(a,b)=\nu_*+b$.  Restricting to $U\cap\mathcal Z_\delta$, where the interior moment map vanishes, we obtain coordinates $(a,b)$ for $a\in\mathscr H$ and $b\in T_{\nu_*}\Lambda_\delta\subset\mathscr I_{x_*}$ sufficiently small.  These coordinates satisfy $\mu(a,b)=\nu_*+b\in\Lambda_\delta$.

\begin{claim}\label{claim: slice homeomorphism}
    Let $p'$ be the map taking $U\cap\mathcal Z_\delta$ to its image $V\subseteq\mathscr M_\delta$ under the quotient map $p\from\mathcal Z_\delta\to\mathscr M_\delta$.  After shrinking $U$ and $V$ if needed, $p'$ is a homeomorphism.
\end{claim}
\begin{proof}[Proof of \Cref{claim: slice homeomorphism}]
To invert $p'$, we show that every orbit near $[x_*]$ meets the slice $U\cap\mathcal Z_\delta$ in exactly one point.  For simplicity we now work with the effective gauge group $G^\mathrm{eff}=G/C_{r_\star}$ (recall $C_{r_\star}$ is the $r_\star$th roots of unity).  Consider the map
\begin{equation}\label{eq: def of N}
    \mathcal N_x\from\mathfrak g\to \mathcal H_{\delta-2}^{k-1}(C^\circ,\mathfrak{su}(E,h_0)),\qquad\xi\mapsto d_0^*\bigl(e^{\xi}\cdot x-x_*\bigr)
\end{equation}
for $x\in\mathscr C_\delta$.  Since $\mathcal C=x_*+\ker d_0^*$, a solution to $\mathcal N_x(\xi)=0$ produces a gauge transformation $e^\xi$ taking $x$ into $\mathcal C$.  The linearization $\de\mathcal N_{x_*}$ at $\xi=0$ is the gauge Laplacian, which is an isomorphism by \Cref{lem: gauge laplacian}.  The implicit function theorem provides a unique smooth family of gauge corrections $e^{\xi(x)}$ as $x$ varies in a neighborhood $\tilde V\subset\mathscr C_\delta$ of $x_*$.  Shrink the neighborhoods so that $U\subset\tilde V$, $p(V\cap\mathcal Z_\delta)=V$, and $e^{\xi(x)}\cdot x\in U$ for every $x\in\tilde V$.  Since $G$ preserves $\mathcal Z_\delta$, this shows every $x\in\tilde V\cap\mathcal Z_\delta$ has a conjugate point in $U\cap\mathcal Z_\delta$.

For the proposed $x\mapsto e^{\xi(x)}\cdot x$ to descend to an inverse $V\to U\cap\mathcal Z_\delta$ of $p'$, we must prove invariance on group orbits.  Concretely, this can fail if a pair $x,y\in\tilde V$ is equivalent under a gauge transformation far from the identity.  Allowing ourselves to shrink $U$, $\tilde V$, $V$ further, we must show that if $y_n=g_n\cdot x_n$ for sequences $x_n,y_n\to x_*$ and $g_n\in G^\mathrm{eff}$, then $g_n\to1$.  In the fixed-parameter setting, Biquard \cite[Lemma~1.7]{Biq91} obtains the necessary weighted Sobolev bounds to prove $g_n\to1$.\footnote{
    Using the effective gauge group avoids needed to pass to convergent subsequences.
}  To apply this to our case, we note that the parabolic weight gaps are uniformly bounded below by $\delta$, and the leading order Levi parts of $x_n,y_n$ are constrained to a small compact set near $\nu_*$, so the Sobolev bounds of \cite{Biq91} still hold.  This proves the claim.
\end{proof}

The map $p'$ takes the coordinates $(a,b)$ on $U\cap\mathcal Z_\delta$ to coordinates on $V$, and the transition maps are smooth since $\xi(x)$ varies smoothly.  The map $\pi\from\mathscr M_\delta\to\Lambda_\delta$ descends from $\mu$, so we have $\pi\circ h=\mu|_{U\cap\mathcal Z_\delta}$ and $\pi(a,b)=\nu_*+b$.  Thus $\pi$ is a submersion at $x_*$.

To show $q\from\mathscr E_\delta^\mathrm{irr}\to\mathscr M_\delta^\mathrm{irr}$ is a principal $B$-bundle, recall the neighborhood $\tilde V\subset\mathscr C_\delta$ of $x_*$ from the proof of the claim.  The implicit function theorem argument used to produce gauge corrections $\xi(x)$ gives the homeomorphism
\[
    O\times(U\cap\mathcal Z_\delta)\to\tilde V\cap\mathcal Z_\delta,
        \qquad(g,y)\mapsto g\cdot y
\]
where $O\subset G^\mathrm{eff}$ is a neighborhood of the identity.  Since $C_{r_\star}\subset H$, the quotient of $\mathcal O\times(U\cap\mathcal Z_\delta)$ by $H$-equivalence gives
\[
    q^{-1}(V)\simeq B\times V.
\]
The $B$-action is multiplication on the first factor, so $q\from\mathscr E_\delta^\mathrm{irr}\to\mathscr M_\delta^\mathrm{irr}$ is a principal $B$-bundle.
\end{proof}

We now construct the extended moduli family $\mathscr M$ over all of $\Lambda$ by varying $\delta$, with the Sobolev order $k\geq2$ fixed.  We write $\mathscr C_\delta$ for the configuration space and $G_\delta,H_\delta,B_\delta$ for the groups at weight $\delta$.
\begin{proposition}\label{prop: universal moduli space}
    As $\delta\to0$, the families $\mathscr M_\delta$, $\mathscr E_\delta$, and $B_\delta$ form directed systems.  Denote their direct limits
    \begin{equation}\label{eq: universal moduli limit}
        \mathscr M=\lim_{\delta\to0}\mathscr M_\delta,\qquad
        \mathscr E=\lim_{\delta\to0}\mathscr E_\delta,\qquad
        B=\lim_{\delta\to0}B_\delta.
    \end{equation}
    The irreducible locus $\mathscr M^\mathrm{irr}$ has a smooth structure, and $\pi\from\mathscr M^\mathrm{irr}\to\Lambda$ is a submersion.  The preimage $\mathscr E^\mathrm{irr}=q^{-1}(\mathscr M^\mathrm{irr})$ is smooth, and $q\from\mathscr E^\mathrm{irr}\to\mathscr M^\mathrm{irr}$ is a principal $B$-bundle.
\end{proposition}
\begin{proof}
    For $0<\delta'<\delta$, the Sobolev inclusions $\mathcal H_\delta^k\subset \mathcal H_{\delta'}^k$, $\mathcal H_\delta^{k+1}\subset \mathcal H_{\delta'}^{k+1}$ from \eqref{eq: weighted Sobolev inclusion} and $\Lambda_\delta\subset\Lambda_{\delta'}$ give inclusions of the configuration spaces, solution spaces, and gauge groups, inducing maps on the quotients.  Every $B_\delta=G_\delta/H_\delta$ is isomorphic to the same boundary torus $(T_1\times\cdots\times T_n)/C_{r_\star}$.

    To prove injectivity of $\mathscr M_\delta\to\mathscr M_{\delta'}$, we must show that any $g\in G_{\delta'}$ relating two points in $\mathscr C_\delta$ in fact belongs to $G_\delta$.  Denote their connection 1-forms by $A$ and $B=g\cdot A=gAg^{-1}-(\de g)g^{-1}$.  Rearranging to $\de g=gA-Bg$ and comparing the $r^{-1}\de r$ and $\de\theta$ coefficients gives
    \begin{equation}\label{eq: derivatives of g}
        r\partial_r g=gA_r-B_rg,\qquad\partial_\theta g=gA_\theta-B_\theta g.
    \end{equation}
    Since $g_i$ commutes with the common leading-order coefficient $\I\alpha_i$, we have $r\partial_r(g-g_i)\in\mathcal H_\delta^k$ and $\partial_\theta(g-g_i)\in\mathcal H_\delta^k$.  By the fundamental theorem of calculus,
    \begin{equation}
        g(r,\theta)-g_i
        =\int_0^r\bigl(gA_r-B_rg\bigr)(s,\theta)\,\frac{\de s}{s}.
    \end{equation}
    Now $A_r,B_r\in \mathcal H_\delta^k$, and the integration preserves their decay weight by the estimate \eqref{eq: radial integration}, so $g-g_i\in\mathcal H_\delta^k$.  Combining, we get $g-g_i\in\mathcal H_\delta^{k+1}=G_\delta$ as required.  The same argument applies to $H_\delta$.  The compositions $\mathscr M_\delta\xrightarrow{\pi}\Lambda_\delta\into\Lambda$ induce $\pi\from\mathscr M\to\Lambda$, and similarly the maps $\mathscr E_\delta\xrightarrow{q}\mathscr M_\delta\into\mathscr M$ induce $q\from\mathscr E\to\mathscr M$.

    On the irreducible locus, the local charts from \Cref{prop: local regularity} are compatible as the sufficiently small Sobolev weight varies (elliptic regularity identifies the fixed-parameter harmonic deformations), so they define a smooth structure on $\mathscr M^\mathrm{irr}$ compatible with each inclusion $\mathscr M_\delta\into\mathscr M$.  In particular, $\pi\from\mathscr M^\mathrm{irr}\to\Lambda$ is a submersion, and $\mathscr E^\mathrm{irr}\xrightarrow{q}\mathscr M^\mathrm{irr}$ is a smooth principal $B$-bundle.
\end{proof}

\subsection{The Algebraic Moduli Family}\label{sec: The Algebraic Moduli Family}

For the homology transport of \Cref{sec: Homology Transport} and the characteristic class computation of \Cref{sec: affine fiber alternative}, use will the algebraic description of parabolic Higgs bundles to study the family $\mathscr M\to\Lambda$.  The existence of harmonic metrics adapted to the parabolic structure gives the translation from the algebraic construction to the analytic construction of the Higgs bundle moduli spaces when the boundary data is held fixed.  We must show that the family of harmonic metrics---and therefore the conversions to the unitary formulation---depend smoothly on the parameters.  In this subsection, we use $h_\mathrm{bg}$ for the background metric from before to distinguish from a family of metrics $h_t$.

Fix a surface $C$, a divisor $D$, and a degree 0 bundle $E\to C$ with fixed determinant structure.  The algebraic data of a Higgs bundle is a holomorphic structure $\delbar_E$, a holomorphic Higgs field $\varphi\in H^0(C,\End(E)\otimes K_C(D))$, and a parabolic filtration $\mathcal F$ over the marked points in $C$.  The Higgs field residues are required to preserve the filtration and act on $\gr F_i^{(k)}$ by $\sigma_i^{(k)}$.  To define the moduli space $\mathcal M(\alpha,\sigma)$, one fixes the flag dimensions and associated parabolic weights $(\alpha_i^{(k)})_{i,k}$, and considers the space of stable Higgs bundles up to smooth complex gauge transformations in $\SL(E)=C^\infty(C,\SL(r_\star,\C))$.  Unlike the unitary formulation, the flag $\mathcal F$ is part of the varying moduli data and the gauge group can change the flag.  Stability here is given in terms of the parabolic degree of $\varphi$-invariant subbundles, where parabolic degree is the sum of the regular degree and all the induced parabolic weights.

Consider a local family $(\delbar_{E_t},\varphi_t,\mathcal F_t)$, $t\in(-\eps,\eps)$ of stable parabolic Higgs bundles with varying parabolic weights $\alpha(t)$ and residue blocks $\sigma(t)$.  For each $t$, stability implies there exists a \emph{harmonic metric} $h_t$ solving Hitchin's equations and adapted to the parabolic structure (see \cite{Sim90}), which is unique up to unitary equivalence.  To translate to the unitary description with background metric $h_\mathrm{bg}$, let $B_t$ be the matrix satisfying $h_t(v,w)=h_\mathrm{bg}(B_tv,w)$, and let $B_t^{1/2}$ be the unique positive Hermitian root with $(B_t^{1/2})^{\dagger}B_t^{1/2}=B_t$.  Such a $B_t^{1/2}$ is unique, and gives an isometry $(E_t,h_t)\to(E,h_\mathrm{bg})$.  Define $(A_t,\Phi_t)$ by
\begin{equation}\label{eq: harmonic metric identification}
    \de_{A_t}=B_t^{1/2}D_{h_t}B_t^{-1/2},\qquad
    \Phi_t=B_t(\varphi_t^{\dagger_{h_t}}-\varphi_t)B_t^{-1}.
\end{equation}
Since $h_t$ unique up to unitary transformations, this produces a well-defined element of $\mathscr M$.  Indeed, in light of our comparison to Biquard--Boalch's function spaces in \Cref{lem: gauge laplacian}, the weighted correspondence of \cite[Theorem~6.1 and Lemma~7.5]{BB04} gives harmonic representatives in $\mathscr C_\delta$ for sufficiently small $\delta>0$, hence the need for the direct limit in \Cref{prop: universal moduli space}.

To prove smooth dependence of $(A_t,\Phi_t)$ on $t$ near $0$, we give a modification of \cite[Theorem~4.1]{KW18}.  Choose smoothly varying gauge transformations to hold the flags fixed, as in \cite[Section~3]{Yae26}.  In a trivialization respecting the flag, define the model metrics $h_{\mathrm{mod},t}=\bigoplus_k r^{2\alpha_i^{(k)}(t)}\Id$ near $p_i$, which are adapted to the parabolic structure for all $t$.  Choose bump functions $\chi_i$ supported near each puncture, and define $S_t$ by
\[
    S_t=\bigoplus_k r^{\chi_i(r)(\alpha_i^{(k)}(t)-\alpha_i^{(k)}(0))}\Id,
\]
These transformations take the harmonic metric $h_0$ to nearby approximate solutions $h_t^\mathrm{app}$ by setting
\[
    h_t^\mathrm{app}(v,w)=h_0(S_tv,S_tw).
\]
Note $h_0^\mathrm{app}=h_0$, and each $h_t$ is adapted to $\alpha(t)$.  Set $h_{t,u}=h_t^\mathrm{app}e^u$ for a $u\in\I\mathfrak g$.  Shrinking $\eps$ if needed, choose a common sufficiently small $\delta>0$ for all $t\in(-\eps,\eps)$.  Using positive Hermitian roots as before to pick unique isometries $(E,h_{t,u})\to(E,h_\mathrm{bg})$, we define
\begin{equation}
    \mathcal K_t\from\I\mathfrak g\longrightarrow
        \mathcal H_{\delta-2}^{k-1}(C^\circ,\I\mathfrak{su}(E,h_\mathrm{bg})),\qquad
    \mathcal K_t(u)=\I\star\bigl(F_{h_{t,u}}+[\varphi_t^{\dagger_{h_{t,u}}}\wedge\varphi_t]\bigr).
\end{equation}
By weighted Sobolev multiplication \eqref{eq: weighted Sobolev multiplication} and the definition of $\mathfrak g$, this $\mathcal K$ is smooth and satisfies $\mathcal K(t_*,0)=0$.  The linearization of $\mathcal K_0$ at $u=0$ is a multiple of the gauge Laplacian, which is invertible by \Cref{lem: gauge laplacian}.  The implicit function theorem therefore gives a smooth correction $u(t)$; uniqueness with fixed determinant identifies $h_{t,u(t)}$ with $h_t$.

We therefore have continuous bijections between the algebraic moduli family $\mathscr M^\mathrm{alg}$ and the analytic family $\mathscr M$.  On the smooth stable locus, the linearization of this correspondence identifies the tangent spaces computed by the respective deformation complexes, so the correspondence is a local diffeomorphism.

\section{Homology Transport and Period Variation}\label{sec: Homology Transport}

We now establish a homology local system for $\pi\from\mathscr M\to\Lambda$ by locally trivializing the family over generic parameters, proving the topological hypothesis of \Cref{thm: DH after partial reduction}.  We first fix generic parabolic weights and vary the Higgs field residues parameters, using the algebraic description from \Cref{sec: The Analytic Family}.  We then allow the weights and residues parameters to vary together.

\subsection{Variation of Higgs Field Residues}\label{sec: Variation of Higgs Field Residues}

For a fixed flag arrangement, let
\[
    V=\left\{\sigma=(\sigma_i^{(k)}):
        \sum_km_i^{(k)}\sigma_i^{(k)}=0
        \text{ for every }i\right\}
\]
be the vector space of Levi-central residue parameters.  Choose parabolic weights $\alpha$ for which semistability equals stability throughout this residue family, and suppose the moduli space is nonempty for at least one $\sigma\in V$.  These hypotheses hold for the exceptional-type flag arrangements below.  More generally, a sufficient condition for the existence of such weights is $\gcd(r_\star,\deg E,\{m_i^{(k)}\}_{i,k})=1$, by a similar argument as in \cite[Proposition~3.2]{BY99}.  We continue to assume $\deg E=0$, so this condition reduces to $\gcd(r_\star,\{m_i^{(k)}\}_{i,k})=1$.  The local trivialization argument below also applies to other fixed degrees, provided the corresponding moduli family is smooth and its joint Hitchin map is proper.

Write the algebraic moduli family as
\[
    \pi_\alpha\from\mathscr M_\alpha^{\mathrm{alg}}\longrightarrow V,
    \qquad \pi_\alpha^{-1}(\sigma)=\cM(\alpha,\sigma).
\]
Using the trace-free deformation complex gives the fixed-determinant version of the smooth family in \cite{LL26}.  The Hitchin fibration extends to the full family.  Let
\[
    \mathcal B_j=H^0\bigl(C,K_C(D)^{\otimes j}\bigr),
    \qquad \mathcal B=\bigoplus_{j=2}^{r_\star}\mathcal B_j,
\]
be the Hitchin base, and define
\begin{equation}
    \operatorname{Hit}\from\mathscr M_\alpha^{\mathrm{alg}}\longrightarrow\mathcal B,
    \qquad \operatorname{Hit}(\delbar_E,\varphi,\mathcal F)
        =\bigl(\operatorname{tr}(\wedge^2\varphi),\ldots,
            \operatorname{tr}(\wedge^{r_\star}\varphi)\bigr).
\end{equation}
Write $a_j(\varphi)$ for the higher traces $\operatorname{tr}(\wedge^j\varphi)\in\mathcal B_j$.  The joint map
\[
    (\pi_\alpha,\operatorname{Hit})
        \from\mathscr M_\alpha^{\mathrm{alg}}\longrightarrow V\times\mathcal B
\]
is proper by the argument in the proof of \cite[Proposition~4.19]{LL25}, restricted to the closed locus with fixed determinant and $\operatorname{tr}\varphi=0$.

Following \cite{FMSW26}, we consider the $\R_{>0}$-actions on the family and on $V\times\mathcal B$ given by
\begin{equation}
    \lambda\cdot(\delbar_E,\varphi,\mathcal F)=(\delbar_E,\lambda\varphi,\mathcal F),
    \qquad
    \lambda\cdot\bigl(\sigma,(a_j)_j\bigr)
        =\bigl(\lambda\sigma,(\lambda^j a_j)_j\bigr).
\end{equation}
The first action preserves stability, and $(\pi_\alpha,\operatorname{Hit})$ is equivariant for these actions.  Fix a point $x\in\mathscr M_\alpha^{\mathrm{alg}}$ and a sequence $\lambda_n\to0$.  The images of $\lambda_n\cdot x$ under $(\pi_\alpha,\operatorname{Hit})$ converge to $(0,0)$, so they lie in a compact subset of $V\times\mathcal B$.  Properness gives a convergent subsequence of $\lambda_n\cdot x$, with limit in $\pi_\alpha^{-1}(0)$.  Since $\pi_\alpha$ is a submersion, its image contains a neighborhood of zero.  Scaling then shows that $\pi_\alpha$ is surjective.

\bigskip

To transport homology, we seek horizontal lifts of paths in $V$ whose flows identify the fibers of $\pi_\alpha$.  Choosing a right inverse of $\de\pi_\alpha$ lifts vectors in $V$ locally, but the resulting path could escape to infinity in finite time.  Since the fibers of $(\pi_\alpha,\mathrm{Hit})$ are compact, we only need to prevent escape along the Hitchin base.  We address this by making a different choice of lift away from a compact set in $\mathscr M_\alpha^\mathrm{alg}$ which preserves a certain norm on $\mathcal B$.

\begin{theorem}\label{thm: scaling transport}\label{prop: family topology}
    For fixed generic ordered weights as above, the residue family $\pi_\alpha\from\mathscr M_\alpha^{\mathrm{alg}}\to V$ is locally trivial with smooth trivializations.  Its integral homology local systems have trivial monodromy.
\end{theorem}
\begin{proof}
    Choose norms on $\mathcal B_j$, set $L=\prod_{j=2}^{r_\star}j$, and define
    \begin{equation}\label{eq: Hitchin control function}
        \rho\from\mathscr M_\alpha^{\mathrm{alg}}\longrightarrow\R_{\geq0},\qquad
        \rho(x)=\sum_{j=2}^{r_\star}\|a_j(x)\|^{2L/j}.
    \end{equation}
    This function factors through $\mathrm{Hit}$, say $\rho=\bar\rho\circ\mathrm{Hit}$.  It is smooth and semi-invariant of degree $2L$, and $(\pi_\alpha,\rho)$ is proper because a bound on $\rho$ bounds every $a_j$.  Let $\zeta_1$ be the vector field on $\mathscr M_\alpha^\mathrm{alg}$ generated by the action of $1\in\mathrm{Lie}(\R_{>0})$.  On $\pi_\alpha^{-1}(0)$,
    \[
        \de\pi_\alpha(\zeta_1)=0,\qquad \de\rho(\zeta_1)=2L\rho.
    \]
    Since $\pi_\alpha$ is a submersion, these identities show that $(\pi_\alpha,\rho)$ is a submersion on $\pi_\alpha^{-1}(0)\sm\rho^{-1}(0)$.

    Let $Z$ be the critical locus of $(\pi_\alpha,\rho)$.  We claim $\rho$ is bounded on $Z\cap\pi_\alpha^{-1}(K)$ for any compact $K\subset V$.  For suppose $x_n$ is a sequence in $Z\cap\pi_\alpha^{-1}(K)$ with $\rho(x_n)\to\infty$, and set $\lambda_n=\rho(x_n)^{-1/(2L)}$.  Semi-invariance gives $\rho(\lambda_n\cdot x_n)=1$ and by equivariance $\pi_\alpha(\lambda_n\cdot x_n)=\lambda_n\pi_\alpha(x_n)\to0$.  The $\R_{>0}$-action is smooth and takes critical points to critical points, so properness of $(\pi_\alpha,\rho)$ gives a subsequence converging to a critical point in $\pi_\alpha^{-1}(0)\cap\rho^{-1}(1)$.  Since $Z$ is closed, this contradicts the previous paragraph.

    Fix $\sigma_*\in V$ and a neighborhood $U$ of $\sigma_*$ with compact closure.  The above says that for $R>0$ large enough, $(\pi_\alpha,\rho)$ is a submersion on $\pi_\alpha^{-1}(\overline U)\cap\{\rho>R\}$.  Choose a smooth right inverse $r_0$ to $\de\pi_\alpha$ on all of $\pi_\alpha^{-1}(U)$, and on $\pi_\alpha^{-1}(U)\cap\{\rho>R\}$ choose a smooth right inverse $r_1$ of $(\de\pi_\alpha,\de\rho)$.  Use cutoff functions to glue $r_0$ and $r_1\big|_{TV\times0}$ to get a smooth family of sections $s_x\from T_{\pi_\alpha(x)}V\to T_x\mathscr M_\alpha^{\mathrm{alg}}$ of $\de\pi_\alpha$ satisfying
    \begin{equation}\label{eq: controlled horizontal lifts}
        \de(\pi_\alpha)_x\circ s_x=\mathrm{id},\qquad
        \de\rho_x\circ s_x=0\quad\text{if }\rho(x)\geq R+1.
    \end{equation}

    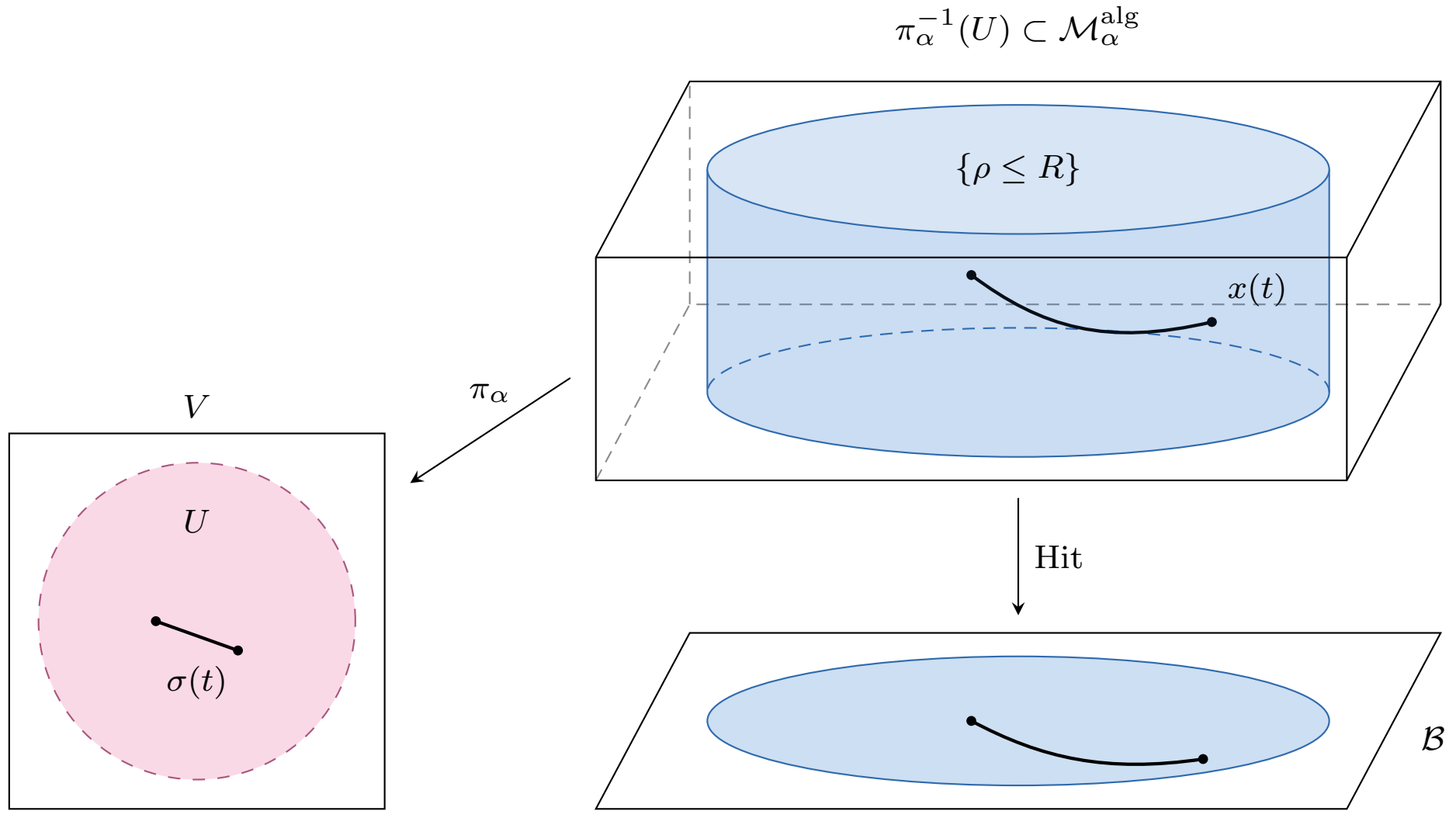
\begin{figure}[htbp]
        \centering
        \begin{tikzpicture}[>=stealth,font=\scriptsize]
            \definecolor{transportBlue}{RGB}{55,125,205}
            \begin{scope}[yshift=1.4cm]
                \draw (-5,-1) rectangle (-1.8,2.2);
                \node[above] at (-3.4,2.2) {$V$};
                \filldraw[dashed,fill=magenta!18,draw=magenta!65!black] (-3.4,0.6) circle (1.35cm);
                \node at (-3.4,1.45) {$U$};
                \draw[thick] (-3.75,0.6) -- (-3.05,0.35);
                \fill (-3.75,0.6) circle (1.2pt);
                \fill (-3.05,0.35) circle (1.2pt);
                \node[below] at (-3.4,0.35) {$\sigma(t)$};
            \end{scope}

            \draw[densely dashed,black!45] (0,3.2) -- (0.8,4.7) -- (7.2,4.7);
            \draw[densely dashed,black!45] (0.8,4.7) -- (0.8,6.6);
            \node[above] at (3.6,6.75) {$\pi_\alpha^{-1}(U)\subset\mathcal M_\alpha^{\mathrm{alg}}$};

            \path[fill=transportBlue,fill opacity=0.16]
                (0.95,3.95) arc (180:360:2.65cm and 0.55cm)
                -- (6.25,5.85) arc (0:180:2.65cm and 0.55cm) -- cycle;
            \draw[transportBlue!85!black,densely dashed]
                (6.25,3.95) arc (0:180:2.65cm and 0.55cm);
            \draw[thick] (3.2,4.95) .. controls (3.8,4.5) and (4.35,4.34) .. (5.25,4.55);
            \fill (3.2,4.95) circle (1.2pt);
            \fill (5.25,4.55) circle (1.2pt);
            \path[fill=transportBlue,fill opacity=0.12]
                (0.95,3.95) arc (180:360:2.65cm and 0.55cm)
                -- (6.25,5.85) arc (360:180:2.65cm and 0.55cm) -- cycle;
            \draw[transportBlue!85!black]
                (0.95,5.85) -- (0.95,3.95) arc (180:360:2.65cm and 0.55cm) -- (6.25,5.85);
            \filldraw[fill=transportBlue!20,draw=transportBlue!85!black]
                (3.6,5.85) ellipse (2.65cm and 0.55cm);
            \node at (3.6,5.85) {$\{\rho\leq R\}$};
            \node[above right] at (5.25,4.55) {$x(t)$};

            \draw (0,3.2) -- (6.4,3.2) -- (7.2,4.7) -- (7.2,6.6)
                -- (0.8,6.6) -- (0,5.1) -- cycle;
            \draw (0,5.1) -- (6.4,5.1) -- (7.2,6.6);
            \draw (6.4,5.1) -- (6.4,3.2);
            \draw[->,shorten <=4pt,shorten >=4pt] (-0.1,4.15) -- (-1.7,3.1) node[midway,above=3pt] {$\pi_\alpha$};

            \draw[->] (3.6,3.05) -- (3.6,2.05) node[midway,right] {$\mathrm{Hit}$};
            \draw (0,0.4) -- (6.4,0.4) -- (7.2,1.9) -- (0.8,1.9) -- cycle;
            \node[right] at (6.9,1) {$\mathcal B$};
            \filldraw[fill=transportBlue!25,draw=transportBlue!85!black]
                (3.6,1.15) ellipse (2.65cm and 0.55cm);
            \draw[thick] (3.2,1.15) .. controls (3.8,0.83) and (4.35,0.705) .. (5.175,0.825);
            \fill (3.2,1.15) circle (1.2pt);
            \fill (5.175,0.825) circle (1.2pt);
        \end{tikzpicture}
        \caption{The lift of a path $\sigma(t)$ in $U\subset V$ to a path $x(t)$, which must remain in the region $\{\rho\leq R\}$.  This prevents finite-time escape along $\mathcal B$.}
        \label{fig: controlled transport}
    \end{figure}
    Let $\sigma(t)$ be a smooth path in $U$ defined for $t\in[0,1]$.  A lift $x(t)$ must solve $\dot x(t)=s_{x(t)}(\dot\sigma(t))$.  Short-time existence is automatic, and by construction, $\rho(x(t))$ stays bounded where $x(t)$ is defined.  By properness of $(\pi_\alpha,\rho)$, the lift $x(t)$ therefore stays in compact set, and can therefore be extended until $t=1$.  For convex $U$, lifting the radial vector field out of $\sigma_*$ gives a smooth trivialization $\pi_\alpha^{-1}(U)\cong\pi_\alpha^{-1}(\sigma_*)\times U$.
    
    The diffeomorphisms used for transport along a path depend on the choice of horizontal lift, but their maps on homology depend only on the endpoints of the path.  The resulting integral homology local systems have trivial monodromy because $V$ is contractible.
\end{proof}

For the upcoming period calculations, we need the transported class of a smooth Hitchin fiber to have a complex submanifold representative.  We obtain such representatives by choosing the local trivialization to preserve a holomorphic level near the reference fiber.
\begin{corollary}\label{cor: holomorphic level transport}
    In the setting of \Cref{thm: scaling transport}, let $F\subset\cM(\alpha,\sigma_*)$ be a smooth compact regular Hitchin fiber.  Over a sufficiently small neighborhood $U$ of $\sigma_*$, the smooth trivialization can be chosen to take $F$ to a smooth compact complex submanifold $F_\sigma\subset\cM(\alpha,\sigma)$ for every $\sigma\in U$.  In particular, the transported class $[F]$ is represented by $F_\sigma$.
\end{corollary}
\begin{proof}
    Write $b=\mathrm{Hit}(F)$ and let $\mathrm{Hit}_{\alpha,\sigma_*}$ denote the restriction of $\mathrm{Hit}$ to $\mathcal M(\alpha,\sigma_*)$.  By \cite[Theorem~5.4]{LL26}, adapted to the fixed-determinant setting, the image of $\mathrm{Hit}_{\alpha,\sigma_*}$ is smooth and has complex dimension $\dim_{\C}F=\frac12\dim_{\C}\mathcal M(\alpha,\sigma_*)$.  Choose a linear projection $h\from\mathcal B\to T_b\im(\mathrm{Hit}_{\alpha,\sigma_*})$.  Regularity of $F$ and properness of $(\pi_\alpha,\mathrm{Hit})$ give neighborhoods $U$ of $\sigma_*$ and $W$ of $b$ on which $\mathscr F=(h\circ\mathrm{Hit})^{-1}(h(b))$ meets the fibers of $\pi_\alpha$ transversely in complex submanifolds of dimension $\dim_{\C}F$.

    Choose the horizontal lifts from \Cref{thm: scaling transport} tangent to $\mathscr F$ on $\pi_\alpha^{-1}(U)\cap\mathrm{Hit}^{-1}(W)$ (say, using another bump function around sufficiently small $W$).  Since $F$ is compact, for $\sigma$ sufficiently close to $\sigma_*$, the parallel transport of $F$ from $\sigma_*$ to $\sigma$ remains in $\mathrm{Hit}^{-1}(W)$.  The image $F_\sigma$ is a smooth compact submanifold of $\mathscr F\cap\pi_\alpha^{-1}(\sigma)\cap\mathrm{Hit}^{-1}(W)$ with the same real dimension.  Since $F_\sigma$ has no boundary, it is both closed and open in $\mathscr F\cap\pi_\alpha^{-1}(\sigma)\cap\mathrm{Hit}^{-1}(W)$, hence a complex submanifold.
\end{proof}

\subsection{Global Period Variation}\label{sec: Global Period Variation}

We now wish to vary the full triple $\nu=(\nu_1,\nu_2,\nu_3)\in\Lambda$ under the correspondence \eqref{eq: nu normalized definition} between $\nu$ and $(\alpha,\sigma)$.  Since the family $\pi_\alpha\from\mathscr M_\alpha^\mathrm{alg}\to V$ only depends on the stability chamber $W$ of $\alpha$, our topological transport extends to compact subsets of $W\times V$.  To allow the weights to cross chamber walls, we exclude the simultaneous weight and residue relations associated with possible stable factors, following Nakajima's definition of genericity \cite[Section~2]{Nak94}.

Let $g$ be the genus of $C$ and $n=|D|$.  If $(E,\varphi)$ is strictly semistable, choose a proper stable Jordan--Hölder factor $S$.  Write $d=\deg S$, $s=\rk(S)$, and $s_i^{(k)}$ for its graded multiplicities.  Since $S$ is stable, its deformation theory gives the necessary bound
\begin{equation}\label{eq: stable factor bound}
    (2-2g-n)s^2+\sum_{i,k}\bigl(s_i^{(k)}\bigr)^2\leq 2.
\end{equation}
This bound is automatic when $g\geq 1$, but can impose an additional restriction when $g=0$.  Since $S$ is a proper factor of $E$, its numerical data also satisfy
\begin{equation}\label{eq: rank requirements}
    0<s<r_\star,\qquad
    0\leq s_i^{(k)}\leq m_i^{(k)},\qquad
    \sum_k s_i^{(k)}=s\quad\text{for every }i.
\end{equation}

The factor $S$ has parabolic degree zero, and $\tr(\varphi|_S)$ is a meromorphic one-form on $C$.  The residue theorem therefore gives the simultaneous relations
\begin{equation}\label{eq: simultaneous walls}
    d+\sum_{i,k}s_i^{(k)}\alpha_i^{(k)}=0,
    \qquad
    \sum_{i,k}s_i^{(k)}\sigma_i^{(k)}=0.
\end{equation}

\begin{definition}[Generic parameter locus]\label{def: generic locus}
    A parameter $\nu\in\Lambda$ is \emph{generic} if the two relations \eqref{eq: simultaneous walls} do not hold simultaneously for any integers $d,s,(s_i^{(k)})$ satisfying \eqref{eq: rank requirements} and \eqref{eq: stable factor bound}.  We denote the \emph{generic parameter locus} by $\Lambda_0\subset\Lambda$.
\end{definition}


Only finitely many integers $d$ can occur in the first relation.  Under the coprimality condition \eqref{eq: coprime hypothesis}, the excluded locus is a finite union of intersections of $\Lambda$ with real codimension-3 affine subspaces, so $\Lambda_0$ is connected.

\begin{corollary}\label{cor: local triviality and homology transport}
    Assume the parabolic flag multiplicities satisfy
    \begin{equation}\label{eq: coprime hypothesis}
        \gcd\bigl(r_\star,\{m_i^{(k)}\}_{i,k}\bigr)=1,
    \end{equation}
    and assume $\mathscr M_0:=\pi^{-1}(\Lambda_0)$ is nonempty.  Then the family $\pi\from\mathscr M_0\to\Lambda_0$ is locally topologically trivial.  Its integral homology local systems have trivial monodromy.
\end{corollary}
\begin{proof}
    Fix $(\alpha_*,\sigma_*)\in\Lambda_0$ and write $\Lambda=P\times V$ where $P$ is the space of parabolic weights and $V$ is the space of Levi-central Higgs field residues.  Consider the pairs of (real and complex) hyperplanes in $P$ and $V$ defined by the equations in \eqref{eq: simultaneous walls}, indexed by numerical data satisfying \eqref{eq: stable factor bound} and \eqref{eq: rank requirements}.  Choose a neighborhood $W\times U\subset\Lambda_0$ of $(\alpha_*,\sigma_*)$ where $W\subset P$ is convex and only meets any hyperplanes which already pass through $\alpha_*$, and where $U\subset V$ remains disjoint from the corresponding hyperplanes in $V$.  Then by construction, for $\sigma\in U$ the stable Higgs bundle representatives are independent of $\alpha\in W$.

    By the hypothesis \eqref{eq: coprime hypothesis} on the GCD, there exists a generic rational $\alpha'\in W$ so that (shrinking $U$ if needed) \Cref{thm: scaling transport} gives the local trivialization of the family $\mathscr M_{\alpha'}^\mathrm{alg}$ over $\{\alpha'\}\times U$, which extends to $W\times U$ since stability is preserved.  The harmonic metrics of \Cref{sec: The Analytic Family} gives the same trivialization for the analytic family over $W\times U$.
    
    To identify the homology globally, we observe that $\Lambda_0$ is simply connected: $\Lambda$ is convex, and the excluded simultaneous relations form a finite union of affine subspaces of real codimension 3.  Thus, the homology local systems have trivial monodromy.
\end{proof}

We can now apply the affine variation theorem to the family over $\Lambda_0$.
\begin{theorem}[Affine variation of periods]\label{thm: general period variation}
    Let $\Lambda_0\subset\Lambda$ be the generic parameter locus.  Assume the coprimality condition \eqref{eq: coprime hypothesis} and that $\mathscr M_0=\pi^{-1}(\Lambda_0)$ is nonempty.  Fix $\nu_*\in\Lambda_0$ and identify cohomology throughout $\Lambda_0$ by parallel transport.  Let $\mathscr E_0=\tilde\pi^{-1}(\Lambda_0)$, let $i_{\nu_*}\from\cM_{\nu_*}\into\mathscr M_0$ be the fiber inclusion, and set $c=i_{\nu_*}^*c(\mathscr E_0)\in H_2(\cM_{\nu_*},\C)$.  Then, for every $\nu\in\Lambda_0$,
    \begin{equation}\label{eq: general period variation}
        [\omega_{s,\nu}]-[\omega_{s,\nu_*}]=-\langle c,\nu_s-\nu_{*,s}\rangle,\qquad s=1,2,3.
    \end{equation}
\end{theorem}
\begin{proof}
    By \Cref{prop: universal moduli space}, $q\from\mathscr E_0\to\mathscr M_0$ is a principal $B$-bundle.  The forms constructed in \Cref{sec: Universal Moduli Space} and the homology transport of \Cref{cor: local triviality and homology transport} satisfy the hypotheses of \Cref{thm: DH after partial reduction}, which gives the formula.
\end{proof}

\section{Characteristic Classes and Periods}
\label{sec: affine fiber alternative}

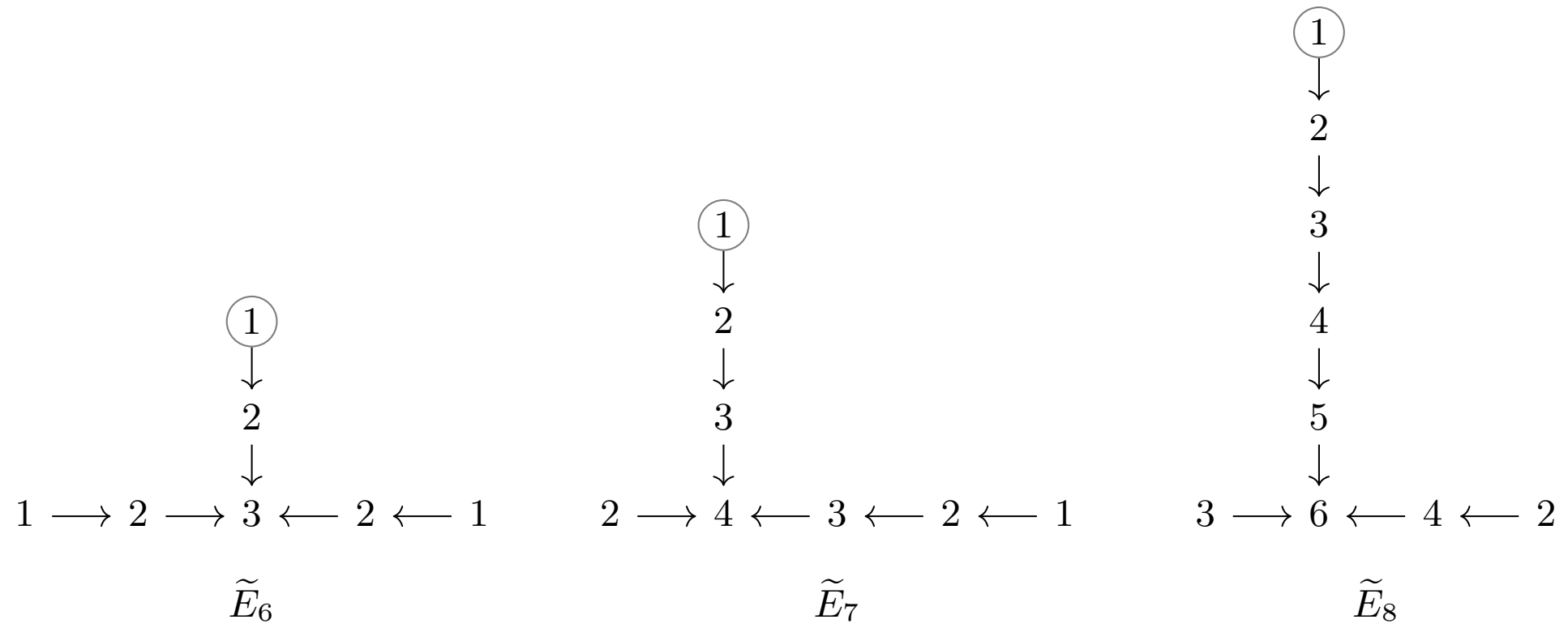
\begin{figure}[htbp]
    \centering
    \begin{minipage}[b]{0.29\textwidth}
        \centering
        \footnotesize
        \begin{tikzcd}[column sep=14pt,row sep=10pt]
            && |[circle,draw=gray,thin,inner sep=1.5pt]| 1 \arrow[d] \\
            && 2 \arrow[d] \\
            1 \arrow[r] & 2 \arrow[r] & 3 & 2 \arrow[l] & 1 \arrow[l]
        \end{tikzcd}
        \par\vspace{8pt} $\tilde E_6$
    \end{minipage}\hspace{0.015\textwidth}
    \begin{minipage}[b]{0.29\textwidth}
        \centering
        \footnotesize
        \begin{tikzcd}[column sep=14pt,row sep=10pt]
            & |[circle,draw=gray,thin,inner sep=1.5pt]| 1 \arrow[d] \\
            & 2 \arrow[d] \\
            & 3 \arrow[d] \\
            2 \arrow[r] & 4 & 3 \arrow[l] & 2 \arrow[l] & 1 \arrow[l]
        \end{tikzcd}
        \par\vspace{8pt} $\tilde E_7$
    \end{minipage}\hspace{0.015\textwidth}
    \begin{minipage}[b]{0.24\textwidth}
        \centering
        \footnotesize
        \begin{tikzcd}[column sep=14pt,row sep=10pt]
            & |[circle,draw=gray,thin,inner sep=1.5pt]| 1 \arrow[d] \\
            & 2 \arrow[d] \\
            & 3 \arrow[d] \\
            & 4 \arrow[d] \\
            & 5 \arrow[d] \\
            3 \arrow[r] & 6 & 4 \arrow[l] & 2 \arrow[l]
        \end{tikzcd}
        \par\vspace{8pt} $\tilde E_8$
    \end{minipage}
    \caption{The exceptional-type flag arrangements as labeled star-shaped quivers.  Deleting the circled vertex gives the corresponding finite Dynkin diagram.}
    \label{fig: E_n flag arrangements}
\end{figure}

In this section, we specialize to three families of parabolic Hitchin moduli spaces of exceptional-type $(E_6,E_7,E_8)$ over $C=\CP^1$ with three marked points $D=\{p_1,p_2,p_3\}$.  The flag dimensions are given by the corresponding affine quivers in \Cref{fig: E_n flag arrangements}, and the multiplicities are
\begin{equation}\label{eq: exceptional multiplicities}
    \begin{array}{c|c|c|c|c}
        \text{type}&r_\star&(m_1^{(k)})_k&(m_2^{(k)})_k&(m_3^{(k)})_k\\[1pt] \hline
        \noalign{\vskip 2pt}
        \tilde E_6&3&(1,1,1)&(1,1,1)&(1,1,1)\\[1pt]
        \tilde E_7&4&(1,1,1,1)&(1,1,1,1)&(2,2)\\[1pt]
        \tilde E_8&6&(1,1,1,1,1,1)&(2,2,2)&(3,3)
    \end{array}
\end{equation}
Due to the repetition, we write $m_i=m_i^{(k)}$.  By \Cref{thm: general period variation}, the period formula has the form $\mathbf A\nu+\mathbf b$ after chooses a cohomology basis, identified throughout by parallel transport.  The matrix $\mathbf A$ is determined by the characteristic class $c=c(\mathscr E)$ restricted to any given fiber, and $\mathbf b$ is a constant cohomology class corresponding to the limit $\nu\to0$.  We work at convenient choices of $\nu$ to carry out these calculations, relying on Groechenig's orbifold parameters for homology information \cite{Gro14}, and the quiver embedding $\mathcal T\from\mathcal X\into\mathcal M$ for evaluating the class $c$.

Our period calculation is motivated by the \emph{Modularization Conjecture}: every ALG or ALG$^*$ gravitational instanton should arise as a Hitchin moduli space with its natural hyperk\"ahler metric \cite[Conjecture~1.2]{FMSW26}.  For the exceptional types $E_6,E_7,E_8$, this requires both the appropriate ALG asymptotics and realization of the corresponding Torelli parameters.  The asymptotic analysis is the subject of upcoming work of Dimakis--Fredrickson--Mazzeo--Weiss.

\subsection{Homology Basis and Intersection Pairing}\label{sec: quiver lattice}

For each family of type $E_6,E_7,E_8$, the integral homology $H_2(\cM(\alpha,\sigma),\Z)$ has a basis of sphere classes corresponding to the vertices of the affine diagram (see \Cref{fig: E_n flag arrangements}).  Let $\widetilde I$ be the vertex set, with central vertex $v_\star$ and remaining vertices $v_i^{(k)}$, matching our flag indexing convention.  At $\sigma=0$ and suitable generic $\alpha$ near a certain orbifold weight discussed below, Groechenig's identification with the resolution of $T^*T/\Z_{r_\star}$ for an elliptic curve $T$ shows that the nilpotent fiber of the Hitchin map has Kodaira type $IV^*,III^*,II^*$, respectively \cite[Theorem~4.1]{Gro14}.\footnote{
    Tensoring by $\mathcal O(1)$ and centering the parabolic weights translates his degree convention into ours.
}
The reduced nilpotent fiber is a deformation retract \cite[Corollary~1.3.6]{HRV13}, so its components give a homology basis $\{[S_v]\mid v\in\widetilde I\}$.  Write $S_\star=S_{v_\star}$ and $S_i^{(k)}=S_{v_i^{(k)}}$.  The homology local system of \Cref{cor: local triviality and homology transport} lets us transport these classes over $\Lambda_0$.

The intersection matrix is the negative affine Cartan matrix $-\widetilde{\mathsf C}$, namely the adjacency matrix minus $2\Id$, and has a 1-dimensional kernel spanned by a weighted sum of the classes corresponding to the diagram vertices.  Specializing a smooth Hitchin fiber $F$ to the nilpotent fiber gives the following relation:
\begin{equation}\label{eq: fiber relation}
    [F]=r_\star[S_\star]
        +\sum_{i=1}^3\sum_{k=0}^{\ell_i-1}r_i^{(k)}[S_i^{(k)}].
\end{equation}
Transporting these classes throughout $\Lambda_0$, the fiber relation holds everywhere.  We call $[F]$ the fiber class of $\mathcal M(\alpha,\sigma)$ for any $(\alpha,\sigma)$.  

For $v\in\widetilde I$, write $r_v$ for the corresponding flag dimension $r_\star$ or $r_i^{(k)}$.  The next definition follows \cite{FMSW26}.
\begin{definition}\label{def: admissible basis}
    An \emph{admissible basis} $\mathfrak B=(b_v)_{v\in\widetilde I}$ for $\mathcal M(\alpha,\sigma)$ is an ordered integral basis satisfying
    \begin{equation}\label{eq: Cartan pairing}
        b_v\cdot b_w=-\widetilde{\mathsf C}_{vw},
        \qquad \sum_{v\in\widetilde I}r_vb_v=[F],
    \end{equation}
    where $\widetilde{\mathsf C}$ is the affine $\widetilde E_n$ Cartan matrix.  We write $\mathsf C$ for the finite $E_n$ Cartan matrix.
\end{definition}
The classes $[S_v]$ above form an admissible basis, and the affine Weyl group $W(\tilde E_n)$ and diagram automorphisms act on the set of admissible bases.  For later calculations, we will construct a different admissible basis extending a homology basis for an embedded quiver variety.

\begin{definition}\label{def: quiver parameters}
    For a pair $(\alpha,\sigma)$, the \emph{corresponding quiver parameters} are
    \begin{equation}\label{eq: beta tau definition}
        \beta_i^{(k)}=\alpha_i^{(k)}-\alpha_i^{(k+1)},
            \qquad \tau_i^{(k)}=\sigma_i^{(k)}-\sigma_i^{(k+1)}
    \end{equation}
    for $k=0,\dots,\ell_i-1$.  The \emph{extension to the central node} also includes
    \begin{equation}\label{eq: beta tau extension}
        \beta_\star=\sum_i\alpha_i^{(\ell_i)},
            \qquad\tau_\star=\sum_i\sigma_i^{(\ell_i)}.
    \end{equation}
    We also write $\beta_{v_i^{(k)}}=\beta_i^{(k)}$, $\beta_{v_\star}=\beta_\star$, and similarly for $\tau$.  The extension is defined to satisfy the trace-free condition
    \begin{equation}\label{eq: beta tau relations}
        \sum_{v\in\widetilde I}r_v\beta_v=0,
        \qquad \sum_{v\in\widetilde I}r_v\tau_v=0.
    \end{equation}
\end{definition}
Choose a small generic $\alpha$ satisfying
\begin{equation}\label{eq: small reference weights}
    \sum_{i=1}^3\sum_{k=0}^{\ell_i-1}r_i^{(k)}\beta_i^{(k)}<1.
\end{equation}
For such $\alpha$, there is an embedding $\cT\from\cX(\beta,\tau)\into\cM(\alpha,\sigma)$ from a star-shaped Nakajima quiver variety with underlying graph the affine Dynkin diagram \cite{Yae26}, whose image is the open locus $\cM_0$ consisting of holomorphically trivial Higgs bundle.

\begin{proposition}\label{prop: quiver fiber splitting}
    For any $(\alpha,\sigma)$ satisfying \eqref{eq: small reference weights}, the inclusion $j\from\cM_0(\alpha,\sigma)\into\cM(\alpha,\sigma)$ gives
    \[
        H_2(\cM,\Z)=j_*H_2(\cM_0,\Z)\oplus\Z[F].
    \]
    Evaluation on $[F]$ gives an isomorphism
    \[
        H^2(\cM,\Z)\longrightarrow H^2(\cM_0,\Z)\oplus\Z,
        \qquad u\longmapsto\bigl(j^*u,\langle u,[F]\rangle\bigr).
    \]
\end{proposition}
\begin{proof}
    Using parallel transport, it suffices to prove the statement for $\sigma=0$.  The deformation retract to the nilpotent fiber, which is a tree of 2-spheres, shows $H_1(\cM,\Z)=0$.  Since the fiber relation \eqref{eq: fiber relation} has a coefficient equal to 1, $[F]$ is primitive.  The affine Cartan pairing therefore has radical $\Z[F]$.  The quotient by $\Z[F]$ has intersection pairing given by the negative Cartan matrix $-\mathsf C$.  The quiver variety has an integral basis with the same intersection pairing \cite{Kro89,Nak94}.

    Inclusion preserves intersection numbers, and the negative definiteness of the intersection pairing implies $j_*$ is injective.  Thus, the composition
    \begin{equation}
        P\from H_2(\cM_0,\Z)\to H_2(\cM,\Z)\to H_2(\cM,\Z)/\Z[F]
    \end{equation}
    remains injective because $\Z[F]$ is the radical of the pairing.  Pick bases on $H_2(\cM_0,\Z)$ and $H_2(\cM,\Z)/\Z[F]$ whose intersection matrices are $-\mathsf C$.  Since $P$ preserves the pairing, its matrix with respect to these bases satisfies $P^\intercal(-\mathsf C)P=-\mathsf C$.  Taking determinants gives $(\det P)^2=1$, so $P$ is an integral isomorphism.  This gives the decomposition.  The cohomology isomorphism follows from $H_1(\cM,\Z)=0$ and the universal coefficient theorem.
\end{proof}

\begin{figure}[htbp]
    \centering
    \begin{tikzpicture}[pics/compacttorus/.style={code={
        \begin{scope}[rotate=-90,yscale=0.45]
            \draw[line width=0.70cm] circle (0.91cm);
            \draw[white,line width=0.67cm] circle (0.91cm);
            \draw[line width=0.70cm] (180:0.91cm) arc (180:360:0.91cm);
            \draw[white,line width=0.67cm,line cap=round] (180:0.91cm) arc (180:360:0.91cm);
        \end{scope}
    }}]
        \draw (-4.1,-1.5) -- (4.1,-1.5);
        \node[right] at (4.1,-1.5) {$\mathcal B$};
        \foreach \x in {-3.15,3.15} {
            \draw (\x,-0.42) -- (\x,-1.5);
            \fill (\x,-1.5) circle (1.2pt);
        }
        \ifFiguresThreeD
            \begin{scope}[shift={(3.15,0.84)}]
                \node[inner sep=0pt,outer sep=0pt,transform shape] at (0,0) {\includegraphics{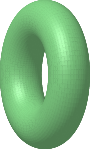}};
            \end{scope}
        \else
            \pic at (3.15,0.84) {compacttorus};
        \fi
        \pic at (-3.15,0.84) {compacttorus};
        \node[right] at (3.95,0.84) {$F$};
        \draw (0,-0.42) -- (0,-1.5);
        \fill (0,-1.5) circle (1.2pt);
        \foreach \x/\y in {-1.68/0,-0.84/0,0/0,0.84/0,1.68/0,0/0.84} {
            \begin{scope}[shift={(\x,\y)}]
                \ifFiguresThreeD
                    \node[inner sep=0pt,outer sep=0pt,transform shape] at (0,0) {\includegraphics{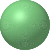}};
                \else
                    \filldraw[fill=green!25] (0,0) circle (0.42cm);
                \fi
            \end{scope}
        }
        \filldraw[fill=white] (0,1.68) circle (0.42cm);
        \draw[densely dashed] (0.42,1.68) arc (0:180:0.42cm and 0.12cm);
        \draw (-0.42,1.68) arc (180:360:0.42cm and 0.12cm);
    \end{tikzpicture}
    \caption{The nilpotent fiber and smooth Hitchin fibers in type $E_6$.  The image of the quiver variety homology basis together with a smooth Hitchin fiber generate $H_2(\mathcal M,\Z)$.}
    \label{fig: E6 Hitchin fibration}
\end{figure}
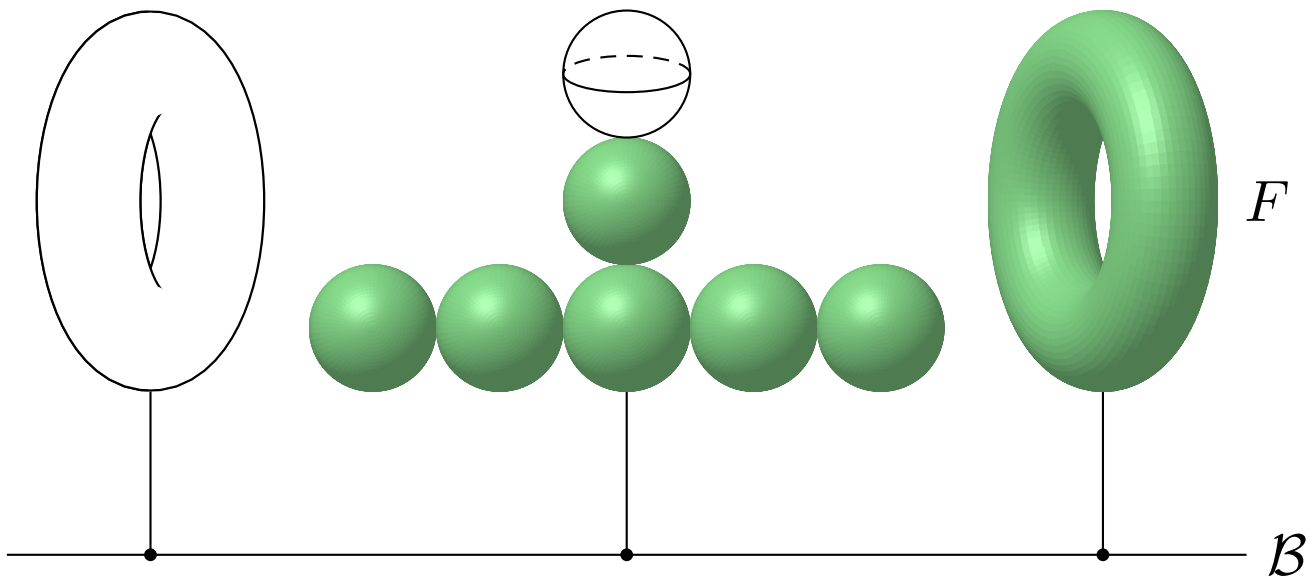

\subsection{The Boundary Characteristic Class}
\label{sec: characteristic class}

Keep $\nu$ and the corresponding $(\alpha,0)$ from above.  We now calculate
\[
    c=i_\nu^*c(\mathscr E)\in H^2(\cM,\C)\otimes_{\C}\fb_{\C},
\]
where $i_\nu\from\pi^{-1}(\nu)=\mathcal M(\alpha,0)\into\mathscr M$.  By \Cref{prop: quiver fiber splitting}, we can determine $c$ by its restriction to the quiver variety $\mathcal X$ and its evaluation on $[F]$.  We define line bundles on $\mathcal M$ relating to Nakajima's vertex bundles on $\mathcal X$.  The latter defines a dual basis for the homology basis for $\mathcal X$ described above.

For each $i,k$, a gauge transformation $g\in G$ acts on the line $\det(\gr F_i^{(k)})$ via multiplication by $\det(g_i^{(k)})$, where $g_i=\lim_{z\to p_i}g$ and $g_i^{(k)}$ is the induced action on $\gr F_i^{(k)}$.  As the Higgs bundle varies, these lines give a $G$-equivariant line bundle $L_i^{(k)}$ on the configuration space.  The ineffective subgroup $C_{r_\star}$ acts trivially on the $r_\star$th tensor power of $L_i^{(k)}$, which therefore descends to $\cM$.  We write
\begin{equation}\label{eq: block determinant classes}
    \lambda_i^{(k)}=\frac{1}{r_\star}c_1\bigl((L_i^{(k)})^{\otimes r_\star}\bigr)\in H^2(\cM,\mathbb Q).
\end{equation}
To express $c$ in terms of these classes, choose a principal $B$-connection $\Theta$ on $\mathscr E\to\mathscr M$ and restrict it to our reference fiber $\pi^{-1}(\nu)$.  Reading off the scalar on each block gives an identification
\[
    \mathfrak b\cong\left\{(a_i^{(k)})\in\bigoplus_{i,k}\I\R:
        \sum_{k=0}^{\ell_i}m_ia_i^{(k)}=0\text{ for each }i\right\}
\]
under which the curvature 2-form $F_\Theta$ can be decomposed into $\I\R$-valued components $F_{\Theta,i}^{(k)}$.  The action of $G$ on $(L_i^{(k)})^{\otimes r_\star}$ has weight $r_\star m_i$, so the induced connection from $\Theta$ is $\de-r_\star m_i\Theta_i^{(k)}$ (the minus sign comes from our left-action convention).
The Chern class of a bundle $V$ is $c_1(V)=\frac\I{2\pi}[\tr(F_V)]$, so
\[
    \lambda_i^{(k)}=\frac{m_i}{2\pi\I}[F_{\Theta,i}^{(k)}],
\]
and
\begin{equation}\label{eq: char class components}
    c=\left([F_{\Theta,i}^{(k)}]\right)_{i,k}
        =\left(\frac{2\pi\I}{m_i}\lambda_i^{(k)}\right)_{i,k}
\end{equation}
The fixed determinant of $E$ gives the relation $\sum_k\lambda_i^{(k)}=0$.

On the quiver side, consider Nakajima's tautological vertex bundles $\mathcal R_v$ over $\mathcal X$ \cite[Section~2]{Nak94}.  We normalize the collection of vertex bundles by $\mathcal R_a\cong\mathcal O$.  Write $t_v=c_1(\det\mathcal R_v)$, with $t_a=0$.

The lift of $\cT$ to configuration spaces identifies the non-central vertex spaces with the flag levels $F_i^{(k)}$ \cite{Yae26}.  As before, although the ineffective group $C_{r_\star}$ may act nontrivially on their determinant bundles, the $r_\star$th powers descend.  We define
\[
    q_i^{(k)}=\frac1{r_\star}c_1\bigl((\det\mathcal R_i^{(k)})^{\otimes r_\star}\bigr).
\]

\begin{lemma}\label{lem: vertex classes}
Under this renormalization, Nakajima's vertex bundle Chern classes $q_i^{(k)}\in H^2(\mathcal X,\Q)$ and our block determinant line bundle Chern classes $\lambda_i^{(k)}\in H^2(\mathcal M,\Q)$ satisfy
\begin{equation}\label{eq: quiver Chern class restriction}
    \cT^*\lambda_i^{(k)}=q_i^{(k)}-q_i^{(k-1)}
\end{equation}
where $q_i^{(-1)}$ and $q_i^{(\ell_i)}$ are understood to be 0.
\end{lemma}
\begin{proof}
The exact sequence
\[
    0\longrightarrow F_i^{(k-1)}\longrightarrow F_i^{(k)}
        \longrightarrow\gr F_i^{(k)}\longrightarrow0
\]
gives the relation
\[
    c_1(L_i^{(k)})=c_1\bigl(\det(\gr F_i^{(k)})\bigr)
        =c_1\bigl(\det\mathcal R_i^{(k)}\bigr)
        -c_1\bigl(\det\mathcal R_i^{(k-1)}\bigr).
\]
Dividing by $r_\star$ proves the lemma.
\end{proof}

\begin{figure}
    \centering
    \begin{minipage}[c]{0.43\textwidth}
        \centering
        \scriptsize
        \begin{tikzcd}[column sep=14pt,row sep=18pt]
            & |[text=gray!75]| v_1^{(0)} \arrow[d,gray!75] \\
            & v_1^{(1)} \arrow[d] \\
            & v_1^{(2)} \arrow[d] \\
            v_3^{(0)} \arrow[r] & \star & v_2^{(2)} \arrow[l] & v_2^{(1)} \arrow[l] & v_2^{(0)} \arrow[l]
        \end{tikzcd}
    \end{minipage}\hspace{0.01\textwidth}
    \begin{minipage}[c]{0.43\textwidth}
        \centering
        \begin{tikzpicture}[scale=1.2]
            \path[use as bounding box] (-1.26,-0.42) rectangle (2.94,2.94);
            \foreach \x/\y in {-0.84/0,0/0,0.84/0,1.68/0,2.52/0,0/0.84,0/1.68} {
                \begin{scope}[shift={(\x,\y)}]
                    \ifFiguresThreeD
                        \node[inner sep=0pt,outer sep=0pt,transform shape] at (0,0) {\includegraphics{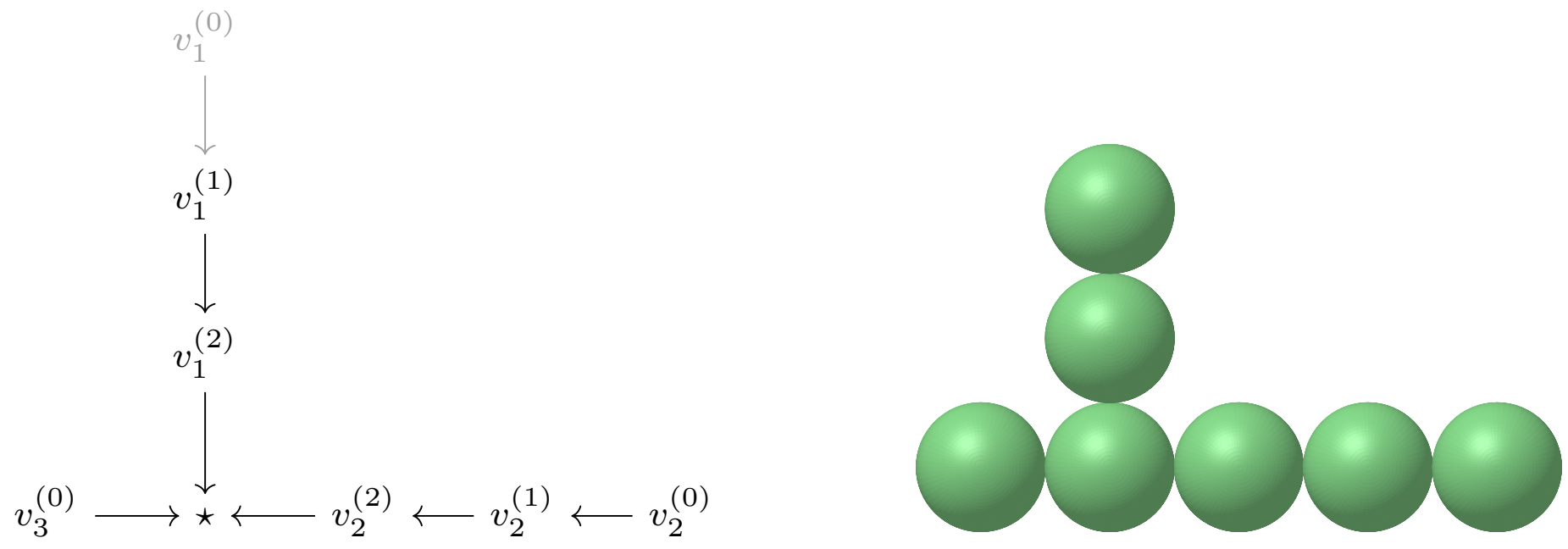}};
                    \else
                        \filldraw[fill=green!25] (0,0) circle (0.42cm);
                    \fi
                \end{scope}
            }
        \end{tikzpicture}
    \end{minipage}
    \caption{Deleting the grey rank-1 vertex from the affine $E_7$ quiver gives the intersection configuration of seven $2$-spheres.}
    \label{fig: E7 sphere configuration}
\end{figure}

Nakajima produces an integral basis $e_w\in H_2(\cX,\Z)$ indexed by the set $I$ obtained by deleting a rank-1 vertex from $\tilde I$ (see \Cref{fig: E7 sphere configuration}).  We will use $I=\widetilde I\sm\{v_1^{(0)}\}$.  The intersection pairing in this basis is given by the negative Cartan matrix as before.  In Nakajima's normalization $\mathcal R_{v_1^{(0)}}\cong\mathcal O$, the Chern classes of the remaining $\mathcal R_v$ form a dual basis $\{t_v\}_{v\in I}$ \cite{Nak01}.  Our normalization $q_v=t_v-(r_v/r_\star)t_\star$ only change on the pairings with $e_\star$.  Writing $e_i^{(k)}=e_{v_i^{(k)}}$, the classes duality pairing for the $q_i^{(k)}$ is
\begin{equation}
    \left\langle q_i^{(k)},e_j^{(h)}\right\rangle
        =\delta_{ij}\delta_{kh},
    \qquad
    \left\langle q_i^{(k)},e_\star\right\rangle
        =-\frac{r_i^{(k)}}{r_\star}.
\end{equation}
Here $0\leq k<\ell_i$ and $v_j^{(h)}\in I\sm\{\star\}$, so the omitted vertex $v_1^{(0)}$ is excluded.  Thus \Cref{lem: vertex classes} gives
\begin{align}\label{eq: ALE duality pairing}
    \bigl\langle\cT^*\lambda_i^{(k)},e_j^{(h)}\bigr\rangle
        &=\delta_{ij}\bigl(\delta_{kh}-\delta_{k-1,h}\bigr),\\
    \bigl\langle\cT^*\lambda_i^{(k)},e_\star\bigr\rangle
        &=\delta_{k,\ell_i}-\frac{m_i^{(k)}}{r_\star}.\notag
\end{align}
Here $0\leq k\leq\ell_i$ and $v_j^{(h)}\in I\sm\{\star\}$.  By \eqref{eq: char class components}, this gives the pairing between $c=i_\nu^*c(\mathscr E)$ and $H_2(\mathcal M_0,\Z)\subset H_2(\mathcal M,\Z)$.  It remains to pair $c$ with a smooth Hitchin fiber $F$.

\begin{lemma}\label{lem: fiber slope}
    The boundary characteristic class $c=i_\nu^*c(\mathscr E)$ satisfies $\langle c,[F]\rangle=0$.
\end{lemma}
\begin{proof}
    Fix a parameter $\nu$ corresponding to a pair $(\alpha,\sigma)$, and let $\nu'$ correspond to a pair $(\alpha,\sigma')$ for $\sigma'$ close to $\sigma$.  Write $F$ for a smooth Hitchin fiber for $\nu$.  By \Cref{cor: holomorphic level transport}, the transport of $[F]$ from $\nu$ to $\nu'$ can be represented by a smooth complex curve $F_{\sigma'}$.  This implies the holomorphic symplectic form $\Omega_{\nu'}=\omega_{2,\nu'}+\I\omega_{3,\nu'}$ restricts to zero on $F_{\sigma'}$, so \Cref{thm: general period variation} gives
    \begin{equation}
        0=\int_{F_{\sigma'}}\Omega_{\nu'}
            -\int_F\Omega_\nu
          =-\bigl\langle(\nu'_2-\nu_2)
                +\I(\nu'_3-\nu_3),c(F)\bigr\rangle.
    \end{equation}
    As $\sigma'$ varies near $\sigma$, this forces $c(F)$ to pair trivially with every element of $\mathfrak b^*\otimes\C$, so $c(F)=0$.
\end{proof}

The affine variation formula and \Cref{lem: fiber slope} implies that throughout $\Lambda_0$,
\begin{equation}\label{eq: constant fiber periods}
    \left(\int_F\omega_1,\int_F\Omega\right)=(A,0)
\end{equation}
for a constant $A>0$.  Indeed, at the reference parameter $F$ is a compact complex curve, so its volume with respect to $\omega_1$ is positive and its holomorphic symplectic period vanishes.

\subsection{Slope of Period Variation}\label{sec: period variation}

We now use the characteristic class to compute period variation in an admissible homology basis extending the image of Nakajima's basis under $\mathcal T_*$.  Call the deleted vertex $a=v_1^{(0)}\in\tilde I\sm I$, and set $\eta_v=\cT_*e_v$ for $v\in I$.  Since $r_a=1$, we can complete these classes by
\begin{equation}\label{eq: affine marked completion}
    \eta_a=[F]-\sum_{v\in I}r_v\eta_v.
\end{equation}
By \Cref{prop: quiver fiber splitting}, $\mathfrak B=(\eta_v)_{v\in\widetilde I}$ is an integral basis, and the intersection submatrix with $\eta_a$ removed is $-\mathsf C$.  Since $[F]$ is in the radical of the pairing, the fiber relation \eqref{eq: fiber relation} determines the remaining pairings and gives the negative affine Cartan matrix $-\widetilde{\mathsf C}$.  Thus $\mathfrak B$ is admissible in the sense of \Cref{def: admissible basis}.  We transport this basis throughout $\Lambda_0$ using \Cref{cor: local triviality and homology transport}.

For a deformation $(\dot\alpha,\dot\sigma)$ of $(\alpha,\sigma)$, consider the variation of quiver parameters extended to the central node:
\begin{equation}\label{eq: beta tau variations}
    \begin{aligned}
        \dot\beta_i^{(k)}&=\dot\alpha_i^{(k)}-\dot\alpha_i^{(k+1)},
        &\dot\beta_\star&=\sum_i\dot\alpha_i^{(\ell_i)},\\
        \dot\tau_i^{(k)}&=\dot\sigma_i^{(k)}-\dot\sigma_i^{(k+1)},
        &\dot\tau_\star&=\sum_i\dot\sigma_i^{(\ell_i)}.
    \end{aligned}
\end{equation}
As in \Cref{def: quiver parameters}, we write $\dot\beta_{v_i^{(k)}}=\dot\beta_i^{(k)}$, $\dot\beta_{v_\star}=\dot\beta_\star$, and similarly for $\dot\tau$.  We still have
\begin{equation}\label{eq: affine parameter relations}
    \sum_{v\in\widetilde I}r_v\dot\beta_v=0,
    \qquad \sum_{v\in\widetilde I}r_v\dot\tau_v=0.
\end{equation}
Although these numbers come from the quiver embedding, the following proposition applies throughout $\Lambda_0$ by parallel transport.

\begin{proposition}\label{prop: affine period slopes}
    Let $D$ denote the Gauss--Manin connection on the degree-two cohomology bundle over $\Lambda_0$.  Let $\dot\nu\in T_\nu\Lambda_0$ correspond to variations $(\dot\alpha,\dot\sigma)$ of the weights and residues, and let $(\dot\beta,\dot\tau)$ be the variations in \eqref{eq: beta tau variations}.  Let $\Omega=\omega_2+\I\omega_3$.  In the admissible homology basis $\mathfrak B$,
    \begin{equation}\label{eq: alternative period slopes}
            \bigl\langle D_{\dot\nu}[\omega_1],\eta_v\bigr\rangle=-4\pi^2\dot\beta_v,
                \qquad\bigl\langle D_{\dot\nu}[\Omega],\eta_v\bigr\rangle=-8\pi^2\dot\tau_v
    \end{equation}
    for all vertices $v\in\tilde I$.  Moreover,
    \begin{equation}\label{eq: fiber period slopes}
        \bigl\langle D_{\dot\nu}[\omega_1],[F]\bigr\rangle=0,
        \qquad\bigl\langle D_{\dot\nu}[\Omega],[F]\bigr\rangle=0.
    \end{equation}
\end{proposition}
\begin{proof}
    By \eqref{eq: general period variation}, the characteristic class $c$ is constant under transport and determines the slopes throughout $\Lambda_0$ by $D_{\dot\nu}[\omega_s]=-\langle\dot\nu_s,c\rangle$.  It therefore suffices to use the pairings computed at the reference parameter used above.  The variation $u=\dot\nu_s$ decomposes into $u_i=\bigoplus_k\I u_i^{(k)}\Id$ with the trace-free condition $\sum_k m_iu_i^{(k)}=0$.  The pairing formulas \eqref{eq: ALE duality pairing} give
    \begin{equation}
        \begin{aligned}
            \bigl\langle D_{\dot\nu}[\omega_s],\eta_{v_i^{(k)}}\bigr\rangle
                &=-2\pi\bigl(u_i^{(k)}-u_i^{(k+1)}\bigr),\\
            \bigl\langle D_{\dot\nu}[\omega_s],\eta_{v_\star}\bigr\rangle
                &=-2\pi\sum_i u_i^{(\ell_i)}.
        \end{aligned}
    \end{equation}
    The normalization \eqref{eq: nu normalized definition} identifies $u_i^{(k)}$ with $2\pi\dot\alpha_i^{(k)}$, $4\pi\Re\dot\sigma_i^{(k)}$, and $4\pi\Im\dot\sigma_i^{(k)}$ for $s=1,2,3$, respectively, giving the desired formulas.  The pairings with $[F]$ follow from \eqref{eq: affine marked completion} and \eqref{eq: affine parameter relations}.
\end{proof}

\subsection{Explicit Period Formula}
\label{sec: quiver calibration}

The numbers $\beta_i^{(k)}$ and $\tau_i^{(k)}$ give coordinates on $\Lambda_0$.  In a fixed cohomology basis, the affine period map has the form $\mathbf A\nu+\mathbf b$, where $\mathbf A$ encodes the characteristic classes calculated above, and $\mathbf b$ is a fixed cohomology class representing the limit as $\beta\to0,\tau\to0$.  It suffices to calculate the period for any fixed reference parameter.  First, we study the limit $\beta\to0$ directly with $\sigma=0$ by comparing to the quiver variety, whose period goes to 0 in the limit.  Taking the limit along such a ray is particularly convenient, since the complex period $[\omega_2]+\I[\omega_3]$ vanishes.  To evaluate $[\omega_1]$ on the additional homology class, we choose a different reference parameter to evaluate the period on the Hitchin fiber class $[F]$, which is constant throughout $\Lambda_0$ by \eqref{eq: constant fiber periods}.

Set $\sigma=0$, and choose a small generic $\alpha$ for which the quiver embedding $\mathcal T$ exists and stable parabolic bundles with zero Higgs field exist.  Such weights exist by \cite[Proposition~4.4]{King94} and the embedding $\mathcal T$.  For all $t\in(0,1]$, the second cohomology of $\mathcal M_0(t\alpha,0)$ is generated by the spheres in the nilpotent fiber.  Under the $\U(1)$-action $\lambda\cdot(\delbar_E,\varphi,\mathcal F)=(\delbar_E,\lambda\varphi,\mathcal F)$, the central sphere is pointwise fixed, and each remaining sphere has two fixed points \cite[Section~2]{Hit21}.

For $t\in(0,1]$, write $\omega_{1,t}$ for the K\"ahler form of $\mathcal M(t\alpha,0)$.  The $\U(1)$-action is Hamiltonian for $\omega_{1,t}$, and admits a moment map $\Psi_t\from\mathcal M(t\alpha,0)\to\R$ given by
\begin{equation}
    \Psi_t(\delbar_E,\varphi,\mathcal F)=\frac\I2\int_C\tr(\varphi^{\dagger_{h_t}}\wedge\varphi).
\end{equation}
Here $h_t$ is the harmonic metric for $(\delbar_E,\varphi,\mathcal F)$, adapted to $t\alpha$.  This is a Morse--Bott function by \cite{GGM04}.

\begin{lemma}\label{lem: vanishing quiver periods}
    Let $\alpha$ be the parabolic weights chosen above, and let $\eta_v\in H_2(\cM_0,\C)$ for $v\in I$ be the classes described in \Cref{sec: period variation}.  Using parallel transport to identify homology,
    \[
        \lim_{t\to0}\bigl\langle[\omega_{1,t}],\eta_v\bigr\rangle=\lim_{t\to0}\int_{\eta_v}\omega_{1,t}=0.
    \]
\end{lemma}
\begin{proof}
    Our choice of $\alpha$ contains the $\varphi=0$ locus as the $\U(1)$-fixed sphere, which we address last.  The remaining spheres each have two fixed points.  At these points, \cite[Proposition~3.4]{GGM04} gives a splitting as parabolic holomorphic bundles
    \[
        E=\bigoplus_{j=1}^{r_\star} E_j,
        \qquad \varphi_j:=\varphi\big|_{E_j}\from E_j\longrightarrow E_{j+1}\otimes K(D),\quad 1\leq j<r_\star.
    \]
    Here zero summands are allowed and $\varphi_{r_\star}=0$.  Consider the action of $\U(1)$ on $E=\bigoplus_j E_j$ given by multiplication by $\lambda^j$ on $E_j$.  This action sends $\varphi$ to $\lambda\varphi$, so uniqueness of the harmonic metric implies that it preserves $h_t$.  The summands $E_j$ have distinct weights under this action and are therefore orthogonal for $h_t$.

    This decomposition lets us evaluate $\Phi$ using the curvature of each summand.  Define the $j$th component
    \begin{equation}\label{eq: components of Morse Bott function}
        \Psi_{j,t}=\frac{\I}{2}\int_{C\setminus D}
            \tr\bigl(\varphi_j^{\dagger_{h_t}}\wedge\varphi_j\bigr),
        \qquad \Psi_{0,t}=\Psi_{r_\star,t}=0.
    \end{equation}
    Since the splitting is orthogonal and the determinant metric is flat, the component of Hitchin's equation on $E_j$ is
    \[
        F_{j,t}
            +\varphi_j^{\dagger_{h_t}}\wedge\varphi_j
            +\varphi_{j-1}\wedge\varphi_{j-1}^{\dagger_{h_t}}=0,
    \]
    where $\varphi_0=\varphi_{r_\star}=0$.  Taking traces,
    \[
        \tr F_{j,t}
            =\tr\bigl(\varphi_{j-1}^{\dagger_{h_t}}\wedge\varphi_{j-1}\bigr)
            -\tr\bigl(\varphi_j^{\dagger_{h_t}}\wedge\varphi_j\bigr).
    \]
    Substituting \eqref{eq: components of Morse Bott function} and applying parabolic Chern--Weil theory yields
    \begin{equation}
        \pi\operatorname{pardeg}_{t\alpha}(E_j)
            =\frac{\I}{2}\int_{C\setminus D}\tr F_{j,t}
            =\Psi_{j-1,t}-\Psi_{j,t}.
    \end{equation}
    The parabolic degree is in terms of the induced parabolic structure on $E_j$.  Taking the weighted sum and telescoping,
    \begin{equation}
        \pi\sum_{j=1}^{r_\star} j\operatorname{pardeg}_{t\alpha}(E_j)
            =\sum_{j=1}^{r_\star}j(\Psi_{j-1,t}-\Psi_{j,t})
            =\sum_{j=1}^{r_\star-1}\Psi_{j,t}
            =\Psi_t.
    \end{equation}

    Since $E$ is holomorphically trivial, each $E_j$ must have topological degree 0.  Therefore the parabolic degrees are just the sums of the weights, so $\Psi_t\in O(t)$ at the fixed point.
    
    Let $x_-,x_+$ be the $\U(1)$-fixed points of a sphere $S$ ordered by increasing $\Psi_t$.  Let $d_S$ be the order of the stabilizer of $S$.  Integrating in action-angle coordinates gives
    \[
        \int_S\omega_{1,t}
            =\frac{2\pi}{d_S}\bigl(\Psi_t(x_+)-\Psi_t(x_-)\bigr)\in O(t).
    \]

    It remains to address the $\U(1)$-fixed sphere, where $\varphi=0$ and Hitchin's equations reduce to the flat unitary connection equation.  Jeffrey's small-weight comparison theorem \cite[Theorem~6.6]{Jef94} identifies these moduli spaces, for sufficiently small $t$, with a fixed symplectic manifold whose symplectic form is multiplied by $t$.  Hence the sphere's volume with respect to $\omega_{1,t}$ tends to zero linearly in $t$.
\end{proof}

To compute the remaining fiber period $\langle[\omega_1],[F]\rangle$, we use the orbifold weights associated with the elliptic quotients of orders $3,4,6$ for types $E_6,E_7,E_8$, respectively (see \cite[Section~4]{Gro14}).  Here $\sigma=0$ and, translating to our degree-zero and trace-free convention, the parabolic weights are
\begin{equation}\label{eq: orbifold weight table}
    \begin{array}{c|c|c|c}
        \text{type}&\alpha_1&\alpha_2&\alpha_3\\ \hline
        E_6&\frac13(1,0,-1)&\frac13(1,0,-1)&\frac13(1,0,-1)\\[3pt]
        E_7&\frac18(3,1,-1,-3)&\frac18(3,1,-1,-3)&\frac14(1,-1)\\[3pt]
        E_8&\frac1{12}(5,3,1,-1,-3,-5)&\frac13(1,0,-1)&\frac14(1,-1).
    \end{array}
\end{equation}

\begin{lemma}\label{lem: Hitchin fiber volume}
    The Hitchin fiber class satisfies
    \[
        \bigl\langle[\omega_1],[F]\bigr\rangle=4\pi^2,
        \qquad
        \bigl\langle[\Omega],[F]\bigr\rangle=0.
    \]
\end{lemma}
\begin{proof}
    The complex period was shown to be zero in \eqref{eq: constant fiber periods}.  For the real period, let $T\to C$ be the elliptic cover and $T^\vee=\operatorname{Pic}^0(T)$.  As the parameter approaches the orbifold parameter in \cite[Section~4]{Gro14}, the exterior spheres collapse to zero volume and the central sphere survives as $S_\star=T^\vee/\Z_{r_\star}$.  The harmonic bundles descend from orthogonal sums of flat line bundles on $T$, so the induced $L^2$ metric is the flat quotient metric whose pullback to $T^\vee$ has volume $4\pi^2$.
 
    These orbifold parameters are not in $\Lambda_0$, but the points in a smooth Hitchin fiber $F$ remain stable here, so smoothness of $\mathscr M^\mathrm{irr}$ identifies the periods with the limits from nearby admissible parameters.  Thus, the fiber relation \eqref{eq: fiber relation} shows $F$ has volume $4\pi^2$ as well.
\end{proof}

\begin{theorem}\label{thm: periods}
    Let $\nu\in\Lambda_0$ encode weights and residues $(\alpha,\sigma)$, and define $(\beta,\tau)$ as in \Cref{def: quiver parameters}.  Let $\mathfrak B$ be the admissible basis for $H_2(\mathcal M(\alpha,\sigma),\C)$ from \Cref{sec: period variation}, indexed by the affine Dynkin vertex set $\tilde I$ with distinguished vertex $a=v_1^{(0)}$.  For all $v\in\tilde I$,
    \begin{equation}\label{eq: calibrated quiver periods}
        \int_{\eta_v}\omega_1=4\pi^2\delta_{va}-4\pi^2\beta_v,
        \qquad
        \int_{\eta_v}\Omega=-8\pi^2\tau_v.
    \end{equation}
\end{theorem}
\begin{proof}
    We know that for $\tau=0$ the complex period vanishes, as the nilpotent fiber generates degree-2 homology and consists of complex curves.  Combining the period slopes of \Cref{prop: affine period slopes} with the reference periods of \Cref{lem: vanishing quiver periods} and \Cref{lem: Hitchin fiber volume} proves the theorem.
\end{proof}

\begin{remark}\label{rem: D4}
    Type $D_4$ parabolic Hitchin moduli spaces parameterize rank 2 bundle with four marked points and full flags.  Its period calculation is already contained in \cite{FMSW26}, so we do not cover it in depth here.  Nevertheless, our same construction applies and the analogous formulas hold.
\end{remark}

The affine $\tilde E_n$ diagrams have $n+1$ vertices.  The real numbers $\beta_v$ and complex numbers $\tau_v$ for $v\in\tilde I$ are only subject to the trace-free constraints \eqref{eq: beta tau relations}, so $\beta$ has $n$ real degrees of freedom and $\tau$ has $n$ complex degrees of freedom.  Thus the linearized period map has full real rank $3n$.  Thus, the marked period map is a local affine diffeomorphism onto an open subset of the space of period triples with fiber period $(4\pi^2,0,0)$. 

Using \Cref{thm: periods}, local surjectivity can almost be upgraded to global surjectivity.  Following the strategy of \cite[Section~9]{FMSW26}, we vary the admissible marking by the affine Weyl group and use its action to move the real periods into the image of the allowed weight domain.  This potentially excludes a real codimension-2 subset of allowable parameters.
\begin{theorem}\label{thm: almost-global surjectivity}
In the $E_6, E_7, E_8$ cases, as $\alpha, \sigma$ vary, the classes $[\omega_1],[\omega_2],[\omega_3]$ nearly exhaust the allowable Torelli parameters in the corresponding family of ALG spaces.  Namely, the image excludes a set of real codimension at least 2.
\end{theorem}

\begin{remark}
    In the $E_6,E_7,E_8$ cases, within the locus where parabolic weights remain strictly ordered, the genericity condition in \Cref{def: generic locus} corresponds under \eqref{eq: calibrated quiver periods} to
    \[
        \left(\int_\gamma\omega_1,\int_\gamma\Omega\right)
        \ne(0,0)
        \qquad
        \text{for every integral class }\gamma
        \text{ with }\gamma^2=-2.
    \]
    This is precisely the nondegeneracy condition defining the ALG period domain in \cite[Theorem~1.7]{CVZ24}.  In particular, the bound in \eqref{eq: stable factor bound} crucially imposes self-intersection $-2$.
\end{remark}

\begin{remark}\label{rem: proposed full surjectivity}
    We briefly explain the complement of the image of the period map.  At each puncture we take strictly ordered weights of total spread less than one.  Our parabolic weight domain is a product of \emph{open} simplices containing $720$, $10{,}080$, and $483{,}840$ affine Weyl alcoves (i.e. parabolic weight chambers), respectively. The affine Weyl group translates of the  closure of the parabolic weight domain cover a real vector space. For \emph{interior} parabolic weights $\alpha$, $(\alpha,\sigma)$ can be generic even when $(\alpha,0)$ is not.  However, some real parameter orbits never intersect the interior of the parabolic weight domain.  There are the generic pairs $(\alpha, \sigma)$ that we miss.

    The same issue occurs in type $D_4$.  Writing $\mathcal R$ for the generic strict-weight parameters, Fredrickson--Mazzeo--Swoboda--Weiss distinguish
\[
    \widetilde{\mathcal R}
    =\bigcup_{w\in W_{\mathrm{aff}}}w\mathcal R
    \subsetneq
    \widetilde{\mathcal R}^{\mathrm{full}}
\]
    where $\widetilde{\mathcal R}^{\mathrm{full}}$ is the preimage of Chen--Viaclovsky--Zhang's period domain under the Torelli map, while $    \widetilde{\mathcal R}$ is the image of $\mathcal{R}$ under the affine Weyl group $W_{\mathrm{aff}}$.  Their difference has real codimension~$2$; see \cite[Remark~9.25]{FMSW26}.
\end{remark}

\bigskip

For parameters admitting the quiver embedding (see \eqref{eq: small reference weights}), the induced map $\mathcal T_*\from H_2(\mathcal X,\C)\into H_2(\mathcal M,\C)$ is injective, and its image has integral basis $\{\eta_v\mid v\in I\}$.  By \Cref{thm: periods}, the pulled-back K\"ahler forms have the same homology pairings as the forms $\omega_s^\mathrm{quiv}$ on the quiver variety \cite{Kro89}, up to a normalization factor.

\begin{corollary}\label{cor: quiver cohomology comparison}
    For the exceptional types $E_6,E_7,E_8$, let $(\alpha,\sigma)$ be a pair admitting the quiver embedding $\mathcal T\from\mathcal X(\beta,\tau)\into\mathcal M(\alpha,\sigma)$.  In the normalization and hyperk\"ahler conventions of \cite{Yae26}, we have
    \[
        \mathcal T^*[\omega_s]=2\pi[\omega_s^{\mathrm{quiv}}],\qquad s=1,2,3.
    \]
\end{corollary}

\appendix
\section{Quasi-Hamiltonian Reduction}\label{sec: quasi Hamiltonian DH}

One of our ideas along the way to proving affine period variation for parabolic Higgs bundles was to pass through the nonabelian Hodge correspondence to the Betti moduli space, which admits a construction as a quasi-Hamiltonian quotient \cite{Yam08}.  At that time, we proved a corollary of \Cref{thm: DH after partial reduction} which applies to group-valued moment maps.  Restricting to a central logarithmic chart gives the necessary geometric hypothesis, and the topological hypothesis would follow from a similar Laplacian estimate as in \Cref{sec: multi Taub NUT}, which became quite a challenge in this case.  Here, we present just the corollary.

Let $G$ be a finite-dimensional Lie group with an invariant nondegenerate symmetric bilinear form $(\ ,\ )$ on $\mathfrak g$ over $\mathbb F=\R$ or $\C$.  Write $\theta_{\mathrm L}$ and $\theta_{\mathrm R}$ for the Maurer--Cartan forms and $\chi=\frac1{12}(\theta_{\mathrm L},[\theta_{\mathrm L},\theta_{\mathrm L}])$ for the Cartan three-form.  A quasi-Hamiltonian $G$-space $(M,\omega,\mu)$ has an invariant two-form $\omega$ and an equivariant map $\mu\from M\to G$ satisfying
\begin{equation}\label{eq: quasi Hamiltonian axioms}
    \begin{aligned}
        \de\omega&=-\mu^*\chi,\\
        \iota_{\zeta_\xi}\omega&=\frac12\mu^*(\theta_{\mathrm L}+\theta_{\mathrm R},\xi),\\
        \ker\omega_x&=\{\zeta_\xi(x):(\operatorname{Ad}_{\mu(x)}+1)\xi=0\}.
    \end{aligned}
\end{equation}
These are the conventions of \cite[Definition~2.2]{AMM98}.  In the complex case we use holomorphic data as in \cite[Definition~2.11]{Yam08}, with smooth transport on the underlying manifolds.

Fix $z_*\in Z(G)$ and a subspace $\Lambda\subseteq\mathfrak g^G$, which we identify with a subspace of $(\mathfrak g^*)^G$ using the bilinear form.  Choose an open region $\Lambda_0\subseteq\Lambda$ such that
\[
    e_*\from\Lambda_0\longrightarrow Z(G),\qquad
    \alpha\longmapsto z_*\exp(\alpha)
\]
is a diffeomorphism onto its image, which we call $C_0$.  Set $E_0=\mu^{-1}(C_0)$ and $\tilde\pi=e_*^{-1}\circ\mu$.  Let $G^{\mathrm{eff}}$ be the quotient by the kernel of the action on $E_0$, and assume $\Lambda$ annihilates the ineffective Lie algebra, so its elements define covectors on $\mathfrak g^{\mathrm{eff}}$.

\begin{theorem}\label{thm: quasi Hamiltonian DH}
    Suppose $E_0$ is smooth and $q\from E_0\to\mathscr M_0=E_0/G$ is a smooth principal $G^{\mathrm{eff}}$-bundle over a paracompact manifold.  Suppose the induced map $\pi\from\mathscr M_0\to\Lambda_0$ is a smooth surjective submersion with symplectic reductions $\mathcal M_\alpha=\mu^{-1}(e_*(\alpha))/G$ as fibers, and admits a compatible homology local system in the sense of \Cref{def: compatible local system}.  Let $c=c(E_0)\in H^2(\mathscr M_0,\mathbb F)\otimes\Lambda^*$ be the curvature class defined in the proof of \Cref{thm: DH after partial reduction}, and write $c_{\alpha_*}=i_{\alpha_*}^*c$ for its restriction under the inclusion $i_{\alpha_*}\from\mathcal M_{\alpha_*}\into\mathscr M_0$.  Then, for every path $\gamma$ from $\alpha_*$ to $\alpha$ in $\Lambda_0$,
    \begin{equation}\label{eq: quasi Hamiltonian variation}
        \mathrm{PT}_\gamma^*[\omega_\alpha]
        =[\omega_{\alpha_*}]-\langle c_{\alpha_*},\alpha-\alpha_*\rangle.
    \end{equation}
\end{theorem}
\begin{proof}
    Since $C_0\subset Z(G)$, we have $e_*^*\theta_{\mathrm L}=e_*^*\theta_{\mathrm R}=\de\alpha$ and $e_*^*\chi=0$.  For the inclusion $j\from E_0\into M$, the first two identities in \eqref{eq: quasi Hamiltonian axioms} therefore give
    \begin{equation}\label{eq: central logarithmic Hamiltonian identities}
        \de(j^*\omega)=0,\qquad
        \iota_{\zeta_\xi}j^*\omega=\de(\tilde\pi,\xi)
        \quad\forall\xi\in\mathfrak g.
    \end{equation}
    The pairing $(\tilde\pi,\xi)$ depends only on the image of $\xi$ in $\mathfrak g^{\mathrm{eff}}$ because $\Lambda$ annihilates the ineffective Lie algebra.  Thus $\varpi=j^*\omega$ satisfies the geometric hypothesis of \Cref{thm: DH after partial reduction} with $B=G^{\mathrm{eff}}$, and its restrictions descend to the reduced forms $\omega_\alpha$.  This proves the theorem.
\end{proof}

\nocite{*}
\bibliographystyle{main}
\bibliography{main}
\end{document}